%% file: WaveIGA.tex
\documentclass{article}
\usepackage[a4paper,top=3cm,bottom=3cm,left=2cm,right=2cm]{geometry}
\usepackage[bookmarks=true, bookmarksopen=true, bookmarksopenlevel=5]{hyperref}
\hypersetup{pdftitle={An unconditionally stable space-time isogeometric method for the acoustic wave equation}}
\usepackage[english]{babel}
\usepackage{amssymb}
\usepackage{amsfonts}
\usepackage{graphicx}
\usepackage{caption}
\usepackage{bm}
\usepackage{bbm} 
\usepackage{amsthm} 
\usepackage{amsmath}
\usepackage{color}
\usepackage{picture}
\usepackage[subrefformat=parens, labelfont=up]{subcaption} 
\usepackage{adjustbox} 
\usepackage{float} 
\usepackage{mathtools}
\usepackage{braket}
\usepackage{nicefrac}
\usepackage{booktabs}
\usepackage{rotating} 
\PassOptionsToPackage{usenames,dvipsnames}{xcolor}
\usepackage[customcolors,norndcorners]{hf-tikz}
\usepackage{tikz}
\usepackage{pgfplots}
\pgfplotsset{compat=newest
}
\usepgfplotslibrary{patchplots}
\usetikzlibrary{pgfplots.patchplots}
\tikzset{set fill color=yellow, set border color=yellow}

\usepackage[backend=bibtex,style=numeric,indexing,maxbibnames=99]{biblatex}  
\numberwithin{equation}{section}

\def\B{\color{black}}

\usepackage{csquotes}
\definecolor{amcol}{rgb}{0.8,0,0}
\definecolor{flcol}{rgb}{0.13, 0.55, 0.13}

\usepackage{soul}
\usepackage{pgfplotstable} %...

\definecolor{revXcol}{rgb}{1,0,0}
\newcommand{\revX}[1]{{\color{revXcol}{#1}}}

\definecolor{revYcol}{rgb}{0,0,1}

\definecolor{revGcol}{rgb}{0,.5,.7}

\theoremstyle{definition}

\theoremstyle{plain}

\theoremstyle{plain}

\theoremstyle{plain}

\theoremstyle{plain}

\theoremstyle{definition}

\theoremstyle{remark}
\newtheorem{remark}{Remark}
\theoremstyle{definition}

\DeclareGraphicsExtensions{{.png},{.pdf},{.jpg}}

\newcommand{\vect}[1]{\boldsymbol{#1}}
\newcommand{\vecx}{\boldsymbol{x}}

\newcommand{\Trial}{\mathcal{W}}
\newcommand{\Test}{\mathcal{V}}

\newcommand{\igastab}{IGA-Stab}
\newcommand{\femstab}{FEM-Stab}

\begin{document}
\date{}
\title{An unconditionally stable space--time isogeometric method for the acoustic wave equation
}

\author{S. Fraschini\thanks{Fakult\"at f\"ur Mathematik, Universit\"at Wien, Oskar-Morgenstern-Platz 1, 1090 Vienna, Austria.} \and G. Loli\thanks{Universit\`a di Pavia, Dipartimento di Matematica ``F. Casorati'', 
		Via A. Ferrata 1, 27100 Pavia, Italy.}  \and A. Moiola\footnotemark[2] \footnotemark[3] \and
	G. Sangalli\footnotemark[2] \thanks{IMATI-CNR ``Enrico Magenes'',  Pavia, Italy. \vskip 1mm \noindent Emails: 
		{sara.fraschini@univie.ac.at, gabriele.loli@unipv.it, 	andrea.moiola@unipv.it, giancarlo.sangalli@unipv.it}}
}

\maketitle
\input{Chapters/Abstract}

%\tableofcontents
\input{Chapters/Introduction}
\input{Chapters/Space_time_variational}
\input{Chapters/Preliminaries_Iga}
\input{Chapters/Space_time_iga}
\input{Chapters/Numerical_results}
\input{Chapters/Conclusions}

\section*{Acknowledgements}
S.~Fraschini was funded by the Austrian Science Fund (FWF) projects \href{https://doi.org/10.55776/F65}{10.55776/F65} and  \href{https://doi.org/10.55776/P33477}{10.55776/P33477}. G.~Loli, A.~Moiola and G.~Sangalli 
are members of the Gruppo Nazionale Calcolo Scientifico - Istituto
Nazionale di Alta Matematica (GNCS-INDAM).
G.~Loli acknowledges support by GNCS-INdAM through the ``Bando Finanziamento Giovani Ricercatori 2021-2022 GNCS''.
A.~Moiola acknowledges support from the PRIN project ``NA-FROM-PDEs'' and ``ASTICE''.
G.~Sangalli acknowledge the support of the Italian Ministry of University and Research (MUR) through the PRIN 2022 PNRR project NOTES (No. P2022NC97R), funded by the European Union - Next Generation EU.
The research has received financial support from ICSC - Italian Research Center on High Performance Computing, Big Data and Quantum Computing, funded by European Union - NextGenerationEU.\\
\begin{figure}[H]
	{\centering
		\hfill\includegraphics[scale=0.075]{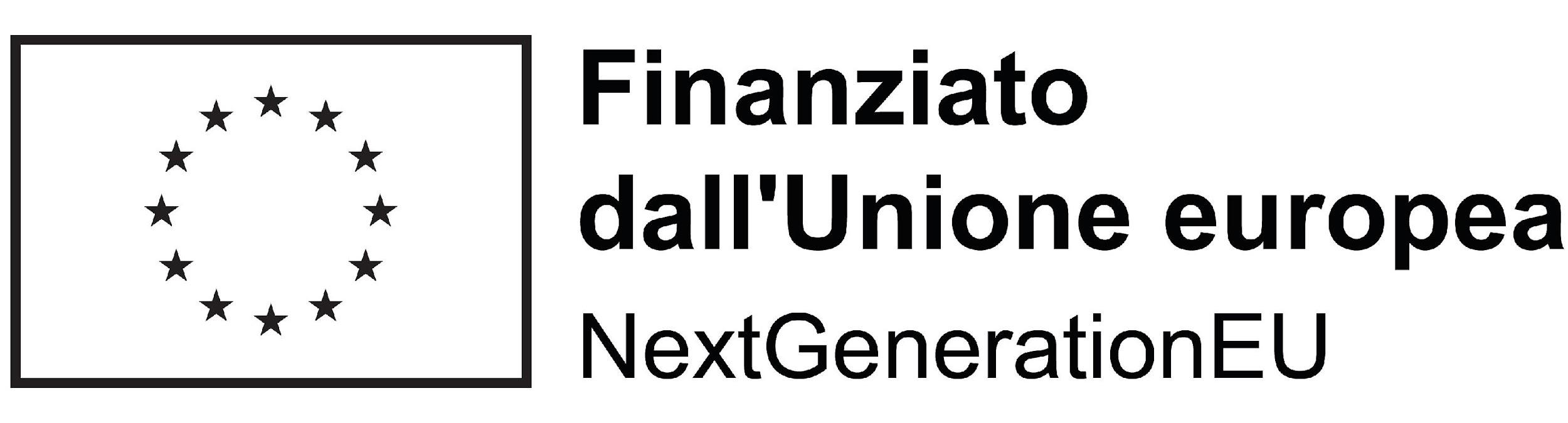}
		\hfill\includegraphics[scale=0.075]{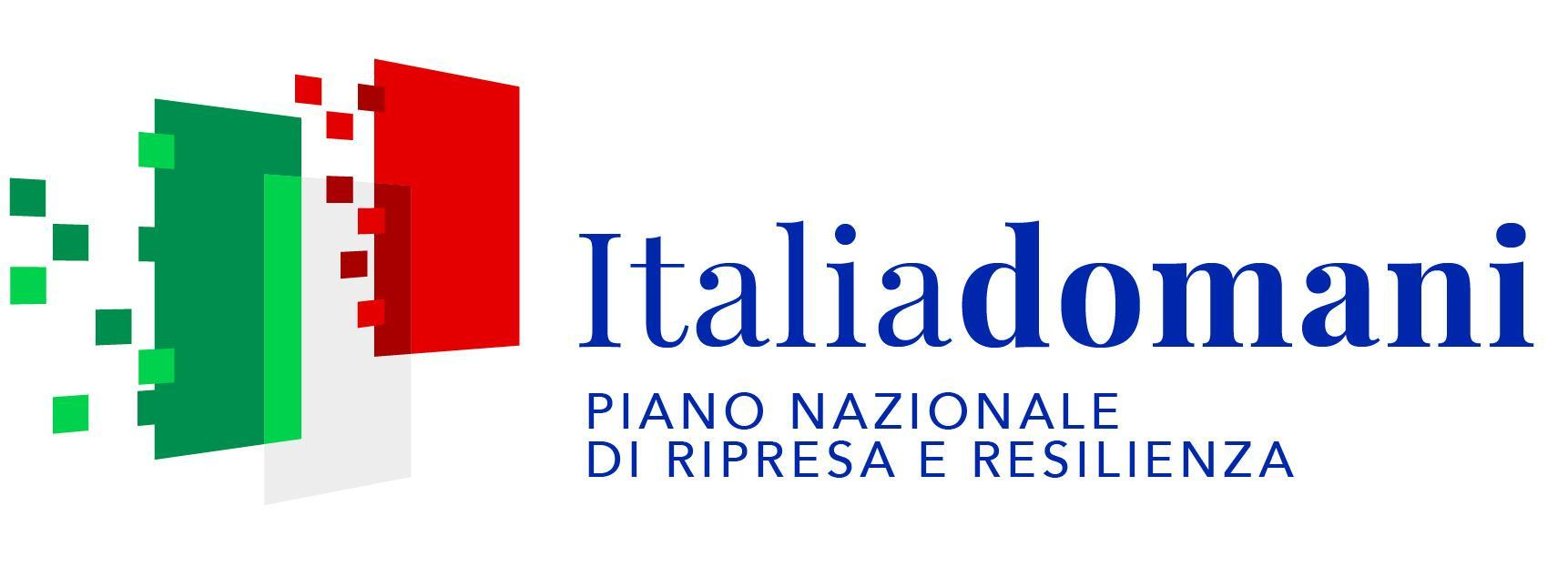}
		\hfill\includegraphics[scale=0.075]{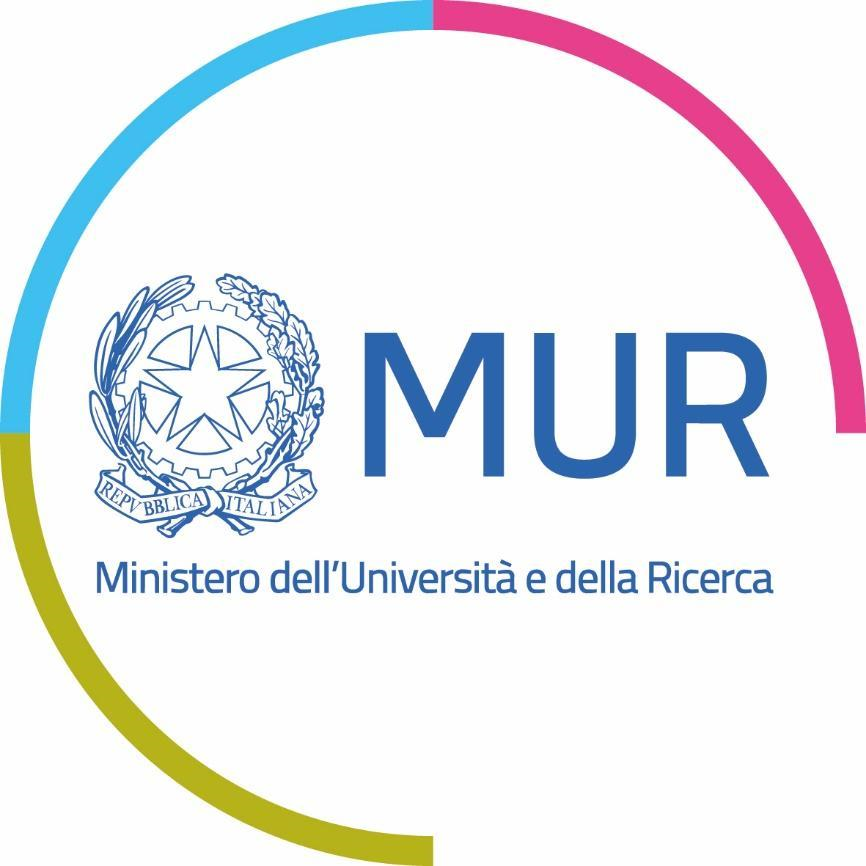}\hfill\mbox{}}
\end{figure}

\printbibliography

\appendix
\input{Chapters/Appendix}
\end{document}

%% file: Chapters/Abstract.tex
\begin{abstract}
%In this paper, we focus on high-order space--time isogeometric discretizations of the linear acoustic wave equation. 
%We deal with smooth approximations in both space and time by employing high-order B-splines of general degree $p$. 
%By exploiting a suitably defined perturbation of order $2p$, we devise a high--order unconditionally stable space--time isogeometric method given by a non-consistent isogeometric formulation.
%To illustrate the effectiveness of this stabilized isogeometric method, we perform several numerical experiments.

We study space--time isogeometric discretizations of the linear acoustic wave equation that use splines of arbitrary degree $p$, both in space and time. 
We propose a space--time variational formulation that is obtained by adding a non-consistent penalty term of order $2p+2$ to the bilinear form coming from integration by parts.
This formulation, when discretized with tensor-product spline spaces
with maximal regularity in time, is unconditionally stable: the mesh size in time is not constrained by the mesh size in space.
We give extensive numerical evidence for the good stability, approximation, dissipation and dispersion properties of the stabilized isogeometric formulation, comparing against stabilized finite element schemes, for a range of wave propagation problems with constant and variable wave speed.

\vskip 1mm
\noindent
\textbf{Keywords:}  Wave equation, isogeometric analysis, space--time Galerkin method, unconditional stability, high-order.
\end{abstract}

%% file: Chapters/Introduction.tex
\section{Introduction}
The principal techniques for the numerical approximation of time-dependent wave phenomena can be categorized into two main classes: separate discretizations of the space and time variables, and the so-called space--time methods. In the former, at first a discretization {either} in space {or} time is used, and then a different discretization technique for the other variable is employed. In the latter, the {simultaneous} discretization of space and time variables in the space--time cylinder is employed. Space--time Finite Element Methods (FEMs) have been  introduced in the seminal papers~\cite{BZ1974,F1969,O1969} and in~\cite{HH1988} for wave propagation problems. Key features of   space--time methods are that they allow: \emph{i)} unstructured meshing of the space--time domain, see, e.g.,~\cite{CKT2009,DKT2007,GSW2017};  \emph{ii)}    efficient treatment of moving boundaries~\cite{TBL1992,TBML1992}; \emph{iii)} multilevel preconditioning~\cite{DFW2016,DFWZ2019}  and  parallelisation~\cite{PSSW2020} in space and time simultaneously.

It is known  (see, e.g.,~\cite{MS2010}) that  high-order methods are crucial to provide accurate numerical approximation of wave-type problems in the presence of high frequencies. There is indeed an ongoing  intense research activity on  high-order discontinuous Galerkin (dG) space-time methods, see for example~\cite{BGL2017,BMPS2021,BCDS2020,DFW2016,DFWZ2019,LSW2014,MP2018,MR2005,PSSW2020}. The dG approach is often considered in this context both for its stability properties, which have been extensively studied for hyperbolic problems, see \cite{shu2009discontinuous} and the references therein, and because it leads, in the temporal discretization, to a splitting of the domain into slabs and then a sequential solution in time, reducing computational costs. In contrast, much less work has been devoted to \textit{high-order} space--time \textit{conforming} FEMs for the linear acoustic wave equation,~\cite{BL1994, BL1995, Z2021,ZL2016}.  In~\cite{BL1994, BL1995}, the proposed position--velocity formulation of the acoustic wave equation does not require any compatibility conditions between the space and time mesh-sizes.
In~\cite{ZL2016}, the authors encompass a broader family of wave phenomena than in~\cite{BL1994, BL1995}, but a restriction on the step-size is required for the well-posedness of the discretization. Subsequent works by Steinbach and Zank \cite{SZ2019,SZ2020,Z2020, Z2021}  have  achieved well-posedness and  unconditional stability using various techniques.
A conforming Petrov--Galerkin discretization of a space--time ultra-weak formulation is presented in \cite{HPSU22}; test functions are splines, while trial functions are coarser in order to ensure discrete inf-sup stability.

In this work, we numerically explore the potential of space--time isogeometric methods for the linear acoustic wave equation.
These methods are based on smooth spline approximations in space and time. In particular, we introduce a novel high-order stabilized formulation with optimal stability and approximation properties.

The isogeometric method, or Isogeometric Analysis (IgA), was introduced in~\cite{HCB2005} as an evolution
of classical FEMs, meant to simplify the interoperability between computer aided design and numerical simulations.
Indeed, IgA employs spline functions, or their generalizations, for both the parametrization of the computational domain and the approximation of the solution of the differential equations of interest. The high-continuity of splines yields  higher accuracy per degree-of-freedom when compared to classical $C^0$ piecewise-polynomial approximation, see, e.g.,~\cite{BS2019,,EBBH2009}, and superior spectral approximation~\cite{HRS2008}.

Indeed, in this work we are interested in the use of spline discretizations in space and time primarily to assess the effect, in wave-propagation problems, of the approximation properties of smooth splines mentioned above.

As for standard continuous FEMs~\cite{SZ2020, Z2021, Z2020}, a plain Galerkin space--time isogeometric discretizations for the acoustic wave equation is stable only under a CFL ({Courant--Friedrichs--Lewy}) condition. Recently, in the finite element case this limitation has been overcome by Steinbach and Zank, who devised in~\cite{SZ2019,Z2020} an unconditionally stable, first-order accurate, space--time FEM, which is extended to high-order continuous FEMs in~\cite{Z2021}. 
However, as shown by our numerical experiments (Figures~\ref{fig:StabZankPConvTot}--\ref{fig:StabZankPRobustTot} of Appendix~\ref{sec:appendix}), the same stabilization technique does not provide unconditional stability of a space--time isogeometric discretization  based on high-order smooth splines. Therefore, taking inspiration by \cite[Lemma 4.2.26]{Z2020}, we introduce a new stabilization that extends~\cite{SZ2019,Z2020} to high-order smooth space--time IgA, and that enjoys stability and optimal convergence (as well as desired energy properties) without any constraint on the time mesh-size.
Differently from \cite{SZ2019,Z2020}, the stabilization we propose in \eqref{equivstabiga} is implemented as a simple penalty term involving high-order derivatives.
For a range of one- and two-dimensional problems with constant and variable wave speed, we study numerically the unconditional stability, the convergence rates, the approximation of high-frequency waves, the dissipation and dispersion properties: under all these respect we observe a computational advantage in terms of accuracy per degree-of-freedom when comparing against high-order FEMs.

The outline of the paper is as follows. We describe the variational setting for the initial boundary value problem under consideration in Section \ref{sec:IBVP} and then present the basics of the space--time B-splines based IgA in Section \ref{sec:iga}.
In Section \ref{sec:st_iga_wave}, we introduce the novel high-order stabilization, and in Section \ref{sec:num_exp}, we present the numerical results illustrating the stability and convergence properties of the proposed method. Finally, in Section \ref{sec:conclusions}, we draw conclusions and highlight some future research directions.
%%% Local Variables:
%%% mode: latex
%%% TeX-master: t
%%% End:

%% file: Chapters/Space_time_variational.tex
\section{Model problem}\label{sec:IBVP}
Let $\Omega \subset \mathbb{R}^d$ be a bounded, open, Lipschitz domain and %let us 
assume that $\partial \Omega$ is partitioned as ${\partial \Omega=\overline{\Gamma_D \cup \Gamma_N \cup \Gamma_R}}$. 
Given a final time $T>0$ and a positive piecewise continuous wave velocity $c$ on $\overline\Omega\times [0,T]$, let us consider the following differential problem:
\begin{equation}\label{iga_wave_strong_eq}
	\begin{cases}
		\begin{aligned}
			\partial_{tt} u - \nabla \cdot (c^2 \nabla u) &= f &&\text{on }Q:=\Omega \times (0,T), \\
			u &= g_D &&\text{on } \Sigma_D:=\Gamma_D \times (0,T), \\
			c^2 \nabla u \cdot \mathbf{n} &= g_N &&\text{on } \Sigma_N:=\Gamma_N \times (0,T), \\
			{\vartheta c \partial_t u + c^2 \nabla u \cdot \mathbf{n}} &= g_R &&\text{on } \Sigma_R:=\Gamma_R \times (0,T),\\
			u &= u_0 &&\text{on } \Sigma_0:=\Omega \times \{0\}, \\
			\partial_t u &= u_1 &&\text{on } \Sigma_0,
		\end{aligned}	
	\end{cases}
\end{equation}
where $\mathbf{n}$ denotes the outward-pointing unit vector normal on $\partial\Omega$, and the impedance parameter $\vartheta>0$ represents how easily the obstacle bounded by $\Gamma_R$ yields to the acoustic wave. 

Let us start with the case of homogeneous initial and Dirichlet boundary condition (i.e. $u_0=g_D=0$). 
To derive a space--time variational formulation of~\eqref{iga_wave_strong_eq}, we introduce the following space--time Hilbert spaces
\begin{equation*}
	\begin{aligned}
		\Trial&:=L^2\big(0,T;H^1_{0,\Gamma_D}(\Omega)\big)\cap H^1_{0,*}\big(0,T;L^2(\Omega)\big), \\
		\Test&:=L^2\big(0,T;H^1_{0,\Gamma_D}(\Omega)\big)\cap H^1_{*,0}\big(0,T;L^2(\Omega)\big),
	\end{aligned}
\end{equation*}
where we have defined
\begin{equation*}
\begin{aligned}
	H^1_{0,\Gamma_D}(\Omega) &:= \left\{ v \in H^1(\Omega) : \ v|_{\Gamma_D}=0 \right\}, \\
	H^1_{0,*}\big(0,T;L^2(\Omega)\big)&:=\{v \in H^1\big(0,T;L^2(\Omega)\big): \ v(0)=0 \ \text{in} \ L^2(\Omega)\},\\
	H^1_{*,0}\big(0,T;L^2(\Omega)\big)&:=\{v \in H^1\big(0,T;L^2(\Omega)\big): \ v(T)=0 \ \text{in} \ L^2(\Omega)\}.
\end{aligned}
\end{equation*}
Following~\cite{Z2020}, we endow both $\Trial$ and $\Test$ with the weighted Sobolev seminorm $|\cdot|_{c,H^1(Q)}$ defined as 
\begin{equation}\label{H1_sem}
	|v|^2_{c,H^1(Q)}:=\int_{Q} \Big(|\partial_t v(x,t)|^2+ c(x,t)^2 |\nabla v(x,t)|^2 \Big) \ dQ \quad \text{for } v \in H^1(Q).
\end{equation}
When the velocity is $c=1$ in $Q$, %and with velocity $c\equiv 1$ 
we denote \eqref{H1_sem} with $|\cdot|_{H^1(Q)}$.
Multiplying by a test function $v \in \Test$ and integrating by parts the first equation in~\eqref{iga_wave_strong_eq}, we get the following space--time variational formulation:
\begin{equation}\label{abstractvar}
	\text{Find } u \in \Trial \text{ such that } a(u,v)=\mathcal{F}(v) \quad \text{for all }v \in \Test,
\end{equation}
where the bilinear form $a:  \Trial \times \Test \rightarrow \mathbb{R}$ and the linear form $\mathcal{F}: \Test \rightarrow \mathbb{R}$ are defined as
\begin{align*}
	a(w,v) &:= \int_{Q}  \big(c^2(x,t) \nabla w(x,t) \cdot \nabla v(x,t) - \partial_t w(x,t) \partial_t v(x,t) \big)\ dQ + \int_{\Sigma_R} \vartheta c(x,t) \partial_t w(x,t)v(x,t) \ d\Sigma_R\\
	\mathcal{F}(v) &:= \int_{Q} f(x,t)v(x,t)\ dQ + \int_{\Sigma_0} u_1(x) v(x,0) \ d\Sigma_0 + \int_{\Sigma_N} g_N(x,t) v(x,t) \ d\Sigma_N + \int_{\Sigma_R} g_R(x,t) v(x,t) \ d\Sigma_R,
\end{align*}
for $w \in \Trial$ and $v \in \Test$.

\begin{remark}
To deal with
inhomogeneous initial and Dirichlet boundary conditions,
we let $\bar u$ be a lifting of $u_0$ and $g_D$, i.e.,~$\bar u \in H^1(Q)$ such that $\bar{u}|_{\Sigma_0}=u_0$ and $\bar{u}|_{\Sigma_D}=g_D$.
Then, the space--time variational formulation becomes:
	\begin{equation*}
		\text{Find } u = \bar{u} + w \text{ with } w \in \Trial \text{ such that } a(w,v)=\mathcal{F}(v)-a(\bar{u},v) \quad \text{for all } v \in \Test.
	\end{equation*}
\end{remark}
The well-posedness of the variational problem~\eqref{abstractvar} is well known under suitable conditions such as $f \in L^2(Q)$ and $\Gamma_N=\Gamma_R=\emptyset$
(see~\cite{L1985,SZ2020,Z2020,Zl1994}).
On the other hand,
the abstract variational formulation~\eqref{abstractvar} with general
right-hand side $\mathcal F\in \mathcal V^*$ does not fit into the setting
of the Banach--Ne\v cas--Babu\v ska Theorem (see~\cite{N1962}) since it does not satisfy the inf-sup condition w.r.t.\ the
norm~\eqref{H1_sem} (see~\cite[Theorem 4.2.24]{Z2020} and~\cite{SZ2022}).

%% file: Chapters/Preliminaries_Iga.tex
\section{Preliminaries on isogeometric analysis}
\label{sec:iga}
\subsection{B-splines}
\label{sec:bsplines}
Given two non-negative integers $p$ and $m$, with $m \geq p+1$, a univariate B-spline basis can be generated starting from an open knot vector in $[0,1]$, which is a sequence of non-decreasing points, also called knots, i.e., 
\begin{equation*}
	\Xi:=\left\{ 0=\xi_1=\dots=\xi_{p+1} \leq \dots \leq \xi_{m}=\dots=\xi_{m+p+1}=1\right\}.
\end{equation*}
For future reference, we also introduce the vector $Z:=\{ \zeta_1,\ldots,\zeta_N\}$ of knots without repetitions, also called breakpoints.
Then, according to Cox--De Boor recursion formulas (see~\cite{DeBoor2001}), univariate B-splines are piecewise polynomials defined for $i=1,\dots,m$ as
\begin{align*} 
	\widehat{b}_{i,0}(\eta) =
	\begin{cases}
	1 &  { \text{if }} \xi_{i}\leq \eta<\xi_{i+1},\\
	0 & \text{otherwise},
\end{cases}
\end{align*}
for $p=0$, while for $p>0$
\begin{align*}
	\widehat{b}_{i,p}(\eta)= \begin{cases}\dfrac{\eta-\xi_{i}}{\xi_{i+p}-\xi_{i}}\widehat{b}_{  i,p-1}(\eta)   +\dfrac{\xi_{i+p+1}-\eta}{\xi_{i+p+1}-\xi_{i+1}}\widehat{b}_{  i+1,p-1}(\eta)   & { \text{if }} \xi_{i}\leq  \eta<\xi_{i+p+1}, \\
	0  & \text{otherwise},
	\end{cases}
\end{align*}
where we adopt the convention $0/0=0$. The univariate spline space is defined as
\begin{equation*}
	\widehat{\mathcal{S}}_h^p : = \mathrm{span}\{\widehat{b}_{i,p}\}_{i = 1}^m,
\end{equation*}
where $h$ denotes the mesh-size, i.e. $h:=\max\{ |\xi_{i+1}-\xi_i| : \ i=1,\dots,m+p \}$. We recall that the smoothness of the obtained B-spline basis is $C^{p-k}$ at every knot, where $k$ denotes the multiplicity of the considered knot, while it is $C^{\infty}$ elsewhere.
For more details on B-splines properties  and their use in IgA we refer to~\cite{Cottrell2009,DeBoor2001,beirao}.

Multivariate B-splines are defined as tensor-product of univariate B-splines. To deal with space--time variational formulation, we consider splines that depend on $d$ space variables and the time variable. Therefore, given non-negative integers $m_\ell, p_\ell$, with $m_\ell \geq p_\ell+1$ for $\ell=1,\dots,d$ and $m_t,p_t$, with $m_t \geq p_t+1$, we  introduce $d+1$ univariate open knot vectors 
\begin{align*}
	\Xi_\ell &:=\left\{ 0=\xi_{\ell,1}=\dots=\xi_{\ell,p_\ell+1} \leq \dots \leq \xi_{\ell,m_\ell}=\dots=\xi_{\ell,m_\ell+p_\ell+1}=1\right\} & & \text{for } \ell=1,\ldots, d, \\
	\Xi_t &:=\left\{ 0=\xi_{t,1}=\dots=\xi_{t,p_t+1} \leq \dots \leq \xi_{t,m_t}=\dots=\xi_{t,m_t+p_t+1}=1\right\}
\end{align*}
and $d+1$ univariate breakpoints vectors
\begin{align*}
	Z_\ell &:=\left\{ \zeta_{\ell,1},\ldots,\zeta_{\ell,N_\ell} \right\} & & \text{for } l=1,\ldots, d, \\
 	Z_t &:=\left\{ \zeta_{t,1},\ldots,\zeta_{t,N_t} \right\}.
\end{align*}
Let $h_\ell$ be the mesh-size associated to the knot vector $\Xi_\ell$ for $\ell=1,\dots,d$, let $h_s:=\max\{h_\ell\ | \ \ell=1,\dots, d\}$ be the maximal mesh-size in all space knot vectors and let $h_t$ be the mesh-size of the time knot vector $\Xi_t$. 
Let also $\boldsymbol{p}$ be the vector that contains the degree indexes, i.e.,~$\boldsymbol{p} :=(\boldsymbol{p}_s,p_t)$, where $\boldsymbol{p}_s:=(p_1,\dots,p_d )$. For simplicity, we assume to have the same polynomial degree in all space directions, i.e., with abuse of notations, we set  $p_s:=p_1=\dots=p_d$.
Then, multivariate B-splines are defined as
\begin{equation*} 
	\widehat{B}_{ \vect{i},\vect{p}}(\vect{\eta},\tau) : =
	\widehat{B}_{\vect{i_s}, \vect{p}_s}(\vect{\eta}) \widehat{b}_{i_t,p_t}(\tau),
\end{equation*}
where 
\begin{equation*}\label{eq:tens-prod-space-parametric}
		\widehat{B}_{\vect{i_s},\vect{p}_s}(\vect{\eta}):=\widehat{b}_{i_1,p_s}(\eta_1) \revX{\cdots} \widehat{b}_{i_d,p_s}(\eta_d),
\end{equation*}
$\vect{i_s}:=(i_1,\dots,i_d)$, $\vect{i}:=(\vect{i_s}, i_t)$ and $\vect{\eta} := (\eta_1, \ldots, \eta_d)$.
The corresponding spline space is defined as
\begin{equation*}
	\widehat{\mathcal{S}} ^{\vect{p}}_{h}:= \mathrm{span}\left\{\widehat{B}_{\vect{i}, \vect{p}} : \ i_\ell = 1,\dots, m_\ell \text{ for } \ell=1,\dots,d; \ i_t=1,\dots,m_t \right\},
\end{equation*} 
where $h:=\max\{h_s,h_t\}$. 
We have that $\widehat{\mathcal{S}} ^{\vect{p}}_{ {h}}  =\widehat{\mathcal{S}} ^{ \vect{p}_s}_{ {h}_s} \otimes
\widehat{\mathcal{S}} ^{p_t}_{h_t}$, where
\begin{equation*}
	\widehat{\mathcal{S}} ^{\vect{p}_s}_{h_s} := \mathrm{span}\left\{\widehat{B}_{\vect{i_s},\vect{p}_s} :  \ i_\ell = 1,\dots, m_\ell; \ \ell=1,\dots,d  \right\}
\end{equation*}
is the space of tensor-product splines on   $\widehat{\Omega}:=(0,1)^d$.
We assume that $p_t, p_s\geq 1$ and that $\widehat{\mathcal{S}} ^{\vect{p}_s}_{h_s}
\subset C^{p_s-1}(\widehat{\Omega}  )$ and  $\widehat{\mathcal{S}} ^{{p}_t}_{h_t}
\subset C^{p_t-1}\left((0,1)\right)$. 
If the space and time meshes are uniform, we have
$\mathrm{dim}(\widehat{\mathcal{S}} ^{\vect{p}}_h)=(p_s+h_s^{-1})^d(p_t+h_t^{-1})$.

\subsection{Space--time isogeometric spaces}
\label{sec:iso_space}
Let us assume that our space computational domain $\Omega$ is parametrized by $\vect{F}: \widehat{\Omega} \rightarrow {\Omega}$, with ${\vect{F}\in {\left[\widehat{\mathcal{S}}^{\vect{p}_s}_{{h}_s}\right]^d}}$. Moreover, we assume that $\vect{F}^{-1}$ has piecewise bounded derivatives of any order.
The space--time computational domain that we consider is then $\Omega\times (0,T)$, where $\Omega\subset\mathbb{R}^d$ and $T>0$ is the final time. 
We define  ${\vecx:=(x_1,\dots,x_d)=  \vect{F}(\vect{\eta})}$ and
$t:=T\tau$. Then, the space--time domain  is given by  the parametrization
$\vect{G}:\widehat{\Omega}\times(0,1)\rightarrow
\Omega\times(0,T)$, such that $ \vect{G}(\vect{\eta}, \tau):=(\vect{F}(\vect{\eta}), T\tau )=(\vecx,t).$

We introduce the spline space with initial, final and Dirichlet boundary conditions, in parametric coordinates, as
\begin{equation}\label{eq:all_basis}
\begin{aligned}
	\widehat{\mathcal{W}}_{h}&:=\left\{ \widehat{w}_h\in \widehat{\mathcal{S}}^{\vect{p}}_h : \ \widehat{w}_h = 0 \text{ on } \widehat{\Gamma}_D\times (0,1) \text{ and } \widehat{w}_h = 0 \text{ on } \widehat{\Omega}\times\{0\} \right\},\\
	\widehat{\mathcal{V}}_{h}&:=\left\{ \widehat{v}_h\in \widehat{\mathcal{S}}^{\vect{p}}_h : \ \widehat{v}_h = 0 \text{ on } \widehat{\Gamma}_D\times (0,1) \text{ and } \widehat{v}_h = 0 \text{ on } \widehat{\Omega}\times\{1\} \right\},
\end{aligned}
\end{equation}
where $\widehat{\Gamma}_D:=\vect{F}^{-1}(\Gamma_D)$ and, for simplicity, we assume $\widehat{\Gamma}_D \subset \partial \widehat{\Omega}$ is the union of entire faces of $\widehat{\Omega}$. We observe that 
\begin{equation*}
\begin{aligned}
	\widehat{\mathcal{W}}_{h}= \widehat{\mathcal{X}}_{s,h_s}\otimes  \widehat{\mathcal{W}}_{t,h_t} && \text{and} && \widehat{\mathcal{V}}_{h}= \widehat{\mathcal{X}}_{s,h_s}\otimes  \widehat{\mathcal{V}}_{t,h_t},
\end{aligned}
\end{equation*}
where 
\begin{equation*}
\begin{aligned}
	\widehat{\mathcal{X}}_{s,h_s}   & := \left\{ \widehat{B}_{h_s}\in \widehat{\mathcal{S}}^{\vect{p}_s}_{h_s}   : \ \widehat{B}_{h_s} = 0 \text{ on } \widehat{\Gamma}_D \right\}
	,\\ 
	\widehat{\mathcal{W}}_{t,h_t} & := \left\{ \widehat{w}_{h_t}\in \widehat{\mathcal{S}}^{ p_t}_{h_t} :  \ \widehat{w}_{h_s}( 0)=0 \right\},\\
	\widehat{\mathcal{V}}_{t,h_t} & := \left\{ \widehat{v}_{h_t}\in \widehat{\mathcal{S}}^{ p_t}_{h_t} :  \ \widehat{v}_{h_t}( 1)=0 \right\}.
	\end{aligned}
\end{equation*}
We highlight that $\widehat{\mathcal{W}}_{h}$ and $\widehat{\mathcal{V}}_{h}$ have the same dimension, that is
\begin{equation*}
	\dim (\widehat{\mathcal{W}}_{h}) = \dim (\widehat{\mathcal{V}}_{h})= N_{\mathrm{dof}},
\end{equation*}
where we defined $N_{\mathrm{dof}}:=N_s n_t$, with $N_s:=\prod_{\ell=1}^dn_{s,\ell}$ and $n_t:=m_t-1$.

Finally, the isogeometric spaces we consider are the isoparametric push-forward of \eqref{eq:all_basis} through the geometric map $\vect{G}$, i.e.,
\begin{equation}\label{eq:iga_spaces}
\begin{aligned}
	\mathcal{W}_{h} &:= \left\{  w_h:=\widehat{w}_h\circ \vect{G}^{-1} \ \middle| \ \widehat{w}_h \in \widehat{\mathcal{W}}_{h} \right\} \subset \Trial, \\
	\mathcal{V}_{h} &:= \left\{  v_h:=\widehat{v}_h\circ \vect{G}^{-1} \ \middle| \ \widehat{v}_h \in \widehat{\mathcal{V}}_{h} \right\} \subset \Test.
\end{aligned}	
\end{equation}
We also have that
\begin{equation*}
\begin{aligned}
	\mathcal{W}_{h}=\mathcal{X}_{s,h_s}\otimes \mathcal{W}_{t,h_t} && \text{and} && \mathcal{V}_{h}= \mathcal{X}_{s,h_s}\otimes  \mathcal{V}_{t,h_t},
\end{aligned}	
\end{equation*}
where
\begin{equation*}
	\begin{aligned}
		\mathcal{X}_{s,h_s}&:=\left\{B_{h_s}:= \widehat{B}_{h_s}\circ \vect{F}^{-1} \ \middle| \ \widehat{B}_{h_s} \in \widehat{\mathcal{X}}_{s,h_s} \right\},\\
		\mathcal{W}_{t,h_t}&:=\left\{  {w}_{{h}_t}:= \widehat{w}_{{h}_t}( \cdot /T) \ \middle| \ \widehat{w}_{{h}_t} \in \widehat{\mathcal{W}}_{t,h_t} \right\},\\
		\mathcal{V}_{t,h_t}&:=\left\{  {v}_{{h_t}}:= \widehat{v}_{{h}_t}( \cdot /T) \ \middle| \ \widehat{v}_{{h}_t} \in \widehat{\mathcal{V}}_{t,h_t} \right\}.
	\end{aligned}	
\end{equation*}

%% file: Chapters/Space_time_iga.tex
\section{Space--time isogeometric method}\label{sec:st_iga_wave}
Recalling the definition of the isogeometric trial and test spaces in~\eqref{eq:iga_spaces}, the Petrov--Galerkin isogeometric discretization of the variational formulation \eqref{abstractvar} reads as
\begin{equation}\label{igawave}
	\text{Find } u_h \in \mathcal{W}_h \text{ such that } a(u_h,v_h)=\mathcal{F}(v_h) \quad \text{for all } v_h \in \mathcal{V}_h.
\end{equation}

As stated in~\cite{SZ2020, Z2020} for space--time finite element discretizations of \eqref{abstractvar}, a CFL condition needs to be satisfied to guarantee the stability of the numerical scheme. An immediate extension of the argument of~\cite{SZ2020, Z2020} shows that the same result holds true for the space--time isogeometric method \eqref{igawave}. 
Specifically, for a uniform discretization with time mesh-size $h_t$ and space mesh-size $h_s$, which is assumed to be uniform in all space directions of the parametric domain, a CFL condition
\begin{equation*}
	h_t \leq C h_s
\end{equation*} 
is required for stability of the space--time isogeometric method \eqref{igawave}. The stability estimate above follows from the results shown in~\cite{F2021t} concerning the splines discretization of a suitable ordinary differential, and the inverse inequalities for splines.

In order to guarantee unconditional stability, i.e.,~the absence of a CFL condition, taking inspiration from~\cite[Lemma 4.2.26]{Z2020} (equivalently,~\cite[Lemma~17.7]{SZ2019}), we propose the following stabilized space--time isogeometric formulation:

\medskip

\noindent\fbox{\parbox{0.99\textwidth}{
\begin{equation}\label{stabfemwave_iga}
	\text{Find } u_h \in \mathcal{W}_h \text{ such that } a_{\mathrm{IGA}}(u_h,v_h)=\mathcal{F}(v_h) \quad \text{for all } v_h \in \mathcal{V}_h,
\end{equation}
where the bilinear form $a_{\mathrm{IGA}}(\cdot,\cdot):  \mathcal{W}_h \times \mathcal{V}_h \rightarrow \mathbb{R}$ is defined as
\begin{equation}\label{equivstabiga}
a_{\mathrm{IGA}}(w_h,v_h) := a(w_h,v_h)- \delta \sum_{k=1}^{N_t-1} (T\zeta_{t,k+1}-T\zeta_{t,k})^{2p_t}\int_{\Omega \times (T\zeta_{t,k},T\zeta_{t,k+1})} c^2(x,t) \ \partial^{p_t}_t\nabla w_h(x,t) \cdot \partial^{p_t}_t\nabla v_h(x,t) \ d\Omega \hspace{0.05cm} dt 
\end{equation}
for $w_h \in \mathcal{W}_h$, $v_h \in \mathcal{V}_h$.}}

\medskip

The choice of the stabilization parameter %constant 
$\delta>0$ will be investigated in Section~\ref{s:ExSmooth}. 
Note that formulation \eqref{equivstabiga} is non-consistent: $a_{\mathrm{IGA}}(u,v_h)\ne\mathcal F(v_h)$ for $u$ solution of \eqref{iga_wave_strong_eq} and a general $v_h\in \mathcal V_h$.

\begin{remark}
By considering globally continuous, tensor-product space--time piecewise-polynomial functions of arbitrary degree $p$, the unconditionally stable space--time finite element method of \cite{Z2021} for $c =1$ is based on the following discrete bilinear form
	\begin{equation}\label{stabZank}
		a_{\mathrm{FEM}}(w_h,v_h) := \int_{Q} \Big[- \partial_t w_h(x,t) \partial_t v_h(x,t) + \nabla w_h(x,t) \cdot \big(\mathcal Q_{h_t}^{p-1,\text{disc}}\nabla v_h\big)(x,t)\Big] \ dQ,
	\end{equation}
	where the operator
	\begin{equation*}
		\mathcal Q_{h_t}^{p-1,\text{disc}}: L^2(Q) \longrightarrow L^2(\Omega) \otimes \mathcal{S}_{h_t}^{p-1,\text{disc}}(0,T)
	\end{equation*}
	is the $L^2$-orthogonal projection onto the space of piecewise-polynomial, discontinuous functions with respect to the time variable, which acts componentwise in \eqref{stabZank}.
	As a consequence of~\cite[Lemma 4.2.26]{Z2020}, if $c = 1$, $\Gamma_N=\Gamma_R=\emptyset$,  $g_D=u_0=u_1 =0$, $p_s=p_t=1$, 
	the two bilinear forms coincide: $a_{\mathrm{FEM}}(w_h,v_h)=a_{\mathrm{IGA}}(w_h,v_h)$ 
	with $\delta = \frac{1}{12}$.
\end{remark}

In the following, we name ``\igastab'' the stabilized formulation \eqref{stabfemwave_iga}, and ``\femstab'' the same formulation with $a_{\mathrm{FEM}}$ in place of $a_{\mathrm{IGA}}$.

%%% Local Variables:
%%% mode: latex
%%% TeX-master: "../WaveIGA"
%%% End:

%% file: Chapters/Numerical_results.tex
\section{Numerical experiments}\label{sec:num_exp}
In this section, we provide numerical experiments illustrating the convergence properties of the \igastab\ method \eqref{equivstabiga}.
We show three sets of experiments.
In the first, we show the convergence rates for a range of problems (constant, smooth, and discontinuous coefficients, one and two space dimensions, Dirichlet, Neumann, and mixed boundary conditions);
in the second, we test how well the \igastab\ method~\eqref{equivstabiga} preserves the solution energy; and, in the third, we study the dispersion properties of \igastab\ method. 
From now on, $p$ denotes the spline degree in both space and time.
All the tests are performed with Matlab R2022b and GeoPDEs toolbox \cite{defalco}, and the corresponding codes are available in the GitHub repository~\cite{XTIgAWaves}.
The Matlab direct solver is used for all the numerical tests except for
Example~\ref{ex:scattering}, where an iterative solver is employed.

\subsection{Unconditional stability and accuracy}
\subsubsection{Smooth solution: stability and accuracy of IGA-Stab and FEM-Stab}
\label{s:ExSmooth}
As in \cite[p.~367]{SZ2019}, we consider the following exact solution of the acoustic wave equation \eqref{iga_wave_strong_eq} in one dimension:
\begin{equation}\label{u1_ex} 
u(x,t)= \sin(\pi x)\sin^2\left(\tfrac{5}{4} \pi t\right) \quad \text{for } (x,t) \in \Omega \times (0,T):=(0,1) \times (0,10),
\end{equation}
with velocity $c= 1$, $\Gamma_N = \Gamma_R = \emptyset$, $g_D = u_0 = u_1 = 0$, and $f$ appropriately chosen.

We demonstrate the effectiveness of \igastab\ method~\ref{equivstabiga} with $\delta=10^{-p}$.
This choice is motivated by the results in Figure~\ref{fig:ErrorStabDelta}, where for $\delta<10^{-p}$ we observe that instability kicks in and all relative errors quickly grow, whereas for $\delta>10^{-p}$ the errors slowly grow because of the larger inconsistent term in \eqref{equivstabiga}.

Figure~\ref{fig:StabRobust} demonstrates that the \igastab\ method is stable and no CFL condition is needed.

In Table~\ref{table:stab_comparison} we summarize the stability properties obtained with different or no stabilization (see Figure~\ref{fig:NoStabRobustTot} of Appendix~\ref{sec:appendix}) and with splines of different regularity, as observed in numerical experiments.
The \igastab\ formulation is unconditionally stable also if we decrease the spline regularity in space, but requires a CFL condition if we reduce the spline regularity in time, see Figures~\ref{fig:StabConvTot} and~\ref{fig:StabRobustTot} of Appendix~\ref{sec:appendix}.
Furthermore, there is numerical evidence that \femstab~\eqref{stabZank} provides unconditional stability even if applied to continuous piecewise-polynomial functions in time and splines of arbitrary regularity in space, but a CFL condition kicks in as soon as the basis functions are at least $C^1$ in time, see Figures~\ref{fig:StabZankPConvTot} and~\ref{fig:StabZankPRobustTot} of Appendix~\ref{sec:appendix}.  
This suggests that, to obtain an unconditionally stable formulation in the ``intermediate smoothness'' case, i.e.\ with $C^{q_t}$-continuous splines in time of degree $p_t$ and $1\le q_t\le p_t-2$, one needs a combination of a volume stabilization similar to \eqref{equivstabiga} and a jump penalization term at the times $\zeta_{t,k}$; this is the subject of current investigation.
In the rest of this paper, we only consider \igastab\ with maximal-degrees splines in both space and time, and \femstab\ with $C^0$ finite elements in both space and time.

In Figure~\ref{fig:StabComparisonHP} we compare the errors committed by
the stable methods, i.e.\ \igastab\ \eqref{equivstabiga} for splines
of maximal regularity, and \femstab\ \eqref{stabZank} for $C^0$
splines. Both methods converge with optimal rates, with  \igastab\
showing better accuracy per degree-of-freedom, with greater advantage  for larger values
of $p$. 

\begin{table}%[h]                                                     
	\centering                                         
	\begin{tabular}{cccc}   
		\toprule                                           
		Stabilization & Spline--regularity in time & Spline--regularity in space & Stability \\
		\midrule
		None & $\forall$ & $\forall$ & CFL \\	
		
		\femstab~\eqref{stabZank} & $C^{0}$ & $\forall$  & Unconditional \\
		
		\femstab~\eqref{stabZank} & $C^{q_t}, q_t>0$ & $\forall$ &  CFL \\
		
		\igastab~\eqref{equivstabiga}  & $C^{q_t}$, $q_t<p_t-1$ & $\forall$ & CFL \\		             
		
		\igastab~\eqref{equivstabiga} & $C^{p_t-1}$ & $\forall$ &  Unconditional \\		
		\bottomrule
	\end{tabular} 
	\caption{Comparison of the stability of different stabilized formulations and different spline regularities. Here $p_t$~denotes the spline degree in time.
	\femstab\ and \igastab\ are unconditionally stable if and only if the spline regularity in time is the minimal ($C^0(0,T)$) and maximal ($C^{p_t-1}(0,T)$) possible, respectively.} 
	\label{table:stab_comparison}           
\end{table}   
\input{Graphics/ErrorStabDeltaFine}
\input{Graphics/ErrorStabRobust}
\input{Graphics/ErrorZankOurCompHP}

\subsubsection{High-frequency oscillations}

\input{Graphics/GraficoErroreLdof_p_1234_comparison}

We study the robustness of the method~\eqref{equivstabiga} with respect to the frequencies of oscillatory solutions.
We approximate~\eqref{iga_wave_strong_eq} on a unidimensional space domain $\Omega:=(0,1)$, with exact solution
\begin{equation}\label{u2_ex}
%u_{\mathrm{ex}}
u (x,t)= \sin(k \pi x) \sin(k \pi t) \quad \text{for } (x,t) \in \Omega \times (0,T):=(0,1) \times (0,2),
\end{equation}
for different values of $k\in\mathbb N$, $c=1$, $\Gamma_N = \Gamma_R = \emptyset$, $g_D=u_0=f=0$ and $u_1(x)=k\pi\sin(k\pi x)$.
As before, we choose $\delta = 10^{-p}$, for $p=1,2,3,4$, and we set $h_t=h_s$. 
Having defined the number of space wavelengths in $\Omega$ as $\sharp \lambda := k/2$, Figure~\ref{fig:Ldof1234} shows the comparison between the
relative errors of the \igastab\ and \femstab\ methods in the space--time $L^2$-norm and $H^1$-seminorm \eqref{H1_sem}, plotted against $N_{s} /\sharp \lambda$, where $N_{s}$ denotes the number of degrees of freedom in space. 
We observe that, for $p>1$, the number of degrees of freedom per
wavelength needed to obtain a given accuracy is roughly independent of
$k$, and that \igastab\ reaches high-accuracy with less
degrees-of-freedom per wavelength, compared to    \femstab\ .

\subsubsection{Scattering problem}
\label{ex:scattering}
\begin{figure}[htb]
	\centering %\hfill
	\includegraphics[width=\linewidth]{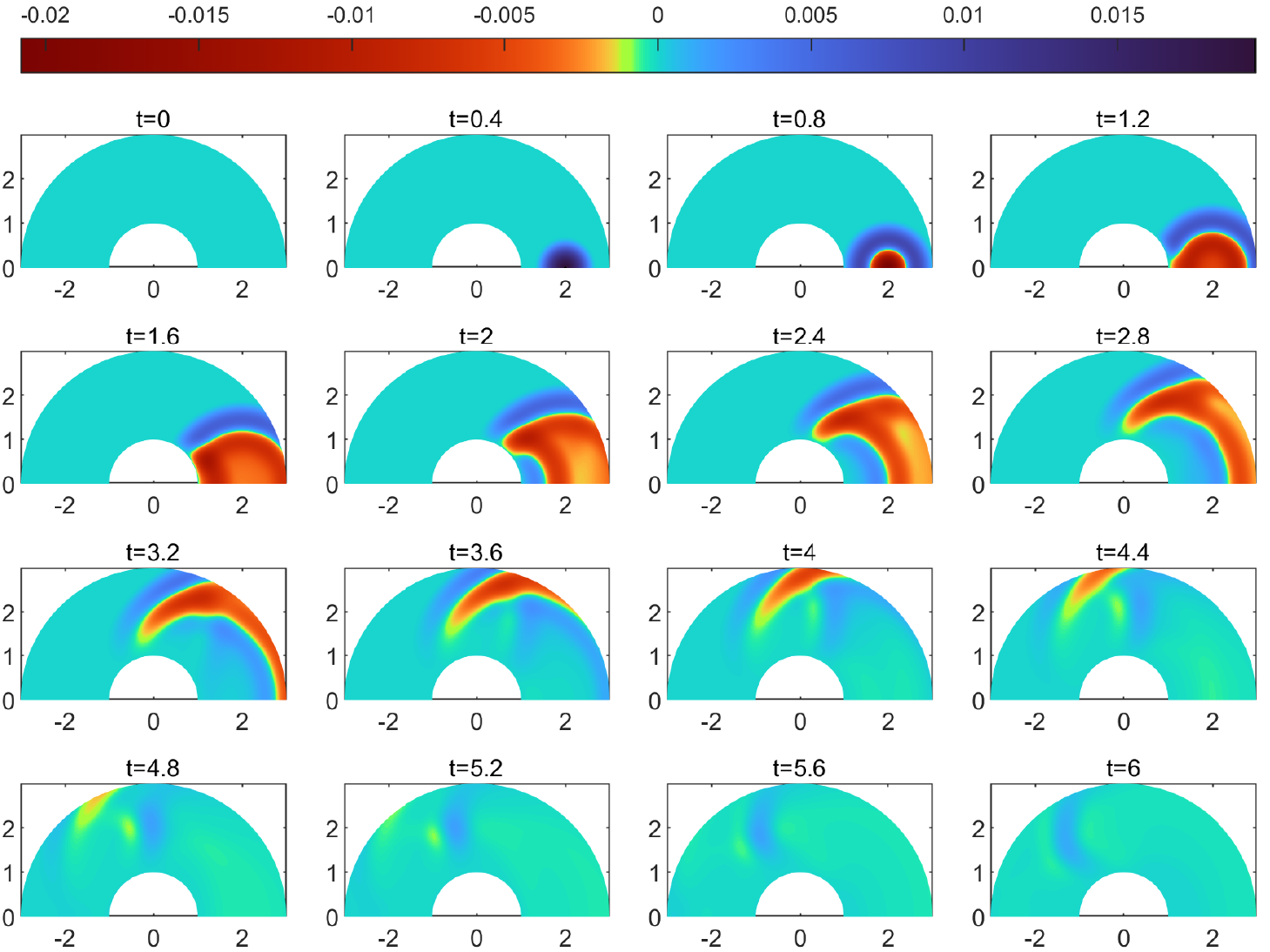}
	\caption{Snapshots of the solution of the bidimensional scattering problem \ref{ex:scattering}, obtained through the \igastab\ method with $p=5$, $h_s=0.0212$ and $h_t=0.0195$. The total number of DOFs is $N_{\mathrm{dof}}=14\,635\,660$.}
	\label{fig:scatteringplotsol}
\end{figure}

As an example of a problem in two space dimensions involving all three types of boundary conditions, we consider~\eqref{iga_wave_strong_eq} on $\Omega \times (0,6)$, with 
\begin{align*}
	\Omega &:= \left\{ (x,y) \in \mathbb{R} \times [0,+\infty) \ : \ 1 \leq x^2+y^2 \leq 9 \right\},\\
	\Gamma_D &:= \left\{ (x,y) \in \partial \Omega \ : \ x^2+y^2=1 \right\}, \\
	\Gamma_N &:= \left\{ (x,y) \in \partial \Omega \ : \ y=0 \right\}, \\
	\Gamma_R &:= \left\{ (x,y) \in  \partial \Omega \ : \ x^2+y^2=9 \right\},\\
	f(x,y,t) &:=\cos(2\pi t) \Psi(t) \Psi\left(\frac{\sqrt{(x-2)^2+y^2}}{0.4}\right),
\end{align*}
where $\Psi:\mathbb{R} \rightarrow \mathbb{R}$ denotes the bump function defined as
\begin{equation}\label{eq:bump}
	\Psi(s)= \begin{cases}
		\begin{aligned}
			&e^{1+\frac{1}{s^2-1}} & & s \in (-1,1),\\
			&0 & & \text{otherwise},
		\end{aligned}
	\end{cases}
\end{equation}
and $c = 1$, $\vartheta=1$ and $g_D=g_N=g_R=u_0=u_1 = 0$.
This problem represents the scattering of a pulse by a unit sound-soft disk; the impedance boundary condition can be thought as a low-order approximation of the radiation condition. The physical mesh is the image of a uniform mesh on the unit square through the geometric parametrization of $\Omega$. 
Figure~\ref{fig:scatteringplotsol} shows the reference numerical solution at different time instants.
\input{Graphics/ErrorScattering}
Figure~\ref{fig:ErrorScattering} shows optimal convergence rates for the errors in the $L^2(Q)$-norm and $H^1(Q)$-seminorm w.r.t.\ the reference numerical solution of Figure~\ref{fig:scatteringplotsol}.
As a consequence of the large number of DOFs ($N_{\mathrm{dof}}=14\,635\,660$), an iterative solver is mandatory.
Therefore, we solve the discrete problems by a preconditioned GMRES method with tolerance $10^{-12}$. 
The preconditioner is constructed similarly to that in~\cite{LMST2020} and it will be investigated in a future work.

\subsubsection{Wavefront in a non-homogenous medium}\label{sec:Csmooth}
\begin{figure}[htbp]
	\centering %\hfill
	\includegraphics[width=\linewidth]{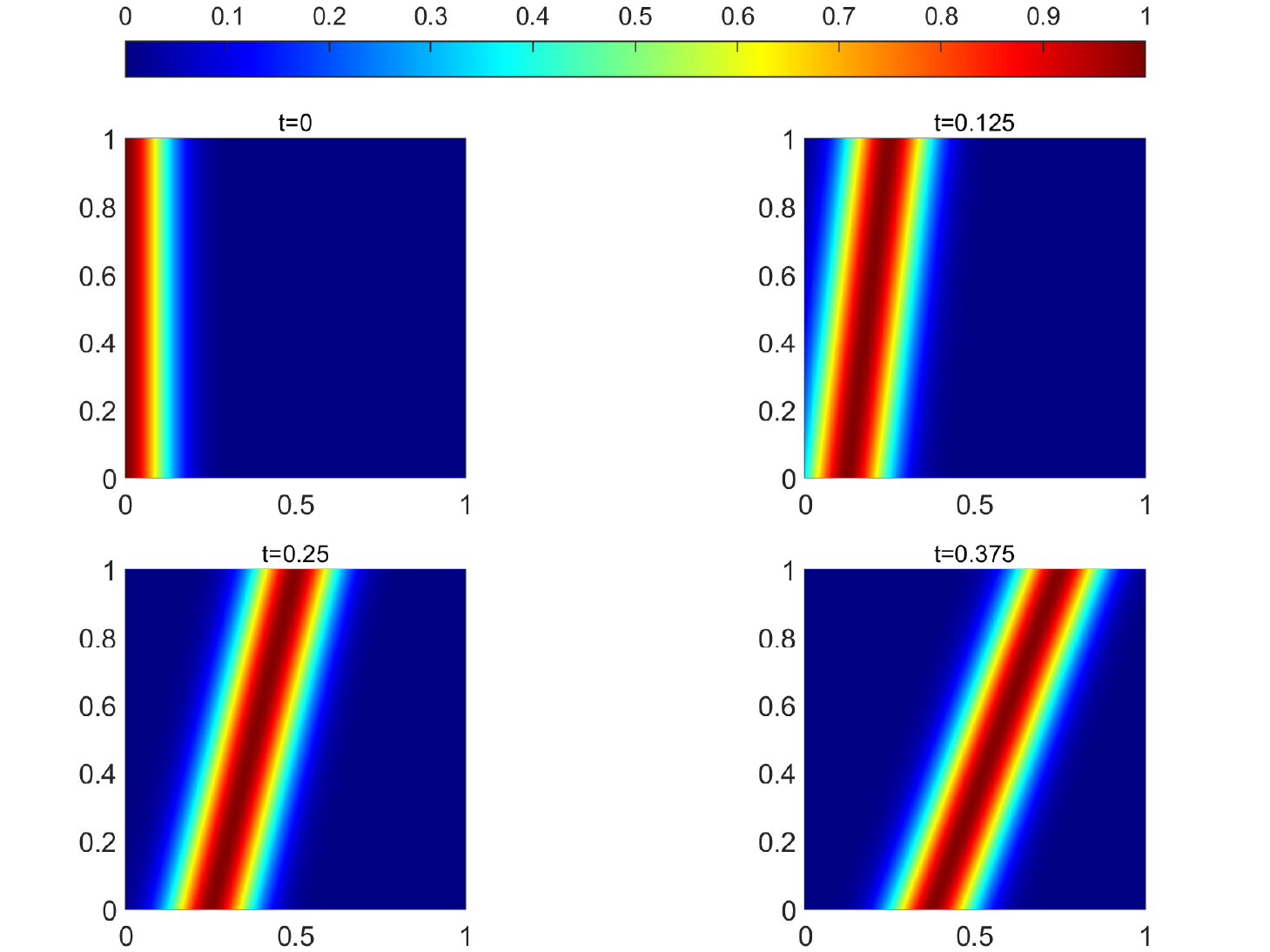}
	\caption{Snapshots of solution~\eqref{ucsmooth_ex} obtained through \igastab\ method with $p=4$, $h_s=h_t=0.0078$.}
	\label{fig:CsmoothplotsolRegEx}
\end{figure}

As an example of a problem with non-constant and smooth velocity, we consider the following exact solution of the acoustic wave equation~\eqref{iga_wave_strong_eq}
\begin{equation}\label{ucsmooth_ex} 
	u(x,y,t)=e^{-64(x-(1+y)t)^2} \quad \text{for } (x,t) \in \Omega \times (0,T):=(0,1)^2 \times (0,0.375),
\end{equation}
with velocity $c(x,y,t)=1+y$, $\Gamma_N = \partial \Omega$, and $g_N$, $u_0$, $u_1$, $f$ appropriately chosen.
This problem represents a vertical wavefront traveling through a medium with velocity varying along the $y$ direction.
Figure~\ref{fig:CsmoothplotsolRegEx} shows the reference numerical solution at different time instants.
As one can observe, the initial condition is constant along the $y$-direction. At $t > 0$, the wave
travels faster at the top of the domain than at the bottom in
accordance with the different values of $c$.
The relative errors in the $L^2(Q)$-norm and weighted $H^1(Q)$-seminorm \eqref{H1_sem} are shown in Figure~\ref{fig:ErrorCsmoothRegEx}, confirming optimal convergence rates of the proposed method.
\input{Graphics/ErrorCsmoothRegEx}

\subsubsection{Discontinuous velocity}\label{sec:Cdisc}
\begin{figure}[htbp]
	\centering
	\begin{subfigure}{0.49\textwidth}
	\includegraphics[width=\textwidth]{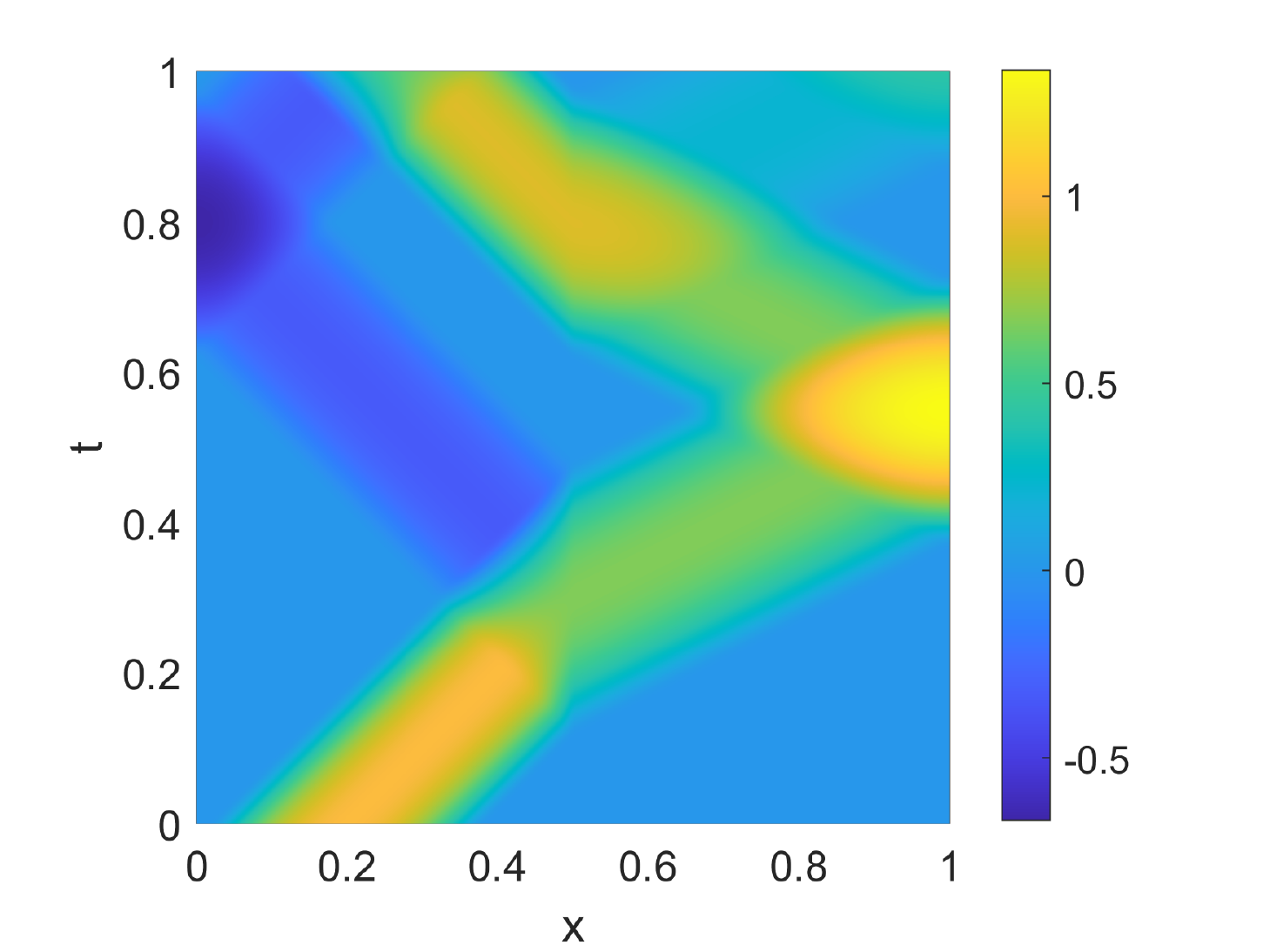}
	\end{subfigure}
	\\
	\begin{subfigure}{0.49\textwidth}
		\includegraphics[width=\textwidth]{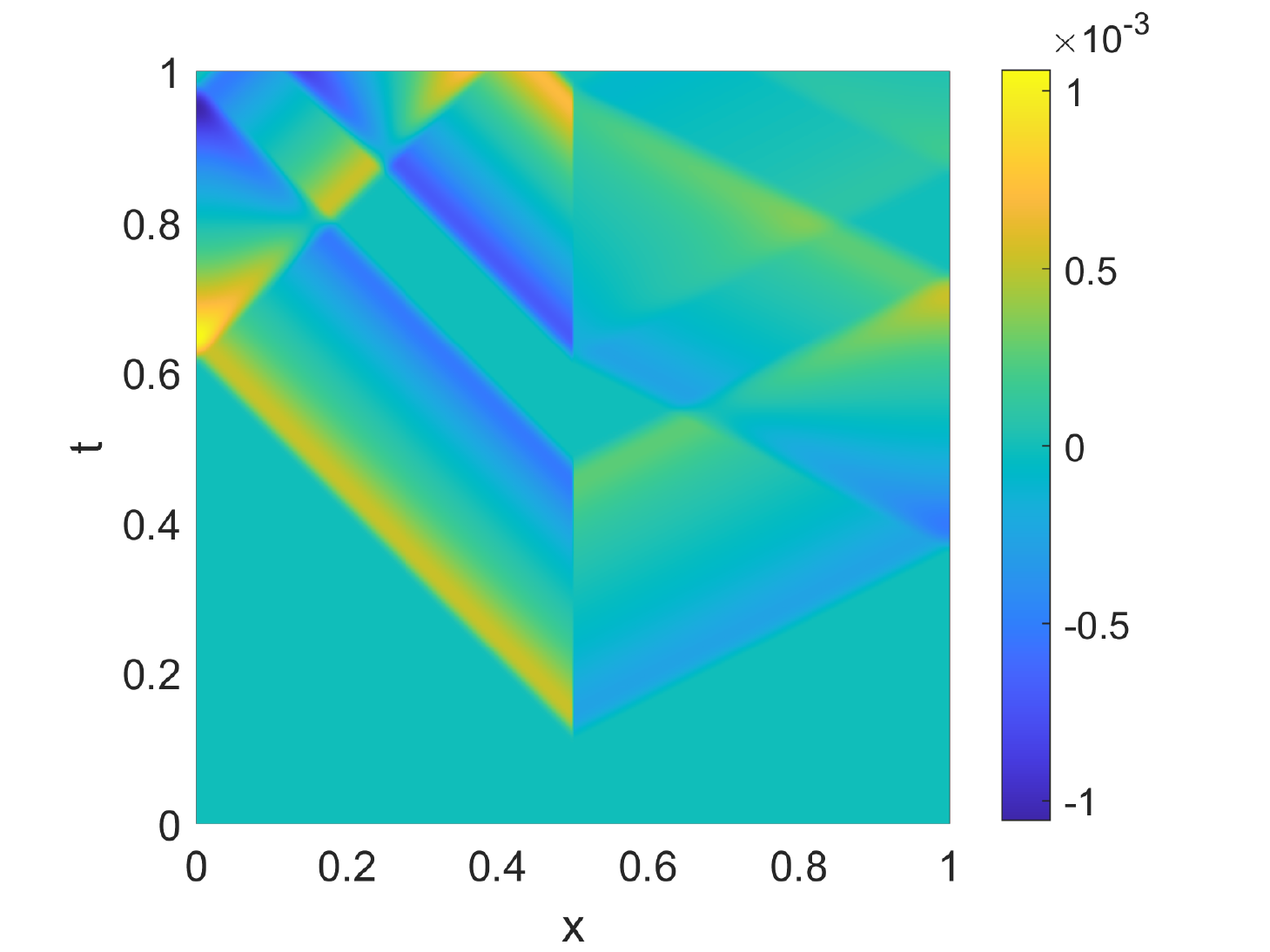}
	\end{subfigure}
	\hfill
	\begin{subfigure}{0.49\textwidth}
		\includegraphics[width=\textwidth]{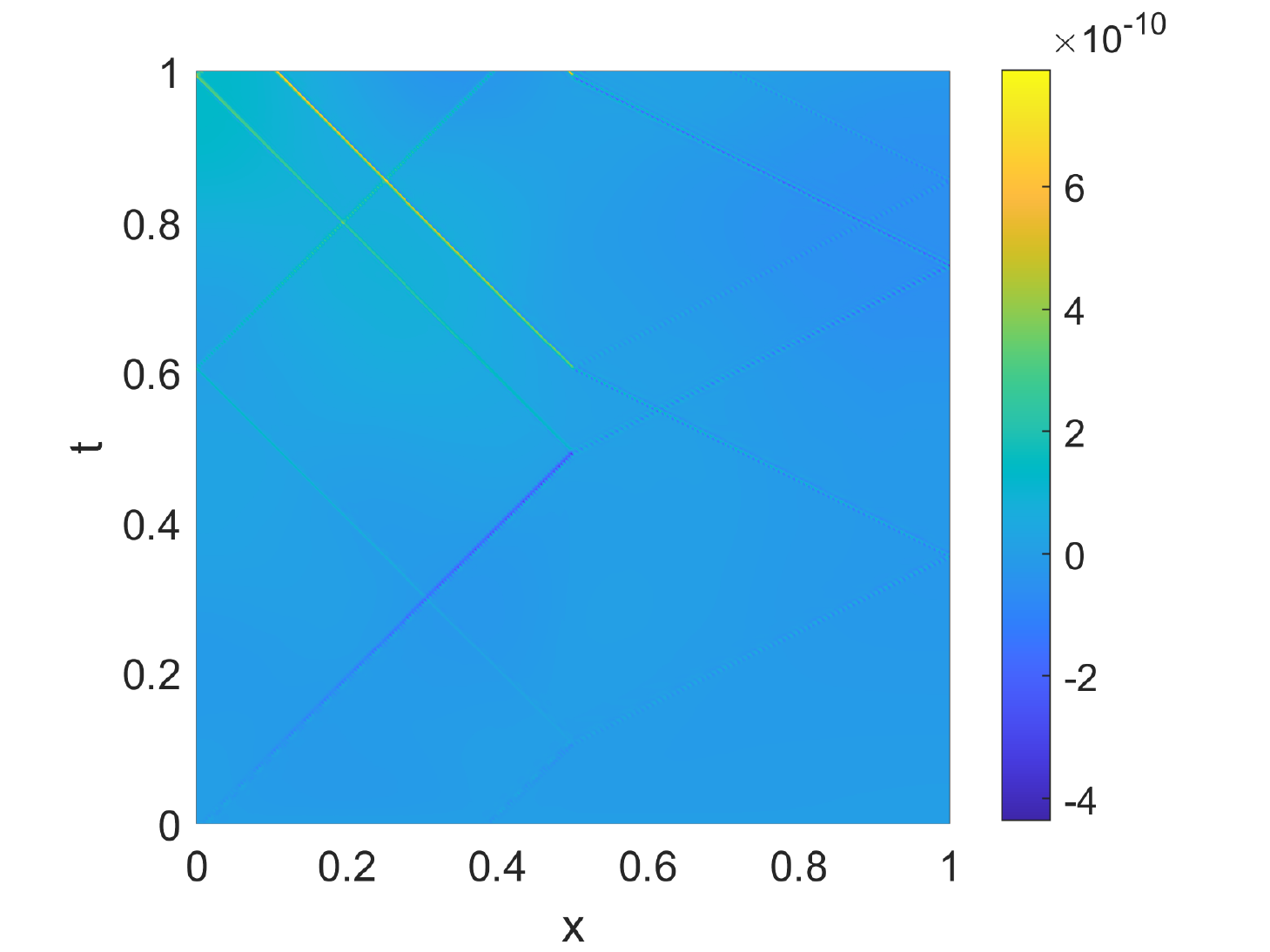}
	\end{subfigure}
	\caption{Top panel: exact solution of the wave propagation problem with piecewise-constant velocity~\eqref{cdisc_exsol}.
	\newline
	Bottom panels: difference between the exact and the numerical solution approximated by the \igastab\ method with $p=4$ and $h_s=h_t=0.000488$.
	Left: splines of maximal regularity in both space and time.
	Right: splines of space--time maximal regularity except at $x=0.5$, where we impose only $C^0$-continuity.
	The error due to the high continuity imposed at $x=0.5$ pollutes the solution in the whole domain.}
	\label{fig:CdiscBumpPlotSol}
\end{figure}

We consider an example of a wave propagating through a medium with discontinuous velocity: problem~\eqref{iga_wave_strong_eq} on $\Omega \times (0,T):=(0,1) \times (0,1)$, with $\Gamma_N = \partial \Omega$, $g_N=0$, $f=0$,
\begin{equation*}
	c(x,t)= \begin{cases}
	1 & 0\le x<\frac12, \\
	2 & \frac12\le x\le 1,
	\end{cases}
\end{equation*}
$u_0(x)=\Psi(5x-1)$ and $u_1(x)=-5 \Psi'(5x-1)$, with $\Psi$ the smooth bump defined %as
in~\eqref{eq:bump}. The exact solution is
\begin{align}\label{cdisc_exsol}
	&u(x,t) = 
	\\&\nonumber
	\begin{cases}
		\Psi \left( 5 (x-t)-1 \right) - \frac{1}{3} \Psi \left( 5 (x+t)- 4 \right) - \frac{1}{3} \Psi \left( 5 (x-t)+ 4 \right) - \frac{8}{9} \Psi \left( 5 (x+t)- \frac{13}{2} \right) & 0\le x<\frac{1}{2}, \\[2mm]
		\frac{2}{3} \left[ \Psi \left( \frac{5}{2}(x-2t)+\frac{1}{4} \right) + \Psi \left( \frac{5}{2}(x+2t)-\frac{21}{4} \right) \right] + \frac{2}{9} \left[\Psi \left( \frac{5}{2}(x-2t)+\frac{11}{4} \right) + \Psi \left( \frac{5}{2}(x+2t)-\frac{31}{4} \right) \right]& \frac12\le x\le1,
	\end{cases}
\end{align}
and its graph is shown in the top panel of Figure~\ref{fig:CdiscBumpPlotSol}. 
The top plots of Figure~\ref{fig:CdiscBumpStabConv2in1} show the $L^2(Q)$-norm and the weighted $H^1(Q)$-seminorm \eqref{H1_sem} relative errors obtained with splines of maximal regularity in both space and time. As we may expect from the low regularity of the solution~\eqref{cdisc_exsol} at $x=0.5$, the convergence rates are suboptimal for splines of degree $p>1$ and space--time regularity $C^{p-1}$. 
The bottom plots show the errors obtained with splines of space--time maximal regularity except at $x=0.5$, where we impose only $C^0$-continuity. As we can see, optimal convergence rates are achieved in this case.
Finally, the lower panels of Figure~\ref{fig:CdiscBumpPlotSol} show the difference between the exact solution~\eqref{cdisc_exsol} and the numerical solutions of the aforementioned methods with $p=4$ and $h_s=h_t=0.000488$.
The possibility to locally reduce the space regularity of the basis functions, without spoiling the unconditional stability, is a favorable feature of the \igastab\ method, in view of the high-order approximation of problems with point or line singularities due to the domain geometry or to the material coefficients.
\input{Graphics/ErrorCdiscBumpConvWide2in1}

\subsection{Energy conservation}\label{s:Energy}
Motivated by the loss of accuracy for long-time computations suffered by numerical schemes that increase or dissipate energy, we test how well the \igastab\ method~\eqref{equivstabiga} preserves the solution energy.  

As in \cite[Remark 4.2.36]{Z2020}, we consider the following exact solution of the acoustic wave equation \eqref{iga_wave_strong_eq} in one dimension:
	\begin{equation}\label{uex_energy}
	\begin{aligned}
	u
	(x,t)&=\big( \cos(\pi t) + \sin(\pi t) \big) \sin(\pi x) \quad \text{for } (x,t) \in Q:=(0,1) \times (0,10),\\
	\end{aligned}
	\end{equation}
	with $c= 1$, $f= 0$, $\Gamma_D = \partial \Omega$ and $g_D = 0$.
	The total energy of the exact solution \eqref{uex_energy}, i.e., 
	\begin{equation*}
	E(t):=\frac{1}{2}\|\partial_t
	 u(\cdot,t)\|^2_{L^2(\Omega)}+\frac{1}{2}\|\nabla u
	(\cdot,t)\|^2_{L^2(\Omega)} \quad \text{for } t \in [0,10]
	\end{equation*}
	 is constant in time and coincides with $\tfrac{\pi^2}{2}$ for all $t \in [0,10]$.
	We then compute the discrete energy 
	\begin{equation*}
	E_h(t):=\frac{1}{2}\|\partial_t u_{h}(\cdot,t)\|^2_{L^2(\Omega)}+\frac{1}{2}\|\nabla u_{h}(\cdot,t)\|^2_{L^2(\Omega)} \quad \text{for } t \in [0,10]
	\end{equation*}
	of the \igastab\ solution with mesh-sizes $h_s = 2^{-7}$, $h_t = h_s$, and stabilization coefficient $\delta = 10^{-p}$ for $p=1,2,3,4$.
	In Figure~\ref{fig:energy}, we can observe that moderately
        large values of $p$ allow to preserve the energy extremely
        well: the relative error $|E(t)-E_h(t)|/E(t)$ does not grow
        with time and is bounded by $10^{-2p}$ at all values of $t$,
        and is well preserved for all degrees $p$.
	
	\input{Graphics/StabEnergy}
	
\subsection{Dispersion properties}\label{s:Dispersion}
In order to study the numerical dispersion of \igastab\ method~\eqref{equivstabiga},
we solve a wave propagation problem across the segment $\Omega=(0,1)$, with periodic boundary conditions, for $T=2$, $c=1$ and $f=0$.
As initial conditions, we consider a $C^0$ tent profile, described by
\begin{equation}\label{eq:tentInitial}
	\begin{cases}
		\begin{aligned}
			u_0(x) &= (1-|4x-1|)\chi_{[0,\nicefrac{1}{2}]}(x) & & \text{for } x \in [0,1], \\
			u_1(x) &= -4\chi_{[0,\nicefrac{1}{4}]}(x) + 4 \chi_{[\nicefrac{1}{4},\nicefrac{1}{2}]}(x) & & \text{for } x \in [0,1],
		\end{aligned}
	\end{cases}	
\end{equation}
and a $C^\infty$ \B bump profile, described by
\begin{equation}\label{eq:bumpInitial}
	\begin{cases}
		\begin{aligned}
			u_0(x) &= \Psi(4x-1)\chi_{[0,\nicefrac{1}{2}]}(x) & & \text{for } x \in [0,1], \\
			u_1(x) &= -4\Psi'(4x-1)\chi_{[0,\nicefrac{1}{2}]}(x) & & \text{for } x \in [0,1],
		\end{aligned}
	\end{cases}	
\end{equation}
where $\Psi$ is defined as in~\eqref{eq:bump}, and both move rightward.

Figures~\ref{fig:dispersionTent} and \ref{fig:dispersionBump} show the discrete solutions and the corresponding errors at final time $T=2$, i.e.\ after two periods, for both the \igastab\ and the \femstab\ methods with equal number of DOFs.
Recall that, for $p=1$, both methods can be written as in \eqref{equivstabiga} and they differ only in the choice of $\delta=1/10$ or $\delta=1/12$.
Figures~\ref{fig:StabComparisonTent} and \ref{fig:StabComparisonBump} show the dependence of the error norm at final time on the polynomial degree $p$. 

The dispersion of a numerical method is related to the phase error (as opposed to the amplitude error measured by the dissipation) of the Fourier components of the solution, defined as
\begin{equation*}
	\left|\arg\left(\frac{c_n}{c_{n,h}} \cdot \frac{\|c_{n,h}\|}{\| c_{n}\|}\right)\right|,
\end{equation*}
where $c_n$ and $c_{n,h}$ denotes the $n$-th complex Fourier coefficients of the exact solution and the numerical one, respectively.
To better understand the dispersion property of the stabilized schemes, in Figures \ref{fig:dispersion1234TentRotByP} and \ref{fig:dispersion1234BumpRotByP}
we plot the time evolution of the phase error of the largest (in magnitude) Fourier coefficients of  the solution.

We observe that for both the \igastab\ and the \femstab\ methods, for $p>1$ the phase error grows moderately in time.
The error in the higher modes is larger for the rougher solution (case \eqref{eq:tentInitial}), but even in this case the \igastab\ is slightly less dispersive than the \femstab.

Finally, we observe that for the case \eqref{eq:tentInitial} the \femstab\ method with an even polynomial degree $p$ seems to perform worse than with an odd $p$. 
This is evident in Figures~\ref{fig:dispersionTent},~\ref{fig:StabComparisonTent} and~\ref{fig:dispersion1234TentRotByP}.

\input{Graphics/Dispersion_p1234Tent}
\input{Graphics/DispersionTestErrorsTent}
\input{Graphics/Dispersion_p1234Bump}
\input{Graphics/DispersionTestErrorsBump}

\input{Graphics/Dispersion_p1234TentRotByP}
\input{Graphics/Dispersion_p1234BumpRotByP}

%% file: Graphics/ErrorStabDeltaFine.tex
\tikzstyle{Linea1}=[dashed]
\tikzstyle{Linea2}=[thick]%,mark=*]
%\pgfplotsset{
	%    tick label style={font=\footnotesize},
	%    label style={font=\footnotesize},
	%    legend style={font=\footnotesize},
	%}
% 
\pgfplotscreateplotcyclelist{Lista1}{%
	{Linea2,Red,mark=*},
	% 	{Linea1,Red},
	{Linea2,Green,mark=triangle*},
	% 	{Linea1,Green},
	{Linea2,Cyan,mark=square*},
	% 	{Linea1,Cyan},
	{Linea2,Violet,mark=diamond*},
	% 	{Linea1,Violet}
	{Linea2,BurntOrange,mark=pentagon*},
	{Linea2,Brown,mark=o},
	{Linea2,Magenta,mark=triangle},
	{Linea2,Blue,mark=square},
	{Linea2,Red,only marks, mark = o, draw opacity=1, mark size = 5pt},
	{Linea2,Green,only marks, mark = o, draw opacity=1, mark size = 5pt},
	{Linea2,Cyan,only marks, mark = o, draw opacity=1, mark size = 5pt},
	{Linea2,Violet,only marks, mark = o, draw opacity=1, mark size = 5pt},
	{Linea2,BurntOrange,only marks, mark = o, draw opacity=1, mark size = 5pt},
	{Linea2,Brown,only marks, mark = o, draw opacity=1, mark size = 5pt},
	{Linea2,Magenta,only marks, mark = o, draw opacity=1, mark size = 5pt},
	{Linea2,Blue,only marks, mark = o, draw opacity=1, mark size = 5pt},
}

\pgfplotscreateplotcyclelist{Lista2}{%
	{Linea2,Red},
	% 	{Linea1,Red},
	{Linea2,Green},
	% 	{Linea1,Green},
	{Linea2,Cyan},
	% 	{Linea1,Cyan},
	{Linea2,Violet},
	% 	{Linea1,Violet}
	{Linea2,BurntOrange},
	{Linea2,Brown},
	{Linea2,Magenta},
	{Linea2,Blue},
	{Linea2,Red,mark = *, draw opacity=1, mark size = 2pt},
	{Linea2,Green,mark = triangle*, draw opacity=1, mark size = 2pt},
	{Linea2,Cyan,mark = square*, draw opacity=1, mark size = 2pt},
	{Linea2,Violet,mark = diamond*, draw opacity=1, mark size = 2pt},
	{Linea2,BurntOrange,mark = pentagon*, draw opacity=1, mark size = 2pt},
	{Linea2,Brown,mark = o, draw opacity=1, mark size = 2pt},
	{Linea2,Magenta,mark = triangle, draw opacity=1, mark size = 2pt},
	{Linea2,Blue,mark = square, draw opacity=1, mark size = 2pt},
}

\def \DATAFILE {Graphics/ErrorTableDeltaFine.csv}
%
%\pgfplotstablegetrowsof{\DATAFILE}
%\pgfmathtruncatemacro{\rownum}{\pgfplotsretval-1}
\newcommand*{\ReadOutElement}[4]{%
	\pgfplotstablegetelem{#2}{#3}\of{#1}%
	\let#4\pgfplotsretval
}

\pgfplotstableread[col sep=comma]{Graphics/ErrorTableDelta.csv}\DataTab

\begin{figure}[htbp]
	\centering
	\begin{subfigure}[t]{0.9\linewidth}
		\begin{tikzpicture}[%font=\footnotesize,
			%trim axis right
			]
			\begin{loglogaxis}[
				cycle list name=Lista2,
				width=.98\linewidth,
				height=0.5\linewidth,
				xlabel={$\delta$},
				ymax=10000,
				xminorticks=false,
				yminorticks=false,
				%max space between ticks=20,
				xtick={1,0.1,0.01,0.001,0.0001,0.00001,0.000001,0.0000001,0.00000001,0.000000001},
				ytick={10000,100,1,0.01,0.0001,0.000001,0.00000001},
				ylabel={$| u - u_h |_{H^1(Q)}/| u |_{H^1(Q)}$},
				legend columns=2,
				%legend pos=south east,
				legend style={at={(0.01,0.99)},anchor=north west},
				xmajorgrids=true,
				ymajorgrids=true,
				legend entries={,,,,,,,,$p=1$,$p=2$,$p=3$,$p=4$,$p=5$,$p=6$,$p=7$,$p=8$}]
				\addplot table [x=DeltaP1, y=H1sP1, col sep=comma] {\DATAFILE};
				\addplot table [x=DeltaP2, y=H1sP2, col sep=comma] {\DATAFILE};
				\addplot table [x=DeltaP3, y=H1sP3, col sep=comma] {\DATAFILE};
				\addplot table [x=DeltaP4, y=H1sP4, col sep=comma] {\DATAFILE};
				\addplot table [x=DeltaP5, y=H1sP5, col sep=comma] {\DATAFILE};
				\addplot table [x=DeltaP6, y=H1sP6, col sep=comma] {\DATAFILE};
				\addplot table [x=DeltaP7, y=H1sP7, col sep=comma] {\DATAFILE};
				\addplot table [x=DeltaP8, y=H1sP8, col sep=comma] {\DATAFILE};
				\ReadOutElement{\DataTab}{11}{Delta}{\xc}
				\ReadOutElement{\DataTab}{11}{H1sP1}{\yc}
				\addplot coordinates{(\xc,\yc)};
				\ReadOutElement{\DataTab}{10}{Delta}{\xc}
				\ReadOutElement{\DataTab}{10}{H1sP2}{\yc}
				\addplot coordinates{(\xc,\yc)};
				\ReadOutElement{\DataTab}{9}{Delta}{\xc}
				\ReadOutElement{\DataTab}{9}{H1sP3}{\yc}
				\addplot coordinates{(\xc,\yc)};
				\ReadOutElement{\DataTab}{8}{Delta}{\xc}
				\ReadOutElement{\DataTab}{8}{H1sP4}{\yc}
				\addplot coordinates{(\xc,\yc)};
				\ReadOutElement{\DataTab}{7}{Delta}{\xc}
				\ReadOutElement{\DataTab}{7}{H1sP5}{\yc}
				\addplot coordinates{(\xc,\yc)};
				\ReadOutElement{\DataTab}{6}{Delta}{\xc}
				\ReadOutElement{\DataTab}{6}{H1sP6}{\yc}
				\addplot coordinates{(\xc,\yc)};
				\ReadOutElement{\DataTab}{5}{Delta}{\xc}
				\ReadOutElement{\DataTab}{5}{H1sP7}{\yc}
				\addplot coordinates{(\xc,\yc)};
				\ReadOutElement{\DataTab}{4}{Delta}{\xc}
				\ReadOutElement{\DataTab}{4}{H1sP8}{\yc}
				\addplot coordinates{(\xc,\yc)};
			\end{loglogaxis}
		\end{tikzpicture}
	\end{subfigure}
	\caption{Relative errors of the \igastab\ method \eqref{equivstabiga} with $h_s=2^{-7}$ and $h_t=5h_s$, plotted against the stabilization parameter $\delta$. The exact solution is defined in \eqref{u1_ex}. The markers correspond to $\delta=10^{-p}$.}
	\label{fig:ErrorStabDelta}
\end{figure}

%% file: Graphics/ErrorStabRobust.tex
\tikzstyle{Linea1}=[dashed]
\tikzstyle{Linea2}=[thick]%,mark=*]

%\pgfplotsset{
	%    tick label style={font=\footnotesize},
	%    label style={font=\footnotesize},
	%    legend style={font=\footnotesize},
	%}
% 
% \pgfplotscreateplotcyclelist{Lista1}{%
	% 	{Linea2,Red},
	% 	%{Linea1,Red},
	% 	{Linea2,Green},
	% 	%{Linea1,Green},
	% 	{Linea2,Cyan},
	% 	%{Linea1,Cyan},
	% 	{Linea2,Violet},
	% 	%{Linea1,Violet}
	% }
\pgfplotscreateplotcyclelist{Lista1}{%
	{Linea2,Red,mark=*},
	% 	{Linea1,Red},
	{Linea2,Green,mark=triangle*},
	% 	{Linea1,Green},
	{Linea2,Cyan,mark=square*},
	% 	{Linea1,Cyan},
	{Linea2,Violet,mark=diamond*},
	% 	{Linea1,Violet}
}

\def \DATAFILE {Graphics/ErrorTableStabRobustRST.csv}

\begin{figure}[htbp]
	\centering
	\hspace*{\fill}
	\begin{subfigure}[t]{0.425\linewidth}
		\centering
		\begin{tikzpicture}[%font=\footnotesize, 
			trim axis left
			]
			\begin{loglogaxis}[
				cycle list name=Lista1,
				width=\linewidth,
				height=\linewidth,
				xlabel={$h_t/h_s$},
				ymin=0.00025,
				ymax=0.5,
				xminorticks=false,
				yminorticks=false,
				%max space between ticks=20,
				%xtick=data,
				%xticklabel={\pgfmathparse{exp(\tick)}\pgfmathprintnumber[fixed relative,precision=3]{\pgfmathresult}},
				%ytick={0.1,0.01,0.001,0.0001}, 
				%yticklabels={$0.4$,$0.6$,$0.8$,$1$}, 
				ylabel={$\| u - u_h \|_{L^2(Q)}/\| u \|_{L^2(Q)}$},
				legend columns=2,
				%legend pos=south east,
				legend style={at={(0.01,0.7)},anchor= north west},
				xmajorgrids=true,
				ymajorgrids=true,
				legend entries={$p=1$,$p=2$,$p=3$,$p=4$}]
				\addplot table [x=htfhs, y=L2P1R, col sep=comma] {\DATAFILE};
				\addplot table [x=htfhs, y=L2P2R, col sep=comma] {\DATAFILE};
				\addplot table [x=htfhs, y=L2P3R, col sep=comma] {\DATAFILE};
				\addplot table [x=htfhs, y=L2P4R, col sep=comma] {\DATAFILE};
			\end{loglogaxis}
		\end{tikzpicture}
	\end{subfigure}
	\hfill
	\begin{subfigure}[t]{0.425\linewidth}
		\centering
		\begin{tikzpicture}[%font=\footnotesize,
			trim axis right
			]
			\begin{loglogaxis}[
				cycle list name=Lista1,
				width=\linewidth,
				height=\linewidth,
				xlabel={$h_t/h_s$},
				ymin=0.00025,
				ymax=0.5,
				xminorticks=false,
				yminorticks=false,
				%max space between ticks=20,
				%xtick=data,
				%xticklabel={\pgfmathparse{exp(\tick)}\pgfmathprintnumber[fixed relative,precision=3]{\pgfmathresult}},
				ytick={0.1,0.01,0.001,0.0001}, 
				ylabel={$| u - u_h |_{H^1(Q)}/| u |_{H^1(Q)}$},
				legend columns=2,
				%legend pos=south east,
				legend style={at={(0.01,0.01)},anchor= south west},
				xmajorgrids=true,
				ymajorgrids=true,
				legend entries={$p=1$,$p=2$,$p=3$,$p=4$}]
				\addplot table [x=htfhs, y=H1P1R, col sep=comma] {\DATAFILE};
				\addplot table [x=htfhs, y=H1P2R, col sep=comma] {\DATAFILE};
				\addplot table [x=htfhs, y=H1P3R, col sep=comma] {\DATAFILE};
				\addplot table [x=htfhs, y=H1P4R, col sep=comma] {\DATAFILE};
			\end{loglogaxis}
		\end{tikzpicture}
	\end{subfigure}
	\hspace*{\fill}
	\caption{Relative errors of the \igastab\ method \eqref{equivstabiga} plotted against the ratio $h_t/h_s$ with fixed $h_t=0.1562$. The exact solution is defined in \eqref{u1_ex}.
	We observe that the method does not require any CFL condition: reducing the space mesh size, the error remains bounded.
	}
	\label{fig:StabRobust}
\end{figure}

%% file: Graphics/ErrorZankOurCompHP.tex
\tikzstyle{Linea1}=[thick,dashed, mark options={solid}]
\tikzstyle{Linea2}=[thick]%,mark=*]
\tikzstyle{Linea3}=[thick,dotted]

%\pgfplotsset{
	%    tick label style={font=\footnotesize},
	%    label style={font=\footnotesize},
	%    legend style={font=\footnotesize},
	%}

\pgfplotscreateplotcyclelist{Lista1}{%
	{Linea1,Red,mark=*},
	{Linea2,Red,mark=*},
	{Linea3,Red},
	{Linea1,Green,mark=triangle*},
	{Linea2,Green,mark=triangle*},
	{Linea3,Green},
	{Linea1,Cyan,mark=square*},
	{Linea2,Cyan,mark=square*},
	{Linea3,Cyan},
	{Linea1,Violet,mark=diamond*},
	{Linea2,Violet,mark=diamond*},
	{Linea3,Violet},
}

\pgfplotscreateplotcyclelist{Lista2}{%
	{Linea2,Red,mark=*},
	{Linea1,Red,mark=*},
	{Linea3,Red},
	{Linea2,Green,mark=triangle*},
	{Linea1,Green,mark=triangle*},
	{Linea3,Green},
	{Linea2,Cyan,mark=square*},
	{Linea1,Cyan,mark=square*},
	{Linea3,Cyan},
	{Linea2,Violet,mark=diamond*},
	{Linea1,Violet,mark=diamond*},
	{Linea3,Violet},
}

\pgfplotscreateplotcyclelist{Lista1P}{%
	{Linea2,Red,mark=*},
	{Linea1,Red,mark=*},
	%{Linea3,Red},
	{Linea2,Green,mark=triangle*},
	{Linea1,Green,mark=triangle*},
	%{Linea3,Green},
	{Linea2,Cyan,mark=square*},
	{Linea1,Cyan,mark=square*},
	%{Linea3,Cyan},
	{Linea2,Violet,mark=diamond*},
	{Linea1,Violet,mark=diamond*},
	%{Linea3,Violet},
}

\def \DATAFILE {Graphics/ErrorTableZankVsOur.csv}
\def \DATAFILEp {Graphics/ErrorTableZankVsOurPndof.csv}

\begin{figure}[htbp]
	\centering
	\hspace*{\fill}
	\begin{subfigure}[t]{0.425\linewidth}
		\centering
		\begin{tikzpicture}[%font=\footnotesize, 
			trim axis left
			]
			\begin{loglogaxis}[
				cycle list name=Lista1,
				yticklabel style={text width=2.125em, align=right},
				width=\linewidth,
				height=\linewidth,
				xlabel={$N_{\mathrm{dof}}$},
				ymax=10,
				ymin=1e-12,
				xminorticks=false,
				yminorticks=false,
				%max space between ticks=20,
				%xtick=data,
				%xticklabel={\pgfmathparse{exp(\tick)}\pgfmathprintnumber[fixed relative,precision=3]{\pgfmathresult}},
				%ytick={0.4,0.6,0.8,1}, 
				%yticklabels={$0.4$,$0.6$,$0.8$,$1$}, 
				ylabel={$\| u - u_h \|_{L^2(Q)}/\| u \|_{L^2(Q)}$},
				legend columns=2,
				%legend pos=south east,
				legend style={at={(0.01,0.01)},anchor=south west},
				xmajorgrids=true,
				ymajorgrids=true,
				legend entries={,$p=1$,$O(N_{\mathrm{dof}}^{-1})$,,$p=2$,$O(N_{\mathrm{dof}}^{-3/2})$,,$p=3$,$O(N_{\mathrm{dof}}^{-2})$,,$p=4$,$O(N_{\mathrm{dof}}^{-5/2})$}]
				\addplot table [x=NdofP1, y=L2P1Z, col sep=comma] {\DATAFILE};
				\addplot table [x=NdofP1, y=L2P1O, col sep=comma] {\DATAFILE};
				\addplot table [x=NdofP1, y=NdofL2P1O, col sep=comma] {\DATAFILE};
				\addplot table [x=NdofP2, y=L2P2Z, col sep=comma] {\DATAFILE};
				\addplot table [x=NdofP2, y=L2P2O, col sep=comma] {\DATAFILE};
				\addplot table [x=NdofP2, y=NdofL2P2O, col sep=comma] {\DATAFILE};
				\addplot table [x=NdofP3, y=L2P3Z, col sep=comma] {\DATAFILE};
				\addplot table [x=NdofP3, y=L2P3O, col sep=comma] {\DATAFILE};
				\addplot table [x=NdofP3, y=NdofL2P3O, col sep=comma] {\DATAFILE};
				\addplot table [x=NdofP4, y=L2P4Z, col sep=comma] {\DATAFILE};
				\addplot table [x=NdofP4, y=L2P4O, col sep=comma] {\DATAFILE};
				\addplot table [x=NdofP4, y=NdofL2P4O, col sep=comma] {\DATAFILE};
			\end{loglogaxis}
		\end{tikzpicture}
	\end{subfigure}
	\hfill
	\begin{subfigure}[t]{0.425\linewidth}
		\centering
		\begin{tikzpicture}[%font=\footnotesize,
			trim axis right
			]
			\begin{loglogaxis}[
				cycle list name=Lista1,
				yticklabel style={text width=2.125em, align=right},
				width=\linewidth,
				height=\linewidth,
				xlabel={$N_{\mathrm{dof}}$},
				ymax=10,
				ymin=1e-12,
				xminorticks=false,
				yminorticks=false,
				%max space between ticks=20,
				%xtick=data,
				%xticklabel={\pgfmathparse{exp(\tick)}\pgfmathprintnumber[fixed relative,precision=3]{\pgfmathresult}},
				%ytick={0.4,0.6,0.8,1}, 
				%yticklabels={$0.4$,$0.6$,$0.8$,$1$}, 
				ylabel={$| u - u_h |_{H^1(Q)}/| u |_{H^1(Q)}$},
				legend columns=2,
				%legend pos=south east,
				legend style={at={(0.01,0.01)},anchor=south west},
				xmajorgrids=true,
				ymajorgrids=true,
				legend entries={,$p=1$,$O(N_{\mathrm{dof}}^{-1/2})$,,$p=2$,$O(N_{\mathrm{dof}}^{-1})$,,$p=3$,$O(N_{\mathrm{dof}}^{-3/2})$,,$p=4$,$O(N_{\mathrm{dof}}^{-2})$}]
				\addplot table [x=NdofP1, y=H1sP1Z, col sep=comma] {\DATAFILE};
				\addplot table [x=NdofP1, y=H1sP1O, col sep=comma] {\DATAFILE};
				\addplot table [x=NdofP1, y=NdofH1sP1O, col sep=comma] {\DATAFILE};
				\addplot table [x=NdofP2, y=H1sP2Z, col sep=comma] {\DATAFILE};
				\addplot table [x=NdofP2, y=H1sP2O, col sep=comma] {\DATAFILE};
				\addplot table [x=NdofP2, y=NdofH1sP2O, col sep=comma] {\DATAFILE};
				\addplot table [x=NdofP3, y=H1sP3Z, col sep=comma] {\DATAFILE};
				\addplot table [x=NdofP3, y=H1sP3O, col sep=comma] {\DATAFILE};
				\addplot table [x=NdofP3, y=NdofH1sP3O, col sep=comma] {\DATAFILE};
				\addplot table [x=NdofP4, y=H1sP4Z, col sep=comma] {\DATAFILE};
				\addplot table [x=NdofP4, y=H1sP4O, col sep=comma] {\DATAFILE};
				\addplot table [x=NdofP4, y=NdofH1sP4O, col sep=comma] {\DATAFILE};
			\end{loglogaxis}
		\end{tikzpicture}
	\end{subfigure}
	%%%%%%%%%%% p-convergence
\hspace*{\fill}
\\
\hspace*{\fill}
	\begin{subfigure}[t]{0.425\linewidth}
		\centering
		\begin{tikzpicture}[%font=\footnotesize, 
			trim axis left
			]
			\begin{semilogyaxis}[
				yticklabel style={text width=2.125em, align=right},
				cycle list name=Lista2,
				width=\linewidth,
				height=\linewidth,
				xlabel={$p$},
				ymin=1e-7,
				ymax=5e-1,
				xminorticks=false,
				yminorticks=false,
				%max space between ticks=20,
				%xtick=data,
				%xticklabel={\pgfmathparse{exp(\tick)}\pgfmathprintnumber[fixed relative,precision=3]{\pgfmathresult}},
				%ytick={0.4,0.6,0.8,1}, 
				%yticklabels={$0.4$,$0.6$,$0.8$,$1$}, 
				ylabel={$\| u(T) - u_h(T) \|_{L^2(\Omega)}/\| u(T) \|_{L^2(\Omega)}$},
				legend columns=1,
				%legend pos=south east,
				legend style={at={(0.01,0.01)},anchor=south west},
				xmajorgrids=true,
				ymajorgrids=true,
				xtick={1,2,3,4,5,6},
				legend entries={\igastab,\femstab}]
				\addplot table [x=p, y=L2O, col sep=comma] {\DATAFILEp};
				\addplot table [x=p, y=L2Z, col sep=comma] {\DATAFILEp};
				%\addplot table [x=NdofP1, y=NdofL2P1O, col sep=comma] {\DATAFILE};
			\end{semilogyaxis}
		\end{tikzpicture}
	\end{subfigure}
	\hfill
	\begin{subfigure}[t]{0.425\linewidth}
		\centering
		\begin{tikzpicture}[%font=\footnotesize,
			trim axis right
			]
			\begin{semilogyaxis}[
				yticklabel style={text width=2.125em, align=right},
				cycle list name=Lista2,
				width=\linewidth,
				height=\linewidth,
				xlabel={$p$},
				ymin=1e-7,
				ymax=5e-1,
				xminorticks=false,
				yminorticks=false,
				%max space between ticks=20,
				%xtick=data,
				%xticklabel={\pgfmathparse{exp(\tick)}\pgfmathprintnumber[fixed relative,precision=3]{\pgfmathresult}},
				%ytick={0.4,0.6,0.8,1}, 
				%yticklabels={$0.4$,$0.6$,$0.8$,$1$}, 
				ylabel={$| u(T) - u_h(T) |_{H^1(\Omega)}/| u(T) |_{H^1(\Omega)}$},
				legend columns=1,
				%legend pos=south east,
				legend style={at={(0.01,0.01)},anchor=south west},
				xmajorgrids=true,
				ymajorgrids=true,
				xtick={1,2,3,4,5,6},
				legend entries={\igastab,\femstab}]
				\addplot table [x=p, y=H1sO, col sep=comma] {\DATAFILEp};
				\addplot table [x=p, y=H1sZ, col sep=comma] {\DATAFILEp};
			\end{semilogyaxis}
		\end{tikzpicture}
	\end{subfigure}
	\hspace*{\fill}
	\caption{
	Comparison between relative errors of the \igastab\ method~\eqref{equivstabiga}, with splines of maximal regularity (continuous lines), and relative errors of the \femstab\ method~\eqref{stabZank}, with $C^0$ splines (dashed lines).
	\newline
	Top: space--time errors plotted against the total number of DOFs $N_{\mathrm{dof}}$, with $h_t = 5 h_s$.
	\newline
	Bottom: final-time errors plotted against the spline degree $p$, with $N_{\mathrm{dof}}=7\,080$ and $h_t \approx 5 h_s$ for both the stabilizations and all the spline degrees.
	}
	\label{fig:StabComparisonHP}
\end{figure}

%% file: Graphics/GraficoErroreLdof_p_1234_comparison.tex
\tikzstyle{Linea1}=[thick,dashed]
\tikzstyle{Linea2}=[thick,mark=*]

%\pgfplotsset{
	%    tick label style={font=\footnotesize},
	%    label style={font=\footnotesize},
	%    legend style={font=\footnotesize},
	%}

\def\HeiFact{.78}

\pgfplotscreateplotcyclelist{Lista0}{%
	{Linea2,Red},
	{Linea2,Green},
	{Linea2,Cyan},
	{Linea2,Violet},
	{Linea2,BurntOrange},
	{Linea2,Brown},
	{Linea2,Magenta},
	{Linea2,Blue},
}

\pgfplotscreateplotcyclelist{Lista1}{%
	{dashed,thick,mark=*,Red},
	{dashed,thick,mark=triangle*,Green},
	{dashed,thick,mark=square*,Cyan},
	{dashed,thick,mark=diamond*,Violet},
	{dashed,thick,mark=pentagon*,BurntOrange},
	{thick,mark=*,Red},
	{thick,mark=triangle*,Green},
	{thick,mark=square*,Cyan},
	{thick,mark=diamond*,Violet},
	{thick,mark=pentagon*,BurntOrange},
	{dashed,black,thick},
	{Linea2,Brown},
	{Linea2,Magenta},
	{Linea2,Blue},
}

\def \DATAFILEa {Graphics/ErrorTableZankLdofsP1.csv}
\def \DATAFILEb {Graphics/ErrorTableZankLdofsP2.csv}
\def \DATAFILEc {Graphics/ErrorTableZankLdofsP3.csv}
\def \DATAFILEd {Graphics/ErrorTableZankLdofsP4.csv}

\def \DATAFILEai {Graphics/ErrorTableLdofs_p1.csv}
\def \DATAFILEbi {Graphics/ErrorTableLdofs_p2.csv}
\def \DATAFILEci {Graphics/ErrorTableLdofs_p3.csv}
\def \DATAFILEdi {Graphics/ErrorTableLdofs_p4.csv}

\begin{figure}[p]
	\centering 
\hspace*{\fill}
	\begin{subfigure}[t]{0.425\linewidth}
		\centering
		\begin{tikzpicture}[trim axis left] %font=\footnotesize, 
			\begin{loglogaxis}
				[cycle list name=Lista1,
				width=\linewidth,
				height=\HeiFact\linewidth,
				xminorticks=false, yminorticks=false,
				ylabel = {\fbox{$p=1$}},
				title={$\| u - u_h \|_{L^2(Q)}/\| u \|_{L^2(Q)}$},
				legend columns=1,
				legend pos=south east,
				legend style={at={(0.01,0.01)},anchor=south west},
				xmajorgrids=true, ymajorgrids=true,
				ymin={0.000000000001}, ymax=10, 
				xmin=0.05,
				xmax=20000,
				% ytick={1,10^(-3),10^(-6),10^(-9),10^(-12)},  % Works on windows, not on Linux...
				%ytick={1,.001,.000001,.000000001,.000000000001}, 
				legend entries={,,,,,$k=1$,$k=2$,$k=4$,$k=8$,$k=16$,$O(N_{\mathrm{dof}}^{-1})$}    ]
				\addplot table [x=nDofsK1, y=L2P1K1, col sep=comma] {\DATAFILEa};
				\addplot table [x=nDofsK2, y=L2P1K2, col sep=comma] {\DATAFILEa};
				\addplot table [x=nDofsK4, y=L2P1K4, col sep=comma] {\DATAFILEa};
				\addplot table [x=nDofsK8, y=L2P1K8, col sep=comma] {\DATAFILEa};
				\addplot table [x=nDofsK16, y=L2P1K16, col sep=comma] {\DATAFILEa};
				\addplot table [x=nDofsK1, y=L2p1k1, col sep=comma] {\DATAFILEai};
				\addplot table [x=nDofsK2, y=L2p1k2, col sep=comma] {\DATAFILEai};
				\addplot table [x=nDofsK4, y=L2p1k4, col sep=comma] {\DATAFILEai};
				\addplot table [x=nDofsK8, y=L2p1k8, col sep=comma] {\DATAFILEai};
				\addplot table [x=nDofsK16, y=L2p1k16, col sep=comma] {\DATAFILEai};
				\addplot[dashed,black,thick] table [x=nDofsK1, y expr=0.1*\thisrow{h2L2}, col sep=comma] {\DATAFILEa};
				%\addplot[dashed,black,thick] table [x=nDofsK1, y expr=10*\thisrow{L2h2}, col sep=comma] {\DATAFILEa};
			\end{loglogaxis}
		\end{tikzpicture}
	\end{subfigure}\hfill
	\begin{subfigure}[t]{0.425\linewidth}
		\centering
		\begin{tikzpicture}[%font=\footnotesize,
			trim axis right]
			\begin{loglogaxis}[
				cycle list name=Lista1,
				width=\linewidth, 
				height=\HeiFact\linewidth,
				xminorticks=false, yminorticks=false,
				%ylabel
				title={$| u - u_h |_{H^1(Q)}/| u |_{H^1(Q)}$ },
				legend columns=1,
				legend pos=south east,
				legend style={at={(0.01,0.01)},anchor=south west},
				xmajorgrids=true, ymajorgrids=true,
				ymin={0.000000000001}, ymax=10,
				xmin=0.05,
				xmax=20000,
				legend entries={,,,,,$k=1$,$k=2$,$k=4$,$k=8$,$k=16$,$O(N_{\mathrm{dof}}^{-1/2})$}]
				\addplot table [x=nDofsK1, y=H1sP1K1, col sep=comma] {\DATAFILEa};
				\addplot table [x=nDofsK2, y=H1sP1K2, col sep=comma] {\DATAFILEa};
				\addplot table [x=nDofsK4, y=H1sP1K4, col sep=comma] {\DATAFILEa};
				\addplot table [x=nDofsK8, y=H1sP1K8, col sep=comma] {\DATAFILEa};
				\addplot table [x=nDofsK16, y=H1sP1K16, col sep=comma] {\DATAFILEa};
				\addplot table [x=nDofsK1, y=H1p1k1, col sep=comma] {\DATAFILEai};
				\addplot table [x=nDofsK2, y=H1p1k2, col sep=comma] {\DATAFILEai};
				\addplot table [x=nDofsK4, y=H1p1k4, col sep=comma] {\DATAFILEai};
				\addplot table [x=nDofsK8, y=H1p1k8, col sep=comma] {\DATAFILEai};
				\addplot table [x=nDofsK16, y=H1p1k16, col sep=comma] {\DATAFILEai};
				\addplot[dashed,black,thick] table [x=nDofsK1, y expr=0.1*\thisrow{h1H1s}, col sep=comma] {\DATAFILEa};
				%\addplot[dashed,black,thick] table [x=nDofsK1, y expr=10*\thisrow{H1h1}, col sep=comma] {\DATAFILEa};
			\end{loglogaxis}
		\end{tikzpicture}
	\end{subfigure}
	%%%%%%%%%%%%%%%%%%%%%%%%%%%%%%%%%%%%%%%%%%%%%%%%%%%%%%%%%%%%%%%%%%%%%%%%%%%%%%%%%%%%%%%%%%%%%%%%%%%%%%%%%%%%%%%%%%%%%%%%%%%%%
	% P2 
	\hspace*{\fill}
	\\
	\hspace*{\fill}
	\begin{subfigure}[t]{0.425\linewidth}
		\centering
		\begin{tikzpicture}[%font=\footnotesize, 
			trim axis left]
			\begin{loglogaxis}[
				cycle list name=Lista1,
				width=\linewidth,
				height=\HeiFact\linewidth,
				xminorticks=false, yminorticks=false,
				ylabel = {\fbox{$p=2$}},
				legend columns=1,
				legend pos=south east,
				legend style={at={(0.01,0.01)},anchor=south west},
				xmajorgrids=true, ymajorgrids=true,
				ymin={0.000000000001}, ymax=10,
				xmin=0.05,
				xmax=20000, 
				legend entries={,,,,,,,,,,$O(N_{\mathrm{dof}}^{-3/2})$}]
				\addplot table [x=nDofsK1, y=L2P2K1, col sep=comma] {\DATAFILEb};
				\addplot table [x=nDofsK2, y=L2P2K2, col sep=comma] {\DATAFILEb};
				\addplot table [x=nDofsK4, y=L2P2K4, col sep=comma] {\DATAFILEb};
				\addplot table [x=nDofsK8, y=L2P2K8, col sep=comma] {\DATAFILEb};
				\addplot table [x=nDofsK16, y=L2P2K16, col sep=comma] {\DATAFILEb};
				\addplot table [x=nDofsK1, y=L2p2k1, col sep=comma] {\DATAFILEbi};
				\addplot table [x=nDofsK2, y=L2p2k2, col sep=comma] {\DATAFILEbi};
				\addplot table [x=nDofsK4, y=L2p2k4, col sep=comma] {\DATAFILEbi};
				\addplot table [x=nDofsK8, y=L2p2k8, col sep=comma] {\DATAFILEbi};
				\addplot table [x=nDofsK16, y=L2p2k16, col sep=comma] {\DATAFILEbi};
				\addplot[dashed,black,thick] table [x=nDofsK1, y expr=0.01*\thisrow{h3L2}, col sep=comma] {\DATAFILEb};
				%\addplot[dashed,black,thick] table [x=nDofsK1, y expr=10*\thisrow{L2h3}, col sep=comma] {\DATAFILEb};
			\end{loglogaxis}
		\end{tikzpicture}
	\end{subfigure}\hfill
	\begin{subfigure}[t]{0.425\linewidth}
		\centering
		\begin{tikzpicture}[%font=\footnotesize,
			trim axis right]
			\begin{loglogaxis}[
				cycle list name=Lista1,
				width=\linewidth,
				height=\HeiFact\linewidth,
				xminorticks=false, yminorticks=false,
				legend columns=1,
				legend pos=south east,
				legend style={at={(0.01,0.01)},anchor=south west},
				xmajorgrids=true, ymajorgrids=true,
				ymin={0.000000000001}, ymax=10,
				xmin=0.05,
				xmax=20000, 
				legend entries={,,,,,,,,,,$O(N_{\mathrm{dof}}^{-1})$}]
				\addplot table [x=nDofsK1, y=H1sP2K1, col sep=comma] {\DATAFILEb};
				\addplot table [x=nDofsK2, y=H1sP2K2, col sep=comma] {\DATAFILEb};
				\addplot table [x=nDofsK4, y=H1sP2K4, col sep=comma] {\DATAFILEb};
				\addplot table [x=nDofsK8, y=H1sP2K8, col sep=comma] {\DATAFILEb};
				\addplot table [x=nDofsK16, y=H1sP2K16, col sep=comma] {\DATAFILEb};
				\addplot table [x=nDofsK1, y=H1p2k1, col sep=comma] {\DATAFILEbi};
				\addplot table [x=nDofsK2, y=H1p2k2, col sep=comma] {\DATAFILEbi};
				\addplot table [x=nDofsK4, y=H1p2k4, col sep=comma] {\DATAFILEbi};
				\addplot table [x=nDofsK8, y=H1p2k8, col sep=comma] {\DATAFILEbi};
				\addplot table [x=nDofsK16, y=H1p2k16, col sep=comma] {\DATAFILEbi};
				\addplot[dashed,black,thick] table [x=nDofsK1, y expr=0.01*\thisrow{h2H1s}, col sep=comma] {\DATAFILEb};
				%\addplot[dashed,black,thick] table [x=nDofsK1, y expr=10*\thisrow{H1h2}, col sep=comma] {\DATAFILEb};
			\end{loglogaxis}
		\end{tikzpicture}
	\end{subfigure}
	%%%%%%%%%%%%%%%%%%%%%%%%%%%%%%%%%%%%%%%%%%%%%%%%%%%%%%%%%%%%%%%%%%%%%%%%%%%%%%%%%%%%%%%%%%%%%%%%%%%%%%%%%%%%%%%%%%%%%%%%%%%%%%%%%%%%%%%%%%%%
	% P3 
\hspace*{\fill}
\\
\hspace*{\fill}
	\begin{subfigure}[t]{0.425\linewidth}
		\centering
		\begin{tikzpicture}[%font=\footnotesize, 
			trim axis left]
			\begin{loglogaxis}[
				cycle list name=Lista1,
				width=\linewidth,
				height=\HeiFact\linewidth,
				xminorticks=false, yminorticks=false,
				xmin=0.05,
				xmax=20000,
				ylabel = {\fbox{$p=3$}},
				legend columns=1,
				legend pos=south east,
				legend style={at={(0.01,0.01)},anchor=south west},
				xmajorgrids=true, ymajorgrids=true,
				ymin={0.000000000001}, ymax=10, %ytick={1,.001,.000001,.000000001,.000000000001}, 
				legend entries={,,,,,,,,,,$O(N_{\mathrm{dof}}^{-2})$}]
				\addplot table [x=nDofsK1, y=L2P3K1, col sep=comma] {\DATAFILEc};
				\addplot table [x=nDofsK2, y=L2P3K2, col sep=comma] {\DATAFILEc};
				\addplot table [x=nDofsK4, y=L2P3K4, col sep=comma] {\DATAFILEc};
				\addplot table [x=nDofsK8, y=L2P3K8, col sep=comma] {\DATAFILEc};
				\addplot table [x=nDofsK16, y=L2P3K16, col sep=comma] {\DATAFILEc};
				\addplot table [x=nDofsK1, y=L2p3k1, col sep=comma] {\DATAFILEci};
				\addplot table [x=nDofsK2, y=L2p3k2, col sep=comma] {\DATAFILEci};
				\addplot table [x=nDofsK4, y=L2p3k4, col sep=comma] {\DATAFILEci};
				\addplot table [x=nDofsK8, y=L2p3k8, col sep=comma] {\DATAFILEci};
				\addplot table [x=nDofsK16, y=L2p3k16, col sep=comma] {\DATAFILEci};
				\addplot[dashed,black,thick] table [x=nDofsK1, y=h4L2, col sep=comma] {\DATAFILEc};
			\end{loglogaxis}
		\end{tikzpicture}
	\end{subfigure}\hfill
	\begin{subfigure}[t]{0.425\linewidth}
		\centering
		\begin{tikzpicture}[%font=\footnotesize,
			trim axis right]
			\begin{loglogaxis}[
				cycle list name=Lista1,
				width=\linewidth, height=\HeiFact\linewidth,
				xminorticks=false, yminorticks=false,
				legend columns=1,
				legend pos=south east,
				legend style={at={(0.01,0.01)},anchor=south west},
				xmajorgrids=true, ymajorgrids=true,
				xmin=0.05,
				xmax=20000,
				ymin={0.000000000001}, ymax=10, %ytick={1,.001,.000001,.000000001,.000000000001}, 
				legend entries={,,,,,,,,,,$O(N_{\mathrm{dof}}^{-3/2})$}]
				\addplot table [x=nDofsK1, y=H1sP3K1, col sep=comma] {\DATAFILEc};
				\addplot table [x=nDofsK2, y=H1sP3K2, col sep=comma] {\DATAFILEc};
				\addplot table [x=nDofsK4, y=H1sP3K4, col sep=comma] {\DATAFILEc};
				\addplot table [x=nDofsK8, y=H1sP3K8, col sep=comma] {\DATAFILEc};
				\addplot table [x=nDofsK16, y=H1sP3K16, col sep=comma] {\DATAFILEc};
				\addplot table [x=nDofsK1, y=H1p3k1, col sep=comma] {\DATAFILEci};
				\addplot table [x=nDofsK2, y=H1p3k2, col sep=comma] {\DATAFILEci};
				\addplot table [x=nDofsK4, y=H1p3k4, col sep=comma] {\DATAFILEci};
				\addplot table [x=nDofsK8, y=H1p3k8, col sep=comma] {\DATAFILEci};
				\addplot table [x=nDofsK16, y=H1p3k16, col sep=comma] {\DATAFILEci};
				\addplot[dashed,black,thick] table [x=nDofsK1, y expr=10*\thisrow{h3H1s}, col sep=comma] {\DATAFILEc};
				%\addplot[dashed,black,thick] table [x=nDofsK1, y expr=10*\thisrow{H1h3}, col sep=comma] {\DATAFILEc};
			\end{loglogaxis}
		\end{tikzpicture}
	\end{subfigure}
	%%%%%%%%%%%%%%%%%%%%%%%%%%%%%%%%%%%%%%%%%%%%%%%%%%%%%%%%%%%%%%%%%%%%%%%%%%%%%%%%%%%%%%%%%%%%%%%%%%%%%%%%%%%%%%%%%%%%%%%%%%%%%%%%%%%%%%%%%%%%
	% P4
	\hspace*{\fill}
	\\
	\hspace*{\fill}
	\begin{subfigure}[t]{0.425\linewidth}
		\centering
		\begin{tikzpicture}[%font=\footnotesize, 
			trim axis left]
			\begin{loglogaxis}[
				cycle list name=Lista1,
				width=\linewidth, height=\HeiFact\linewidth,
				xlabel={$N_{\mathrm{dof}}/\sharp \lambda$},
				xminorticks=false, yminorticks=false,
				ylabel = {\fbox{$p=4$}},
				legend columns=1,
				legend pos=south east,
				legend style={at={(0.01,0.01)},anchor=south west},
				xmajorgrids=true, ymajorgrids=true,
				ymin={2.5e-13}, ymax=10,
				xmin=0.05,
				xmax=20000,
				legend entries={,,,,,,,,,,$O(N_{\mathrm{dof}}^{-5/2})$}]
				\addplot table [x=nDofsK1, y=L2P4K1, col sep=comma] {\DATAFILEd};
				\addplot table [x=nDofsK2, y=L2P4K2, col sep=comma] {\DATAFILEd};
				\addplot table [x=nDofsK4, y=L2P4K4, col sep=comma] {\DATAFILEd};
				\addplot table [x=nDofsK8, y=L2P4K8, col sep=comma] {\DATAFILEd};
				\addplot table [x=nDofsK16, y=L2P4K16, col sep=comma] {\DATAFILEd};
				\addplot table [x=nDofsK1, y=L2p4k1, col sep=comma] {\DATAFILEdi};
				\addplot table [x=nDofsK2, y=L2p4k2, col sep=comma] {\DATAFILEdi};
				\addplot table [x=nDofsK4, y=L2p4k4, col sep=comma] {\DATAFILEdi};
				\addplot table [x=nDofsK8, y=L2p4k8, col sep=comma] {\DATAFILEdi};
				\addplot table [x=nDofsK16, y=L2p4k16, col sep=comma] {\DATAFILEdi};
				\addplot[dashed,black,thick] table [x=nDofsK1, y=h5L2, col sep=comma] {\DATAFILEd};
			\end{loglogaxis}
		\end{tikzpicture}
	\end{subfigure}\hfill
	\begin{subfigure}[t]{0.425\linewidth}
		\centering
		\begin{tikzpicture}[%font=\footnotesize,
			trim axis right]
			\begin{loglogaxis}[
				cycle list name=Lista1,
				width=\linewidth,
				height=\HeiFact\linewidth,
				xlabel={$N_{\mathrm{dof}}/\sharp \lambda$},
				xminorticks=false,yminorticks=false,				%ylabel={$| u - u_h |_{H^1(Q)}/| u |_{H^1(Q)}$},
				legend columns=1,
				legend pos=south east,
				legend style={at={(0.01,0.01)},anchor=south west},
				xmajorgrids=true, ymajorgrids=true,
				ymin={2.5e-13}, ymax=10, 
				xmin=0.05,
				xmax=20000,
				legend entries={,,,,,,,,,,$O(N_{\mathrm{dof}}^{-2})$}]
				\addplot table [x=nDofsK1, y=H1sP4K1, col sep=comma] {\DATAFILEd};
				\addplot table [x=nDofsK2, y=H1sP4K2, col sep=comma] {\DATAFILEd};
				\addplot table [x=nDofsK4, y=H1sP4K4, col sep=comma] {\DATAFILEd};
				\addplot table [x=nDofsK8, y=H1sP4K8, col sep=comma] {\DATAFILEd};
				\addplot table [x=nDofsK16, y=H1sP4K16, col sep=comma] {\DATAFILEd};
				\addplot table [x=nDofsK1, y=H1p4k1, col sep=comma] {\DATAFILEdi};
				\addplot table [x=nDofsK2, y=H1p4k2, col sep=comma] {\DATAFILEdi};
				\addplot table [x=nDofsK4, y=H1p4k4, col sep=comma] {\DATAFILEdi};
				\addplot table [x=nDofsK8, y=H1p4k8, col sep=comma] {\DATAFILEdi};
				\addplot table [x=nDofsK16, y=H1p4k16, col sep=comma] {\DATAFILEdi};
				\addplot[dashed,black,thick] table [x=nDofsK1, y=h4H1s, col sep=comma] {\DATAFILEd};
			\end{loglogaxis}
		\end{tikzpicture}
	\end{subfigure}
	\hspace*{\fill}
	\caption{Comparison between relative errors of the \igastab\ method \eqref{equivstabiga} (continuous lines) and relative errors of the \femstab\ method \eqref{stabZank} (dashed lines) plotted against the number of space DOFs per wavelength $N_{\mathrm{dof}}/\sharp \lambda$, at different wavenumbers $k$.
		$L^2(Q)$ norms are shown on the left, $H^1(Q)$ seminorms on the right.
		Rows 1 to 4 correspond to $p=1$ to $p=4$.
		%vs temporal mesh-size $h_t = h_s$, for $p=1$. 
		The exact solution is defined in \eqref{u2_ex}.
	}
	\label{fig:Ldof1234}
\end{figure}

%% file: Graphics/ErrorScattering.tex
\tikzstyle{Linea1}=[thick,dashed]
\tikzstyle{Linea2}=[thick]

%\pgfplotsset{
	%    tick label style={font=\footnotesize},
	%    label style={font=\footnotesize},
	%    legend style={font=\footnotesize},
	%}

\pgfplotscreateplotcyclelist{Lista1}{%
%	{Linea2,Red,mark=*},
%	{Linea1,Red},
%	{Linea2,Green,mark=triangle*},
%	{Linea1,Green},
%	{Linea2,Cyan,mark=square*},
%	{Linea1,Cyan},
%	{Linea2,Violet,mark=diamond*},
%	{Linea1,Violet}
%	{Linea2,Red,mark=*},
%	{Linea1,mark=*},
%	{Linea2,Green,mark=triangle*},
%	{Linea1,mark=triangle*},
%	{Linea2,Cyan,mark=square*},
%	{Linea1,mark=square*},
%	{Linea2,Violet,mark=diamond*},
%	{Linea1,mark=diamonds*}
	{Linea2,Red,mark=*},
	{Linea1,Red},
	{Linea2,Green,mark=triangle*},
	{Linea1,Green},
	{Linea2,Cyan,mark=square*},
	{Linea1,Cyan},
	{Linea2,Violet,mark=diamond*},
	{Linea1,Violet}
}

\def \DATAFILE {Graphics/ErrorTableScatteringTer.csv}

\begin{figure}[htbp]
	\centering
	\hspace*{\fill}
\begin{subfigure}[t]{0.425\linewidth}
	\centering
\begin{tikzpicture}[%font=\footnotesize, 
	trim axis left
	]
\begin{loglogaxis}[
	cycle list name=Lista1,
	width=\linewidth,
	height=\linewidth,
	xlabel={$h_t$},
	ymin=5e-6,
	ymax=10,
	xminorticks=false,
	yminorticks=false,
	%max space between ticks=20,
	%xtick=data,
	%xticklabel={\pgfmathparse{exp(\tick)}\pgfmathprintnumber[fixed relative,precision=3]{\pgfmathresult}},
	%ytick={0.4,0.6,0.8,1}, 
	%yticklabels={$0.4$,$0.6$,$0.8$,$1$}, 
	ylabel={$\| u - u_h \|_{L^2(Q)} / \| u \|_{L^2(Q)}$
	},
	legend columns=1,
	%legend pos=south east,
	legend style={at={(0.99,0.01)},anchor=south east},
	xmajorgrids=true,
	ymajorgrids=true,
	legend entries={$p=1$,$O(h^2)$,$p=2$,$O(h^3)$,$p=3$,$O(h^4)$,$p=4$,$O(h^5)$}]
	\addplot table [x=ht, y=L2P1, col sep=comma] {\DATAFILE};
	\addplot table [x=ht, y=h2L2, col sep=comma] {\DATAFILE};
	\addplot table [x=ht, y=L2P2, col sep=comma] {\DATAFILE};
	\addplot table [x=ht, y=h3L2, col sep=comma] {\DATAFILE};
	\addplot table [x=ht, y=L2P3, col sep=comma] {\DATAFILE};
	\addplot table [x=ht, y=h4L2, col sep=comma] {\DATAFILE};
	\addplot table [x=ht, y=L2P4, col sep=comma] {\DATAFILE};
	\addplot table [x=ht, y=h5L2, col sep=comma] {\DATAFILE};
\end{loglogaxis}
\end{tikzpicture}
\end{subfigure}
\hfill
\begin{subfigure}[t]{0.425\linewidth}
	\centering
	\begin{tikzpicture}[%font=\footnotesize,
		trim axis right
		]
	\begin{loglogaxis}[
	cycle list name=Lista1,
	width=\linewidth,
	height=\linewidth,
	xlabel={$h_t$},
	ymin=5e-6,
	ymax=10,
	xminorticks=false,
	yminorticks=false,
	%max space between ticks=20,
	%xtick=data,
	%xticklabel={\pgfmathparse{exp(\tick)}\pgfmathprintnumber[fixed relative,precision=3]{\pgfmathresult}},
	%ytick={0.4,0.6,0.8,1}, 
	%yticklabels={$0.4$,$0.6$,$0.8$,$1$}, 
	ylabel={$| u - u_h |_{H^1(Q)} /| u |_{H^1(Q)}$
	},
	legend columns=1,
	%legend pos=south east,
	legend style={at={(0.99,0.01)},anchor=south east},
	xmajorgrids=true,
	ymajorgrids=true,
	legend entries={$p=1$,$O(h)$,$p=2$,$O(h^2)$,$p=3$,$O(h^3)$,$p=4$,$O(h^4)$}]
	\addplot table [x=ht, y=H1sP1, col sep=comma] {\DATAFILE};
	\addplot table [x=ht, y=h1H1s, col sep=comma] {\DATAFILE};
	\addplot table [x=ht, y=H1sP2, col sep=comma] {\DATAFILE};
	\addplot table [x=ht, y=h2H1s, col sep=comma] {\DATAFILE};
	\addplot table [x=ht, y=H1sP3, col sep=comma] {\DATAFILE};
	\addplot table [x=ht, y=h3H1s, col sep=comma] {\DATAFILE};
	\addplot table [x=ht, y=H1sP4, col sep=comma] {\DATAFILE};
	\addplot table [x=ht, y=h4H1s, col sep=comma] {\DATAFILE};
\end{loglogaxis}
\end{tikzpicture}
\end{subfigure}
\hspace*{\fill}
\caption{Relative errors of the \igastab\ method \eqref{equivstabiga} plotted against the time mesh-size $h_t \approx h_s$ for the scattering problem \ref{ex:scattering}.
}
\label{fig:ErrorScattering}
\end{figure}

%% file: Graphics/ErrorCsmoothRegEx.tex
\tikzstyle{Linea1}=[thick,dashed]
\tikzstyle{Linea2}=[thick]

%\pgfplotsset{
	%    tick label style={font=\footnotesize},
	%    label style={font=\footnotesize},
	%    legend style={font=\footnotesize},
	%}

\pgfplotscreateplotcyclelist{Lista1}{%
	%	{Linea2,Red,mark=*},
	%	{Linea1,Red},
	%	{Linea2,Green,mark=triangle*},
	%	{Linea1,Green},
	%	{Linea2,Cyan,mark=square*},
	%	{Linea1,Cyan},
	%	{Linea2,Violet,mark=diamond*},
	%	{Linea1,Violet}
	%	{Linea2,Red,mark=*},
	%	{Linea1,mark=*},
	%	{Linea2,Green,mark=triangle*},
	%	{Linea1,mark=triangle*},
	%	{Linea2,Cyan,mark=square*},
	%	{Linea1,mark=square*},
	%	{Linea2,Violet,mark=diamond*},
	%	{Linea1,mark=diamonds*}
	{Linea2,Red,mark=*},
	{Linea1,Red},
	{Linea2,Green,mark=triangle*},
	{Linea1,Green},
	{Linea2,Cyan,mark=square*},
	{Linea1,Cyan},
	{Linea2,Violet,mark=diamond*},
	{Linea1,Violet}
}

\def \DATAFILE {Graphics/ErrorTableCsmoothEx.csv}

\begin{figure}[htbp]
	\centering
	\hspace*{\fill}
	\begin{subfigure}[t]{0.425\linewidth}
		\centering
		\begin{tikzpicture}[%font=\footnotesize, 
			trim axis left
			]
			\begin{loglogaxis}[
				cycle list name=Lista1,
				width=\linewidth,
				height=\linewidth,
				xlabel={$h_t$},
				ymin=1e-9,
				ymax=5,
				xminorticks=false,
				yminorticks=false,
				%max space between ticks=20,
				%xtick=data,
				%xticklabel={\pgfmathparse{exp(\tick)}\pgfmathprintnumber[fixed relative,precision=3]{\pgfmathresult}},
				%ytick={0.4,0.6,0.8,1}, 
				%yticklabels={$0.4$,$0.6$,$0.8$,$1$}, 
				ylabel={$\| u - u_h \|_{L^2(Q)} / \| u \|_{L^2(Q)}$
				},
				legend columns=1,
				%legend pos=south east,
				legend style={at={(0.99,0.01)},anchor=south east},
				xmajorgrids=true,
				ymajorgrids=true,
				legend entries={$p=1$,$O(h^2)$,$p=2$,$O(h^3)$,$p=3$,$O(h^4)$,$p=4$,$O(h^5)$}]
				\addplot table [x=ht, y=L2P1, col sep=comma] {\DATAFILE};
				\addplot table [x=ht, y=h2L2, col sep=comma] {\DATAFILE};
				\addplot table [x=ht, y=L2P2, col sep=comma] {\DATAFILE};
				\addplot table [x=ht, y=h3L2, col sep=comma] {\DATAFILE};
				\addplot table [x=ht, y=L2P3, col sep=comma] {\DATAFILE};
				\addplot table [x=ht, y=h4L2, col sep=comma] {\DATAFILE};
				\addplot table [x=ht, y=L2P4, col sep=comma] {\DATAFILE};
				\addplot table [x=ht, y=h5L2, col sep=comma] {\DATAFILE};
			\end{loglogaxis}
		\end{tikzpicture}
	\end{subfigure}
	\hfill
	\begin{subfigure}[t]{0.425\linewidth}
		\centering
		\begin{tikzpicture}[%font=\footnotesize,
			trim axis right
			]
			\begin{loglogaxis}[
				cycle list name=Lista1,
				width=\linewidth,
				height=\linewidth,
				xlabel={$h_t$},
				ymin=1e-9,
				ymax=5,
				xminorticks=false,
				yminorticks=false,
				%max space between ticks=20,
				%xtick=data,
				%xticklabel={\pgfmathparse{exp(\tick)}\pgfmathprintnumber[fixed relative,precision=3]{\pgfmathresult}},
				%ytick={0.4,0.6,0.8,1}, 
				%yticklabels={$0.4$,$0.6$,$0.8$,$1$}, 
				ylabel={$| u - u_h |_{c,H^1(Q)} /| u |_{c,H^1(Q)}$
				},
				legend columns=1,
				%legend pos=south east,
				legend style={at={(0.99,0.01)},anchor=south east},
				xmajorgrids=true,
				ymajorgrids=true,
				legend entries={$p=1$,$O(h)$,$p=2$,$O(h^2)$,$p=3$,$O(h^3)$,$p=4$,$O(h^4)$}]
				\addplot table [x=ht, y=H1sP1, col sep=comma] {\DATAFILE};
				\addplot table [x=ht, y=h1H1s, col sep=comma] {\DATAFILE};
				\addplot table [x=ht, y=H1sP2, col sep=comma] {\DATAFILE};
				\addplot table [x=ht, y=h2H1s, col sep=comma] {\DATAFILE};
				\addplot table [x=ht, y=H1sP3, col sep=comma] {\DATAFILE};
				\addplot table [x=ht, y=h3H1s, col sep=comma] {\DATAFILE};
				\addplot table [x=ht, y=H1sP4, col sep=comma] {\DATAFILE};
				\addplot table [x=ht, y=h4H1s, col sep=comma] {\DATAFILE};
			\end{loglogaxis}
		\end{tikzpicture}
	\end{subfigure}
	\hspace*{\fill}
	\caption{Relative errors of the \igastab\ method \eqref{equivstabiga} plotted against the time mesh-size $h_t \approx h_s$ for the problem with a smooth velocity field $c$ described in Section \ref{sec:Csmooth}. The exact solution is defined in \eqref{ucsmooth_ex}.
	}
	\label{fig:ErrorCsmoothRegEx}
\end{figure}

%% file: Graphics/ErrorCdiscBumpConvWide2in1.tex
\tikzstyle{Linea1}=[thick,dashed]
\tikzstyle{Linea2}=[thick]%,mark=*]

%\pgfplotsset{
	%    tick label style={font=\footnotesize},
	%    label style={font=\footnotesize},
	%    legend style={font=\footnotesize},
	%}

\pgfplotscreateplotcyclelist{Lista1}{%
	{Linea2,Red,mark=*},
	{Linea1,Red},
	{Linea2,Green,mark=triangle*},
	%{Linea1,Green},
	{Linea2,Cyan,mark=square*},
	%{Linea1,Cyan},
	{Linea2,Violet,mark=diamond*},
	{Linea1,Violet},
}
\pgfplotscreateplotcyclelist{Lista2}{%
	{Linea2,Red,mark=*},
	{Linea1,Red},
	{Linea2,Green,mark=triangle*},
	%{Linea1,Green},
	{Linea2,Cyan,mark=square*},
	%{Linea1,Cyan},
	{Linea2,Violet,mark=diamond*},
	%{Linea1,Violet},
}

\pgfplotscreateplotcyclelist{Lista1rs}{%
	{Linea2,Red,mark=*},
	{Linea1,Red},
	{Linea2,Green,mark=triangle*},
	{Linea1,Green},
	{Linea2,Cyan,mark=square*},
	{Linea1,Cyan},
	{Linea2,Violet,mark=diamond*},
	{Linea1,Violet},
}

\def \DATAFILE {Graphics/ErrorTableCdiscBumpConvWide.csv}

\def \DATAFILErs {Graphics/ErrorTableCdiscBumpConvWideRS0.csv}

\begin{figure}[htbp]
	\centering
	\hspace*{\fill}
	\begin{subfigure}[t]{0.425\linewidth}
		\centering
		\begin{tikzpicture}[%font=\footnotesize, 
			trim axis left
			]
			\begin{loglogaxis}[
				yticklabel style={text width=2.5em, align=right},
				cycle list name=Lista1,
				width=\linewidth,
				height=\linewidth,
				ymin=25e-6,
				ymax=5e-1,
				xminorticks=false,
				yminorticks=false,
				%max space between ticks=20,
				%xtick=data,
				%xticklabel={\pgfmathparse{exp(\tick)}\pgfmathprintnumber[fixed relative,precision=3]{\pgfmathresult}},
				%ytick={0.4,0.6,0.8,1}, 
				%yticklabels={$0.4$,$0.6$,$0.8$,$1$}, 
				ylabel={$\| u - u_h \|_{L^2(Q)}/\| u \|_{L^2(Q)}$},
				legend columns=1,
				%legend pos=south east,
				legend style={at={(0.99,0.01)},anchor=south east},
				xmajorgrids=true,
				ymajorgrids=true,
				legend entries={$p=1$,$O(h^2)$,$p=2$,$p=3$,$p=4$}]
				\addplot table [x=ht, y=L2P1, col sep=comma] {\DATAFILE};
				\addplot table [x=ht, y=h2L2, col sep=comma] {\DATAFILE};
				\addplot table [x=ht, y=L2P2, col sep=comma] {\DATAFILE};
				%\addplot table [x=ht, y=h3L2, col sep=comma] {\DATAFILE};
				\addplot table [x=ht, y=L2P3, col sep=comma] {\DATAFILE};
				%\addplot table [x=ht, y=h4L2, col sep=comma] {\DATAFILE};
				\addplot table [x=ht, y=L2P4, col sep=comma] {\DATAFILE};
				%\addplot table [x=ht, y=h5L2, col sep=comma] {\DATAFILE};
			\end{loglogaxis}
		\end{tikzpicture}
	\end{subfigure}
	\hfill
	\begin{subfigure}[t]{0.425\linewidth}
		\centering
		\begin{tikzpicture}[%font=\footnotesize,
			trim axis right
			]
			\begin{loglogaxis}[
				yticklabel style={text width=2.5em, align=right},
				cycle list name=Lista2,
				width=\linewidth,
				height=\linewidth,
				ymin=25e-6,
				ymax=5e-1,
				xminorticks=false,
				yminorticks=false,
				%max space between ticks=20,
				%xtick=data,
				%xticklabel={\pgfmathparse{exp(\tick)}\pgfmathprintnumber[fixed relative,precision=3]{\pgfmathresult}},
				%ytick={0.4,0.6,0.8,1}, 
				%yticklabels={$0.4$,$0.6$,$0.8$,$1$}, 
				ylabel={$| u - u_h |_{c,H^1(Q)}/| u |_{c,H^1(Q)}$},
				legend columns=1,
				%legend pos=south east,
				legend style={at={(0.99,0.01)},anchor=south east},
				xmajorgrids=true,
				ymajorgrids=true,
				legend entries={$p=1$,$O(h)$,$p=2$,$p=3$,$p=4$}]
				\addplot table [x=ht, y=H1sP1, col sep=comma] {\DATAFILE};
				\addplot table [x=ht, y=h1H1s, col sep=comma] {\DATAFILE};
				\addplot table [x=ht, y=H1sP2, col sep=comma] {\DATAFILE};
				%\addplot table [x=ht, y=ht2H1sp2, col sep=comma] {\DATAFILE};
				\addplot table [x=ht, y=H1sP3, col sep=comma] {\DATAFILE};
				%\addplot table [x=ht, y=ht2H1sp3, col sep=comma] {\DATAFILE};
				\addplot table [x=ht, y=H1sP4, col sep=comma] {\DATAFILE};
				%\addplot table [x=ht, y=ht2H1sp4, col sep=comma] {\DATAFILE};
			\end{loglogaxis}
		\end{tikzpicture}
	\end{subfigure}
	\hspace*{\fill}
	\\
	\hspace*{\fill}
	\begin{subfigure}[t]{0.425\linewidth}
		\centering
		\begin{tikzpicture}[%font=\footnotesize, 
			trim axis left
			]
			\begin{loglogaxis}[
				yticklabel style={text width=2.5em, align=right},
				cycle list name=Lista1rs,
				width=\linewidth,
				height=\linewidth,
				xlabel={$h_t$},
				ymin=5e-11,
				ymax=5e-1,
				xminorticks=false,
				yminorticks=false,
				%max space between ticks=20,
				%xtick=data,
				%xticklabel={\pgfmathparse{exp(\tick)}\pgfmathprintnumber[fixed relative,precision=3]{\pgfmathresult}},
				%ytick={0.4,0.6,0.8,1}, 
				%yticklabels={$0.4$,$0.6$,$0.8$,$1$}, 
				ylabel={$\| u - u_h \|_{L^2(Q)}/\| u \|_{L^2(Q)}$},
				legend columns=1,
				%legend pos=south east,
				legend style={at={(0.99,0.01)},anchor=south east},
				xmajorgrids=true,
				ymajorgrids=true,
				legend entries={$p=1$,$O(h^2)$,$p=2$,$O(h^3)$,$p=3$,$O(h^4)$,$p=4$,$O(h^5)$}]
				\addplot table [x=ht, y=L2P1, col sep=comma] {\DATAFILErs};
				\addplot table [x=ht, y=h2L2, col sep=comma] {\DATAFILErs};
				\addplot table [x=ht, y=L2P2, col sep=comma] {\DATAFILErs};
				\addplot table [x=ht, y=h3L2, col sep=comma] {\DATAFILErs};
				\addplot table [x=ht, y=L2P3, col sep=comma] {\DATAFILErs};
				\addplot table [x=ht, y=h4L2, col sep=comma] {\DATAFILErs};
				\addplot table [x=ht, y=L2P4, col sep=comma] {\DATAFILErs};
				\addplot table [x=ht, y=h5L2, col sep=comma] {\DATAFILErs};
			\end{loglogaxis}
		\end{tikzpicture}
	\end{subfigure}
	\hfill
	\begin{subfigure}[t]{0.425\linewidth}
		\centering
		\begin{tikzpicture}[%font=\footnotesize,
			trim axis right
			]
			\begin{loglogaxis}[
				yticklabel style={text width=2.5em, align=right},
				cycle list name=Lista1rs,
				width=\linewidth,
				height=\linewidth,
				xlabel={$h_t$},
				ymin=5e-11,
				ymax=5e-1,
				xminorticks=false,
				yminorticks=false,
				%max space between ticks=20,
				%xtick=data,
				%xticklabel={\pgfmathparse{exp(\tick)}\pgfmathprintnumber[fixed relative,precision=3]{\pgfmathresult}},
				%ytick={0.4,0.6,0.8,1}, 
				%yticklabels={$0.4$,$0.6$,$0.8$,$1$}, 
				ylabel={$| u - u_h |_{c,H^1(Q)}/| u |_{c,H^1(Q)}$},
				legend columns=1,
				%legend pos=south east,
				legend style={at={(0.99,0.01)},anchor=south east},
				xmajorgrids=true,
				ymajorgrids=true,
				legend entries={$p=1$,$O(h)$,$p=2$,$O(h^2)$,$p=3$,$O(h^3)$,$p=4$,$O(h^4)$}]
				\addplot table [x=ht, y=H1sP1, col sep=comma] {\DATAFILErs};
				\addplot table [x=ht, y=h1H1s, col sep=comma] {\DATAFILErs};
				\addplot table [x=ht, y=H1sP2, col sep=comma] {\DATAFILErs};
				\addplot table [x=ht, y=h2H1s, col sep=comma] {\DATAFILErs};
				\addplot table [x=ht, y=H1sP3, col sep=comma] {\DATAFILErs};
				\addplot table [x=ht, y=h3H1s, col sep=comma] {\DATAFILErs};
				\addplot table [x=ht, y=H1sP4, col sep=comma] {\DATAFILErs};
				\addplot table [x=ht, y=h4H1s, col sep=comma] {\DATAFILErs};
			\end{loglogaxis}
		\end{tikzpicture}
	\end{subfigure}
	\hspace*{\fill}
	\caption{Relative errors of the \igastab\ method~\eqref{equivstabiga}, plotted against the temporal mesh-size $h_t = h_s$, for the problem with piecewise-constant velocity described in Section \ref{sec:Cdisc}.
	The relative errors are plotted against the temporal mesh-size $h_t = h_s$.
	\newline
	Top row: splines with maximal regularity both in space and time.
	\newline
	Bottom row: splines of space--time maximal regularity except at $x=0.5$, where only $C^0$-continuity is imposed.}
	\label{fig:CdiscBumpStabConv2in1}
\end{figure}

%% file: Graphics/StabEnergy.tex
\tikzstyle{Linea1}=[thick,dashed]
\tikzstyle{Linea2}=[thick, only marks, mark=oplus]
\tikzstyle{Linea3}=[thick, only marks, mark=halfcircle]

%\pgfplotsset{
	%    tick label style={font=\footnotesize},
	%    label style={font=\footnotesize},
	%    legend style={font=\footnotesize},
	%}

\pgfplotscreateplotcyclelist{Lista1}{%
	{Linea1,Red},
	{Linea1,Green},
	{Linea1,Cyan},
	{Linea1,Violet},
	%{Linea1,BurntOrange},
	%{Linea1,Brown},
	%{Linea1,Magenta},
	%{Linea1,Blue},
	%{Linea2,Goldenrod},
	{Linea2,Red},
	{Linea2,Green},
	{Linea2,Cyan},
	{Linea2,Violet},
	%{Linea2,BurntOrange},
	%{Linea2,Brown},
	%{Linea2,Magenta},
	%{Linea2,Blue},
	%{Linea2,Goldenrod},
	{Linea3,Red},
	{Linea3,Green},
	{Linea3,Cyan},
	{Linea3,Violet},
	%{Linea3,BurntOrange},
	%{Linea3,Brown},
	%{Linea3,Magenta},
	%{Linea3,Blue},
	%{Linea3,Goldenrod},
}

\def \DATAFILE {Graphics/ErrorEnergyTable.csv}

\begin{figure}[htbp]
\centering %\hfill
\begin{tikzpicture}%[font=\footnotesize, 	%trim axis left	]
\begin{semilogyaxis}[scale=1.1,
	cycle list name=Lista1,
	width=0.8\linewidth,
	height=0.4\linewidth,
	xlabel={$t$},
	xmin=0,
	xmax=10,
	ymin=0.0000000000001,
	xminorticks=false,
	yminorticks=false,
	%max space between ticks=20,
	%xtick=data,
	%xticklabel={\pgfmathparse{exp(\tick)}\pgfmathprintnumber[fixed relative,precision=3]{\pgfmathresult}},
	ytick={0.0000000001,0.00000001,0.000001,0.0001,0.01}, 
	ylabel={$|E_h(t)-E(t)|/E(t)$},
	legend columns=4,
	legend pos=south east,
	legend style={at={(0.99,0.01)},anchor=south east},
	grid=both,
	legend entries={$p=1$,$p=2$,$p=3$,$p=4$}]
	\addplot table [x=t, y=p1, col sep=comma] {\DATAFILE};
	\addplot table [x=t, y=p2, col sep=comma] {\DATAFILE};
	\addplot table [x=t, y=p3, col sep=comma] {\DATAFILE};
	\addplot table [x=t, y=p4, col sep=comma] {\DATAFILE};
	\addplot table [x=t, y=p1, col sep=comma,restrict expr to domain={\thisrow{p1Sign}}{1:1}] {\DATAFILE};
	\addplot table [x=t, y=p2, col sep=comma,restrict expr to domain={\thisrow{p2Sign}}{1:1}] {\DATAFILE};
	\addplot table [x=t, y=p3, col sep=comma,restrict expr to domain={\thisrow{p3Sign}}{1:1}] {\DATAFILE};
	\addplot table [x=t, y=p4, col sep=comma,restrict expr to domain={\thisrow{p4Sign}}{1:1}] {\DATAFILE};
	\addplot table [x=t, y=p1, col sep=comma,restrict expr to domain={\thisrow{p1Sign}}{-1:-1}] {\DATAFILE};
	\addplot table [x=t, y=p2, col sep=comma,restrict expr to domain={\thisrow{p2Sign}}{-1:-1}] {\DATAFILE};
	\addplot table [x=t, y=p3, col sep=comma,restrict expr to domain={\thisrow{p3Sign}}{-1:-1}] {\DATAFILE};
	\addplot table [x=t, y=p4, col sep=comma,restrict expr to domain={\thisrow{p4Sign}}{-1:-1}] {\DATAFILE};
\end{semilogyaxis}
\end{tikzpicture}
\caption{Time evolution of the energy relative error for the problem with solution \eqref{uex_energy}.
%Relative errors of the energy
The marker ``$\oplus$'' denotes time instants when $E_h \geq E$, while ``$\ominus$'' stands for $E_h \leq E$.}
\label{fig:energy}
\end{figure}
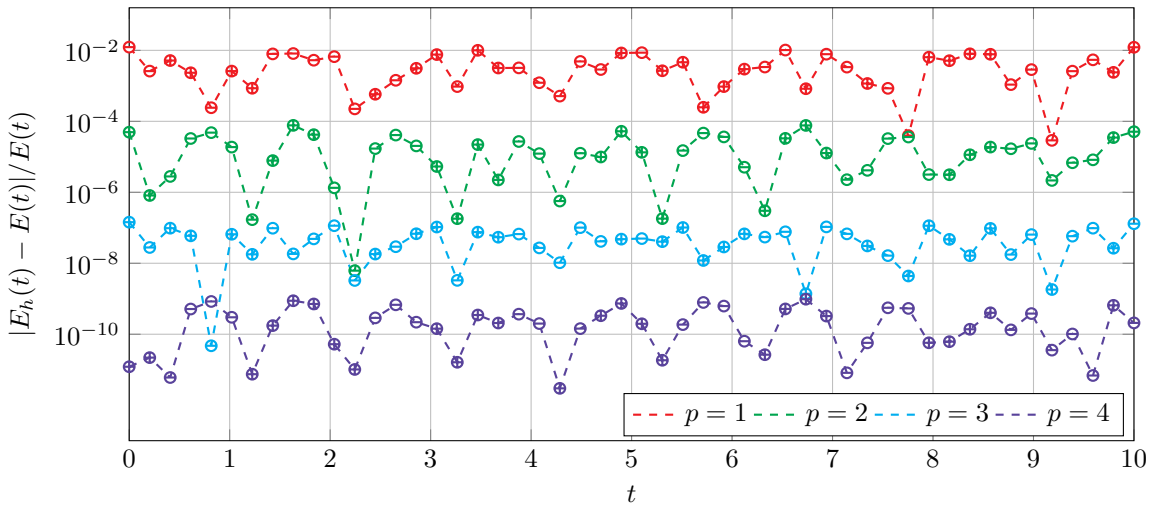

%% file: Graphics/Dispersion_p1234Tent.tex
\tikzstyle{Linea1}=[thick]
\tikzstyle{Linea2}=[thick,mark=*]

%\pgfplotsset{
	%    tick label style={font=\footnotesize},
	%    label style={font=\footnotesize},
	%    legend style={font=\footnotesize},
	%}

\pgfplotscreateplotcyclelist{Lista1}{%
	{Linea1,Black},
	{Linea1,Red},
	{Linea1,Green},
	{Linea1,Cyan},
	{Linea1,Violet},
	{Linea1,BurntOrange},
	{Linea1,Brown},
	{Linea1,Magenta},
	{Linea1,Blue},
}
\pgfplotscreateplotcyclelist{Lista2}{%{Linea1,Red},
	{Linea1,Red},
	{Linea1,Green},
	{Linea1,Cyan},
	{Linea1,Violet},
	{Linea1,BurntOrange},
	{Linea1,Brown},
	{Linea1,Magenta},
	{Linea1,Blue},
}

\def\DATAFILEpAo{Graphics/PlotTableTentOur_p1.csv}
\def\DATAFILEpBo{Graphics/PlotTableTentOur_p2.csv}
\def\DATAFILEpCo{Graphics/PlotTableTentOur_p3.csv}
\def\DATAFILEpDo{Graphics/PlotTableTentOur_p4.csv}

\def\DATAFILEpAz{Graphics/PlotTableTentZank_p1.csv}
\def\DATAFILEpBz{Graphics/PlotTableTentZank_p2.csv}
\def\DATAFILEpCz{Graphics/PlotTableTentZank_p3.csv}
\def\DATAFILEpDz{Graphics/PlotTableTentZank_p4.csv}

\def\HeiFac{0.3}

\begin{figure}[htbp]
\centering
	\hspace*{\fill}
	\begin{subfigure}[t]{0.425\linewidth}
		\centering
		\begin{tikzpicture}[trim axis left]
			\begin{axis}[
				title={\igastab},
				yticklabel style={text width=2.5em, align=right},
				cycle list name=Lista1,
				width=\linewidth,
				height=\linewidth,
				xlabel={$x$},
				xmin=0,
				xmax=1,
				ymax=1.05,
				ymin=-0.2,
				xminorticks=false, yminorticks=false,
				legend style={at={(0.99,0.99)},anchor=north east},
				grid=both,
				legend entries={$u(x,0)$\\ $u_h(x,T),p=1$\\ $u_h(x,T),p=2$\\ $u_h(x,T),p=3$\\ $u_h(x,T),p=4$\\}
				]
				\addplot table [x=phyX, y=uh_start, col sep=comma] {\DATAFILEpAo};
				\addplot table [x=phyX, y=uh_end, col sep=comma] {\DATAFILEpAo};
				\addplot table [x=phyX, y=uh_end, col sep=comma] {\DATAFILEpBo};
				\addplot table [x=phyX, y=uh_end, col sep=comma] {\DATAFILEpCo};
				\addplot table [x=phyX, y=uh_end, col sep=comma] {\DATAFILEpDo};
			\end{axis}
		\end{tikzpicture}
	\end{subfigure}\hfill
	\begin{subfigure}[t]{0.425\linewidth}
		\centering
		\begin{tikzpicture}[trim axis right]
			\begin{axis}[
				title={\femstab},
				yticklabel style={text width=2.5em, align=right},
				cycle list name=Lista1,
				width=\linewidth,
				height=\linewidth,
				xlabel={$x$},
				xmin=0,
				xmax=1,
				ymax=1.05,
				ymin=-0.2,
				xminorticks=false, yminorticks=false,
				legend style={at={(0.99,0.99)},anchor=north east},
				grid=both,
				legend entries={$u(x,0)$\\ $u_h(x,T),p=1$\\ $u_h(x,T),p=2$\\ $u_h(x,T),p=3$\\ $u_h(x,T),p=4$\\}
				]
				\addplot table [x=phyX, y=uh_start, col sep=comma] {\DATAFILEpAz};
				\addplot table [x=phyX, y=uh_end, col sep=comma] {\DATAFILEpAz};
				\addplot table [x=phyX, y=uh_end, col sep=comma] {\DATAFILEpBz};
				\addplot table [x=phyX, y=uh_end, col sep=comma] {\DATAFILEpCz};
				\addplot table [x=phyX, y=uh_end, col sep=comma] {\DATAFILEpDz};
			\end{axis}
		\end{tikzpicture}
	\end{subfigure}
\hspace*{\fill}
\\
\bigskip\bigskip
%\noindent
\hspace*{\fill}
	\begin{subfigure}[t]{0.425\linewidth}
		\centering
		\begin{tikzpicture}[trim axis left]
			\begin{axis}[
				yticklabel style={text width=2.5em, align=right},
				cycle list name=Lista2,
				width=\linewidth,
				height=\linewidth,
				xlabel={$x$},
				xmin=0,
				xmax=1,
				ytick={-0.1,0,0.1}, 
				xminorticks=false, yminorticks=false,
				grid=both,
				legend entries={$p=1$\\ $p=2$\\ $p=3$\\ $p=4$\\}
				]
				%\addplot table [x=phyX, y=eu_start, col sep=comma] {\DATAFILEpA};
				\addplot table [x=phyX, y expr=\thisrow{uh_end} - \thisrow{uh_start}, col sep=comma] {\DATAFILEpAo};
				\addplot table [x=phyX, y expr=\thisrow{uh_end} - \thisrow{uh_start}, col sep=comma] {\DATAFILEpBo};
				\addplot table [x=phyX, y expr=\thisrow{uh_end} - \thisrow{uh_start}, col sep=comma] {\DATAFILEpCo};
				\addplot table [x=phyX, y expr=\thisrow{uh_end} - \thisrow{uh_start}, col sep=comma] {\DATAFILEpDo};
			\end{axis}
		\end{tikzpicture}
	\end{subfigure}
	\hfill
	\begin{subfigure}[t]{0.425\linewidth}
		\centering
		\begin{tikzpicture}[trim axis right]
			\begin{axis}[
				yticklabel style={text width=2.5em, align=right},
				cycle list name=Lista2,
				width=\linewidth,
				height=\linewidth,
				xlabel={$x$},
				xmin=0,
				xmax=1,
				ytick={-0.1,0,0.1}, 
				xminorticks=false, yminorticks=false,
				grid=both,
				legend entries={$p=1$\\ $p=2$\\ $p=3$\\ $p=4$\\}
				]
				%\addplot table [x=phyX, y=eu_start, col sep=comma] {\DATAFILEpA};
				\addplot table [x=phyX, y expr=\thisrow{uh_end} - \thisrow{uh_start}, col sep=comma] {\DATAFILEpAz};
				\addplot table [x=phyX, y expr=\thisrow{uh_end} - \thisrow{uh_start}, col sep=comma] {\DATAFILEpBz};
				\addplot table [x=phyX, y expr=\thisrow{uh_end} - \thisrow{uh_start}, col sep=comma] {\DATAFILEpCz};
				\addplot table [x=phyX, y expr=\thisrow{uh_end} - \thisrow{uh_start}, col sep=comma] {\DATAFILEpDz};
			\end{axis}
		\end{tikzpicture}
	\end{subfigure}
	\hspace*{\fill}
	\caption{Approximation of the periodic problem of Section~\ref{s:Dispersion} with initial data \eqref{eq:tentInitial}.
	Top row: exact and discrete solutions at final time $T=2$; bottom row: error at final time.
	Left column: \igastab\ method; right column: \femstab\ method.
	For both the stabilizations and all the spline degrees, $N_{\mathrm{dof}}=17\,424$ and $h_t \approx 2 h_s$.}
	\label{fig:dispersionTent}
\end{figure}

%% file: Graphics/DispersionTestErrorsTent.tex
\tikzstyle{Linea1}=[thick,dashed, mark options={solid}]
\tikzstyle{Linea2}=[thick]%,mark=*]
\tikzstyle{Linea3}=[thick,dashed]

%\pgfplotsset{
	%    tick label style={font=\footnotesize},
	%    label style={font=\footnotesize},
	%    legend style={font=\footnotesize},
	%}

\pgfplotscreateplotcyclelist{Lista1}{%
	{Linea2,Red,mark=*},
	{Linea1,Red,mark=*},
	%{Linea3,Red},
	{Linea2,Green,mark=triangle*},
	{Linea1,Green,mark=triangle*},
	%{Linea3,Green},
	{Linea2,Cyan,mark=square*},
	{Linea1,Cyan,mark=square*},
	%{Linea3,Cyan},
	{Linea2,Violet,mark=diamond*},
	{Linea1,Violet,mark=diamond*},
	%{Linea3,Violet},
}

\def \DATAFILEo {Graphics/ErrorTableDispersionTentOur.csv}
\def \DATAFILEz {Graphics/ErrorTableDispersionTentZank.csv}

\begin{figure}[htbp]
			\centering
	\hspace*{\fill}
	\begin{subfigure}[t]{0.425\linewidth}
\centering
		\begin{tikzpicture}[%font=\footnotesize, 
		trim axis left
			]
			\begin{semilogyaxis}[
				cycle list name=Lista1,
				width=.95\linewidth,
				height=.95\linewidth,
				xlabel={$p$},
				ymin=5e-5,
				ymax=5e-1,
				xminorticks=false,
				yminorticks=false,
				%max space between ticks=20,
				%xtick=data,
				%xticklabel={\pgfmathparse{exp(\tick)}\pgfmathprintnumber[fixed relative,precision=3]{\pgfmathresult}},
				%ytick={0.4,0.6,0.8,1}, 
				%yticklabels={$0.4$,$0.6$,$0.8$,$1$}, 
				ylabel={$\| u(T) - u_h(T) \|_{L^2(\Omega)}/\| u(T) \|_{L^2(\Omega)}$},
				legend columns=1,
				%legend pos=south east,
				legend style={at={(0.99,0.99)},anchor=north east},
				xmajorgrids=true,
				ymajorgrids=true,
				xtick={1,2,3,4},
				legend entries={\igastab,\femstab}]
				\addplot table [x=p, y=ErrL2, col sep=comma] {\DATAFILEo};
				\addplot table [x=p, y=ErrL2, col sep=comma] {\DATAFILEz};
				%\addplot table [x=NdofP1, y=NdofL2P1O, col sep=comma] {\DATAFILE};
			\end{semilogyaxis}
		\end{tikzpicture}
	\end{subfigure}
	\hfill
	\begin{subfigure}[t]{0.425\linewidth}
		\centering
		\begin{tikzpicture}[%font=\footnotesize,
			trim axis right
			]
			\begin{semilogyaxis}[
				cycle list name=Lista1,
				width=.95\linewidth,
				height=.95\linewidth,
				xlabel={$p$},
				ymin=5e-5,
				ymax=5e-1,
				xminorticks=false,
				yminorticks=false,
				%max space between ticks=20,
				%xtick=data,
				%xticklabel={\pgfmathparse{exp(\tick)}\pgfmathprintnumber[fixed relative,precision=3]{\pgfmathresult}},
				%ytick={0.4,0.6,0.8,1}, 
				%yticklabels={$0.4$,$0.6$,$0.8$,$1$}, 
				ylabel={$| u(T) - u_h(T) |_{H^1(\Omega)}/| u(T) |_{H^1(\Omega)}$},
				legend columns=1,
				%legend pos=south east,
				legend style={at={(0.01,0.01)},anchor=south west},
				xmajorgrids=true,
				ymajorgrids=true,
				xtick={1,2,3,4},
				legend entries={\igastab,\femstab}]
				\addplot table [x=p, y=ErrH1s, col sep=comma] {\DATAFILEo};
				\addplot table [x=p, y=ErrH1s, col sep=comma] {\DATAFILEz};
			\end{semilogyaxis}
		\end{tikzpicture}
	\end{subfigure}
	\hspace*{\fill}
	\caption{
	Comparison between the relative errors at final time of the \igastab\ and the \femstab\ methods, for the periodic tent wave problem \eqref{eq:tentInitial}.
	For both the stabilizations and all the spline degrees, $N_{\mathrm{dof}}=17\,424$ and $h_t \approx 2 h_s$.
	}
	\label{fig:StabComparisonTent}
\end{figure}

%% file: Graphics/Dispersion_p1234Bump.tex
\tikzstyle{Linea1}=[thick]
\tikzstyle{Linea2}=[thick,mark=*]

%\pgfplotsset{
	%    tick label style={font=\footnotesize},
	%    label style={font=\footnotesize},
	%    legend style={font=\footnotesize},
	%}

\pgfplotscreateplotcyclelist{Lista1}{%
	{Linea1,Black},
	{Linea1,Red},
	{Linea1,Green},
	{Linea1,Cyan},
	{Linea1,Violet},
	{Linea1,BurntOrange},
	{Linea1,Brown},
	{Linea1,Magenta},
	{Linea1,Blue},
}
\pgfplotscreateplotcyclelist{Lista2}{%{Linea1,Red},
		{Linea1,Red},
	{Linea1,Green},
	{Linea1,Cyan},
	{Linea1,Violet},
	{Linea1,BurntOrange},
	{Linea1,Brown},
	{Linea1,Magenta},
	{Linea1,Blue},
}

\def\DATAFILEpAo{Graphics/PlotTableBumpOur_p1.csv}
\def\DATAFILEpBo{Graphics/PlotTableBumpOur_p2.csv}
\def\DATAFILEpCo{Graphics/PlotTableBumpOur_p3.csv}
\def\DATAFILEpDo{Graphics/PlotTableBumpOur_p4.csv}

\def\DATAFILEpAz{Graphics/PlotTableBumpZank_p1.csv}
\def\DATAFILEpBz{Graphics/PlotTableBumpZank_p2.csv}
\def\DATAFILEpCz{Graphics/PlotTableBumpZank_p3.csv}
\def\DATAFILEpDz{Graphics/PlotTableBumpZank_p4.csv}

\def\HeiFac{0.3}

\begin{figure}[htbp]
	\centering
	\hspace*{\fill}
	\begin{subfigure}[t]{0.425\linewidth}
		\centering
		\begin{tikzpicture}[trim axis left]
			\begin{axis}[
				title={\igastab},
				cycle list name=Lista1,
				yticklabel style={text width=2.125em, align=right},
				width=\linewidth,
				height=\linewidth,
				xlabel={$x$},
				xmin=0,
				xmax=1,
				ymax=1.1,
				ymin=-0.15,
				xminorticks=false, yminorticks=false,
				legend style={at={(0.99,0.99)},anchor=north east},
				grid=both,
				legend entries={$u(x,0)$\\ $u_h(x,T),p=1$\\ $u_h(x,T),p=2$\\ $u_h(x,T),p=3$\\ $u_h(x,T),p=4$\\}
				]
				\addplot table [x=phyX, y=uh_start, col sep=comma] {\DATAFILEpAo};
				\addplot table [x=phyX, y=uh_end, col sep=comma] {\DATAFILEpAo};
				\addplot table [x=phyX, y=uh_end, col sep=comma] {\DATAFILEpBo};
				\addplot table [x=phyX, y=uh_end, col sep=comma] {\DATAFILEpCo};
				\addplot table [x=phyX, y=uh_end, col sep=comma] {\DATAFILEpDo};
			\end{axis}
		\end{tikzpicture}
	\end{subfigure}\hfill
	\begin{subfigure}[t]{0.425\linewidth}
		\centering
		\begin{tikzpicture}[trim axis right]
			\begin{axis}[
				title={\femstab},
				cycle list name=Lista1,
				yticklabel style={text width=2.5em, align=right},
				width=\linewidth,
				height=\linewidth,
				xlabel={$x$},
				xmin=0,
				xmax=1,
				ymax=1.1,
				ymin=-0.15,
				xminorticks=false, yminorticks=false,
				legend style={at={(0.99,0.99)},anchor=north east},
				grid=both,
				legend entries={$u(x,0)$\\ $u_h(x,T),p=1$\\ $u_h(x,T),p=2$\\ $u_h(x,T),p=3$\\ $u_h(x,T),p=4$\\}
				]
				\addplot table [x=phyX, y=uh_start, col sep=comma] {\DATAFILEpAz};
				\addplot table [x=phyX, y=uh_end, col sep=comma] {\DATAFILEpAz};
				\addplot table [x=phyX, y=uh_end, col sep=comma] {\DATAFILEpBz};
				\addplot table [x=phyX, y=uh_end, col sep=comma] {\DATAFILEpCz};
				\addplot table [x=phyX, y=uh_end, col sep=comma] {\DATAFILEpDz};
			\end{axis}
		\end{tikzpicture}
	\end{subfigure}
\hspace*{\fill}
\\
\bigskip\bigskip
%\noindent
\hspace*{\fill}
	\begin{subfigure}[t]{0.425\linewidth}
		\centering
		\begin{tikzpicture}[trim axis left]
			\begin{axis}[
				cycle list name=Lista2,
				yticklabel style={text width=2.5em, align=right},
				width=\linewidth,
				height=\linewidth,
				xlabel={$x$},
				xmin=0,
				xmax=1,
				ymax=0.12,
				ymin=-0.12,
				ytick={-0.1,0,0.1}, 
				xminorticks=false, yminorticks=false,
				legend style={at={(0.8,0.01)},anchor=south east},
				grid=both,
				legend entries={$p=1$\\ $p=2$\\ $p=3$\\ $p=4$\\}
				]
				%\addplot table [x=phyX, y=eu_start, col sep=comma] {\DATAFILEpA};
				\addplot table [x=phyX, y expr=\thisrow{uh_end} - \thisrow{uh_start}, col sep=comma] {\DATAFILEpAo};
				\addplot table [x=phyX, y expr=\thisrow{uh_end} - \thisrow{uh_start}, col sep=comma] {\DATAFILEpBo};
				\addplot table [x=phyX, y expr=\thisrow{uh_end} - \thisrow{uh_start}, col sep=comma] {\DATAFILEpCo};
				\addplot table [x=phyX, y expr=\thisrow{uh_end} - \thisrow{uh_start}, col sep=comma] {\DATAFILEpDo};
			\end{axis}
		\end{tikzpicture}
	\end{subfigure}\hfill
	\begin{subfigure}[t]{0.425\linewidth}
		\centering
		\begin{tikzpicture}[trim axis right]
			\begin{axis}[
				cycle list name=Lista2,
				yticklabel style={text width=2.5em, align=right},
				width=\linewidth,
				height=\linewidth,
				xlabel={$x$},
				xmin=0,
				xmax=1,
				ymax=0.12,
				ymin=-0.12,
				ytick={-0.1,0,0.1}, 
				xminorticks=false, yminorticks=false,
				legend style={at={(0.8,0.01)},anchor=south east},
				grid=both,
				legend entries={$p=1$\\ $p=2$\\ $p=3$\\ $p=4$\\}
				]
				%\addplot table [x=phyX, y=eu_start, col sep=comma] {\DATAFILEpA};
				\addplot table [x=phyX, y expr=\thisrow{uh_end} - \thisrow{uh_start}, col sep=comma] {\DATAFILEpAz};
				\addplot table [x=phyX, y expr=\thisrow{uh_end} - \thisrow{uh_start}, col sep=comma] {\DATAFILEpBz};
				\addplot table [x=phyX, y expr=\thisrow{uh_end} - \thisrow{uh_start}, col sep=comma] {\DATAFILEpCz};
				\addplot table [x=phyX, y expr=\thisrow{uh_end} - \thisrow{uh_start}, col sep=comma] {\DATAFILEpDz};
			\end{axis}
		\end{tikzpicture}
	\end{subfigure}
	\hspace*{\fill}	
	\caption{Approximation of the periodic problem of Section~\ref{s:Dispersion} with initial data \eqref{eq:bumpInitial}.
	Top row: exact and discrete solutions at final time $T=2$; bottom row: error at final time.
	Left column: \igastab\ method; right column: \femstab\ method.
	For both the stabilizations and all the spline degrees, $N_{\mathrm{dof}}=17\,424$ and $h_t \approx 2 h_s$.}
	\label{fig:dispersionBump}
\end{figure}

%% file: Graphics/DispersionTestErrorsBump.tex
\tikzstyle{Linea1}=[thick,dashed, mark options={solid}]
\tikzstyle{Linea2}=[thick]%,mark=*]
\tikzstyle{Linea3}=[thick,dashed]

%\pgfplotsset{
	%    tick label style={font=\footnotesize},
	%    label style={font=\footnotesize},
	%    legend style={font=\footnotesize},
	%}

\pgfplotscreateplotcyclelist{Lista1}{%
	{Linea2,Red,mark=*},
	{Linea1,Red,mark=*},
	%{Linea3,Red},
	{Linea2,Green,mark=triangle*},
	{Linea1,Green,mark=triangle*},
	%{Linea3,Green},
	{Linea2,Cyan,mark=square*},
	{Linea1,Cyan,mark=square*},
	%{Linea3,Cyan},
	{Linea2,Violet,mark=diamond*},
	{Linea1,Violet,mark=diamond*},
	%{Linea3,Violet},
}

\def \DATAFILEo {Graphics/ErrorTableDispersionBumpOur.csv}
\def \DATAFILEz {Graphics/ErrorTableDispersionBumpZank.csv}

\begin{figure}[htbp]
	\centering
	\hspace*{\fill}
	\begin{subfigure}[t]{0.425\linewidth}
		\centering
		\begin{tikzpicture}[%font=\footnotesize, 
			%trim axis left
			]
			\begin{semilogyaxis}[
				cycle list name=Lista1,
				width=\linewidth,
				height=\linewidth,
				xlabel={$p$},
				ymin=1e-5,
				ymax=5e-1,
				xminorticks=false,
				yminorticks=false,
				%max space between ticks=20,
				%xtick=data,
				%xticklabel={\pgfmathparse{exp(\tick)}\pgfmathprintnumber[fixed relative,precision=3]{\pgfmathresult}},
				%ytick={0.4,0.6,0.8,1}, 
				%yticklabels={$0.4$,$0.6$,$0.8$,$1$}, 
				ylabel={$\| u(T) - u_h(T) \|_{L^2(\Omega)}/\| u(T) \|_{L^2(\Omega)}$},
				legend columns=1,
				%legend pos=south east,
				legend style={at={(0.99,0.99)},anchor=north east},
				xmajorgrids=true,
				ymajorgrids=true,
				xtick={1,2,3,4},
				legend entries={\igastab,\femstab}]
				\addplot table [x=p, y=ErrL2, col sep=comma] {\DATAFILEo};
				\addplot table [x=p, y=ErrL2, col sep=comma] {\DATAFILEz};
				%\addplot table [x=NdofP1, y=NdofL2P1O, col sep=comma] {\DATAFILE};
			\end{semilogyaxis}
		\end{tikzpicture}
	\end{subfigure}
	\hfill
	\begin{subfigure}[t]{0.425\linewidth}
		\centering
		\begin{tikzpicture}[%font=\footnotesize,
			%trim axis right
			]
			\begin{semilogyaxis}[
				cycle list name=Lista1,
				width=\linewidth,
				height=\linewidth,
				xlabel={$p$},
				ymin=1e-5,
				ymax=5e-1,
				xminorticks=false,
				yminorticks=false,
				%max space between ticks=20,
				%xtick=data,
				%xticklabel={\pgfmathparse{exp(\tick)}\pgfmathprintnumber[fixed relative,precision=3]{\pgfmathresult}},
				%ytick={0.4,0.6,0.8,1}, 
				%yticklabels={$0.4$,$0.6$,$0.8$,$1$}, 
				ylabel={$| u(T) - u_h(T) |_{H^1(\Omega)}/| u(T) |_{H^1(\Omega)}$},
				legend columns=1,
				%legend pos=south east,
				legend style={at={(0.01,0.01)},anchor=south west},
				xmajorgrids=true,
				ymajorgrids=true,
				xtick={1,2,3,4},	
				legend entries={\igastab,\femstab}]
				\addplot table [x=p, y=ErrH1s, col sep=comma] {\DATAFILEo};
				\addplot table [x=p, y=ErrH1s, col sep=comma] {\DATAFILEz};
			\end{semilogyaxis}
		\end{tikzpicture}
	\end{subfigure}
	\hspace*{\fill}
	\caption{Comparison between the relative errors at final time of the \igastab\ and the \femstab\ methods, for the periodic bump wave problem \eqref{eq:bumpInitial}.
	For both the stabilizations and all the spline degrees, $N_{\mathrm{dof}}=17\,424$ and $h_t \approx 2 h_s$.
	}
	\label{fig:StabComparisonBump}
\end{figure}

%% file: Graphics/Dispersion_p1234TentRotByP.tex
\tikzstyle{Linea1}=[thick]
\tikzstyle{Linea2}=[thick,mark=*]

%\pgfplotsset{
%    tick label style={font=\footnotesize},
%    label style={font=\footnotesize},
%    legend style={font=\footnotesize},
%}

\pgfplotscreateplotcyclelist{Lista1}{%
{Linea1,Red},
%{Linea2,Red},
{Linea1,Green},
%{Linea2,Green},
{Linea1,Cyan},
%{Linea2,Cyan},
{Linea1,Violet},
%{Linea2,Violet},
{Linea1,BurntOrange},
%{Linea2,BurntOrange},
{Linea1,Brown},
%{Linea2,Brown}
{Linea1,Magenta},
%{Linea2,Magenta}
{Linea1,Blue},
%{Linea2,Goldenrod}
}

\def\DATAFILEcAo{Graphics/CnTableTentOur_p1.csv}
\def\DATAFILEcBo{Graphics/CnTableTentOur_p2.csv}
\def\DATAFILEcCo{Graphics/CnTableTentOur_p3.csv}
\def\DATAFILEcDo{Graphics/CnTableTentOur_p4.csv}
\def\DATAFILEcAz{Graphics/CnTableTentZank_p1.csv}
\def\DATAFILEcBz{Graphics/CnTableTentZank_p2.csv}
\def\DATAFILEcCz{Graphics/CnTableTentZank_p3.csv}
\def\DATAFILEcDz{Graphics/CnTableTentZank_p4.csv}

\def\HeiFac{0.4}

\begin{sidewaysfigure}[htbp]
%\hfill
\begin{subfigure}[t]{0.475\linewidth}
\begin{tikzpicture}[trim axis right]
\begin{semilogyaxis}[
cycle list name=Lista1,
width=.97\linewidth,
height=\HeiFac\linewidth,
xmin=0,	xmax=2,
ymin=0.0000001, ymax=3,
ytick={1,0.01,0.0001,0.000001},
%ytick={1,10^(-3),10^(-6),10^(-9)},
%ytick={1,.01,.0001,.000001},
%title={$|\arg(c_n) - \arg(c_{h_n})|%_{\mathrm{mod}\,\pi}
%	$, \igastab},
title={\igastab},
ylabel={\fbox{\rotatebox{270}{$c_1$}}}, %ylabel style={align=right},
grid=both,
% legend columns=4, %8,
% legend pos=south east,
% legend style={at={(0.99,0.01)},anchor=south east},
% legend entries={$p=1$,$p=2$,$p=3$,$p=4$}
% legend entries={$c_1$,$c_2$,$c_3$,$c_5$,$c_6$,$c_7$,$c_9$,$c_{10}$}
]
\addplot table [x=phyT, y=c1, col sep=comma] {\DATAFILEcAo};
\addplot table [x=phyT, y=c1, col sep=comma] {\DATAFILEcBo};
\addplot table [x=phyT, y=c1, col sep=comma] {\DATAFILEcCo};
\addplot table [x=phyT, y=c1, col sep=comma] {\DATAFILEcDo};
\end{semilogyaxis}
\end{tikzpicture}
\end{subfigure}
\hfill
\begin{subfigure}[t]{0.475\linewidth}
\begin{tikzpicture}[trim axis right]
\begin{semilogyaxis}[
cycle list name=Lista1,
width=.97\linewidth,
height=\HeiFac\linewidth,
xmin=0,	xmax=2,
ymin=0.0000001, ymax=3,
ytick={1,0.01,0.0001,0.000001},
%ytick={1,10^(-2),10^(-4),10^(-6)},  % Works on Windows but not on Linux
%ytick={1,.01,.0001,.000001},
%title={$| \arg(c_n) - \arg(c_{h_n})|%_{\mathrm{mod}\,\pi}
%$, \femstab},
title={\femstab},
grid=both,
% legend columns=8,
% legend pos=south east,
% legend style={at={(0.99,0.01)},anchor=south east},
% legend entries={$c_1$,$c_2$,$c_3$,$c_5$,$c_6$,$c_7$,$c_9$,$c_{10}$}
]
\addplot table [x=phyT, y=c1, col sep=comma] {\DATAFILEcAz};
\addplot table [x=phyT, y=c1, col sep=comma] {\DATAFILEcBz};
\addplot table [x=phyT, y=c1, col sep=comma] {\DATAFILEcCz};
\addplot table [x=phyT, y=c1, col sep=comma] {\DATAFILEcDz};
\end{semilogyaxis}
\end{tikzpicture}
\end{subfigure}
%%%%%%%%%%%%%%%%%%%%%%%%%%%%%%%%%%%%%%%%%%%%%%%%%%%%%%%%%%%%%%%%%%%%%%%%%%%%%%%%%%%%%%%%%%%%%%%%%%%%%%
% p=2
%\hfill
\begin{subfigure}[t]{0.475\linewidth}
\begin{tikzpicture}[trim axis right]
\begin{semilogyaxis}[
cycle list name=Lista1,
width=.97\linewidth,
height=\HeiFac\linewidth,
xmin=0,	xmax=2,
ymin=0.0000001, ymax=3,
ytick={1,0.01,0.0001,0.000001},
ylabel={\fbox{\rotatebox{270}{$c_2$}}},%	ylabel={$| \arg(c_n) - \arg(c_{h_n}) |$},
grid=both,
% legend columns=8,
% legend pos=south east,
% legend style={at={(0.99,0.01)},anchor=south east},
% legend entries={$c_1$,$c_2$,$c_3$,$c_5$,$c_6$,$c_7$,$c_9$,$c_{10}$}
]
\addplot table [x=phyT, y=c2, col sep=comma] {\DATAFILEcAo};
\addplot table [x=phyT, y=c2, col sep=comma] {\DATAFILEcBo};
\addplot table [x=phyT, y=c2, col sep=comma] {\DATAFILEcCo};
\addplot table [x=phyT, y=c2, col sep=comma] {\DATAFILEcDo};
\end{semilogyaxis}
\end{tikzpicture}
\end{subfigure}
\hfill
\begin{subfigure}[t]{0.475\linewidth}
\begin{tikzpicture}[trim axis right]
\begin{semilogyaxis}[
cycle list name=Lista1,
width=.97\linewidth,
height=\HeiFac\linewidth,
xmin=0,	xmax=2,
ymin=0.0000001, ymax=3,
ytick={1,0.01,0.0001,0.000001},
%ytick={1,.01,.0001,.000001}, 
%ylabel={$\fbox{p=2}$},%	ylabel={$| \arg(c_n) - \arg(c_{h_n}) |$},
grid=both,
% legend columns=8,
% legend pos=south east,
% legend style={at={(0.99,0.01)},anchor=south east},
% legend entries={$c_1$,$c_2$,$c_3$,$c_5$,$c_6$,$c_7$,$c_9$,$c_{10}$}
]
\addplot table [x=phyT, y=c2, col sep=comma] {\DATAFILEcAz};
\addplot table [x=phyT, y=c2, col sep=comma] {\DATAFILEcBz};
\addplot table [x=phyT, y=c2, col sep=comma] {\DATAFILEcCz};
\addplot table [x=phyT, y=c2, col sep=comma] {\DATAFILEcDz};
\end{semilogyaxis}
\end{tikzpicture}
\end{subfigure}
%%%%%%%%%%%%%%%%%%%%%%%%%%%%%%%%%%%%%%%%%%%%%%%%%%%%%%%%%%%%%%%%%%%%%%%%%%%%%%%%%%%%%%%%%%%%%%%%%%%%%%
% p=3
%\hfill
\begin{subfigure}[t]{0.475\linewidth}
\begin{tikzpicture}[trim axis right]
\begin{semilogyaxis}[
cycle list name=Lista1,
width=.97\linewidth,
height=\HeiFac\linewidth,
xmin=0,	xmax=2,
ymin=0.0000001, ymax=3,
ytick={1,0.01,0.0001,0.000001},
ylabel={\fbox{\rotatebox{270}{$c_3$}}},%	ylabel={$| \arg(c_n) - \arg(c_{h_n}) |$},
grid=both,
]
\addplot table [x=phyT, y=c3, col sep=comma] {\DATAFILEcAo};
\addplot table [x=phyT, y=c3, col sep=comma] {\DATAFILEcBo};
\addplot table [x=phyT, y=c3, col sep=comma] {\DATAFILEcCo};
\addplot table [x=phyT, y=c3, col sep=comma] {\DATAFILEcDo};
\end{semilogyaxis}
\end{tikzpicture}
\end{subfigure}
\hfill
\begin{subfigure}[t]{0.475\linewidth}
\begin{tikzpicture}[trim axis right]
\begin{semilogyaxis}[
cycle list name=Lista1,
width=.97\linewidth,
height=\HeiFac\linewidth,
xmin=0,	xmax=2,
ymin=0.0000001, ymax=3,
ytick={1,0.01,0.0001,0.000001},
%ylabel={$\fbox{p=3}$},%	ylabel={$| \arg(c_n) - \arg(c_{h_n}) |$},
grid=both,
]
\addplot table [x=phyT, y=c3, col sep=comma] {\DATAFILEcAz};
\addplot table [x=phyT, y=c3, col sep=comma] {\DATAFILEcBz};
\addplot table [x=phyT, y=c3, col sep=comma] {\DATAFILEcCz};
\addplot table [x=phyT, y=c3, col sep=comma] {\DATAFILEcDz};
\end{semilogyaxis}
\end{tikzpicture}
\end{subfigure}
%%%%%%%%%%%%%%%%%%%%%%%%%%%%%%%%%%%%%%%%%%%%%%%%%%%%%%%%%%%%%%%%%%%%%%%%%%%%%%%%%%%%%%%%%%%%%%%%%%%%%%%%
%% p=4
%\hfill
\begin{subfigure}[t]{0.475\linewidth}
\begin{tikzpicture}[trim axis right]\begin{semilogyaxis}[
cycle list name=Lista1,
width=.97\linewidth,
height=\HeiFac\linewidth,
xlabel={$t$},
xmin=0,	xmax=2,
ymin=0.0000001, ymax=3,
ytick={1,0.01,0.0001,0.000001},
ylabel={\fbox{\rotatebox{270}{$c_5$}}},
grid=both,
]
\addplot table [x=phyT, y=c5, col sep=comma] {\DATAFILEcAo};
\addplot table [x=phyT, y=c5, col sep=comma] {\DATAFILEcBo};
\addplot table [x=phyT, y=c5, col sep=comma] {\DATAFILEcCo};
\addplot table [x=phyT, y=c5, col sep=comma] {\DATAFILEcDo};
\end{semilogyaxis}
\end{tikzpicture}
\end{subfigure}
\hfill
\begin{subfigure}[t]{0.475\linewidth}
\begin{tikzpicture}[trim axis right]\begin{semilogyaxis}[
cycle list name=Lista1,
width=.97\linewidth,
height=\HeiFac\linewidth,
xlabel={$t$},
xmin=0,	xmax=2,
ymin=0.0000001, ymax=3,
ytick={1,0.01,0.0001,0.000001},
grid=both,
legend columns=4, 
legend pos=south east,
legend style={at={(0.99,0.01)},anchor=south east},
legend entries={$p=1$,$p=2$,$p=3$,$p=4$}
]
\addplot table [x=phyT, y=c5, col sep=comma] {\DATAFILEcAz};
\addplot table [x=phyT, y=c5, col sep=comma] {\DATAFILEcBz};
\addplot table [x=phyT, y=c5, col sep=comma] {\DATAFILEcCz};
\addplot table [x=phyT, y=c5, col sep=comma] {\DATAFILEcDz};
\end{semilogyaxis}
\end{tikzpicture}
\end{subfigure}
\caption{Phase errors of the largest 4 Fourier coefficients ($c_1,c_2,c_3,c_5$) for the periodic problem with initial conditions \eqref{eq:tentInitial} (tent), approximated with the \igastab\ method~\eqref{equivstabiga} (left) and the \femstab\ method~\eqref{stabZank} (right), $N_{\mathrm{dof}}=17\,424$ and $h_t \approx 2h_s$.
}
\label{fig:dispersion1234TentRotByP}
\end{sidewaysfigure}

%% file: Graphics/Dispersion_p1234BumpRotByP.tex
\tikzstyle{Linea1}=[thick]
\tikzstyle{Linea2}=[thick,mark=*]

%\pgfplotsset{
%    tick label style={font=\footnotesize},
%    label style={font=\footnotesize},
%    legend style={font=\footnotesize},
%}

\pgfplotscreateplotcyclelist{Lista1}{%
{Linea1,Red},
%{Linea2,Red},
{Linea1,Green},
%{Linea2,Green},
{Linea1,Cyan},
%{Linea2,Cyan},
{Linea1,Violet},
%{Linea2,Violet},
{Linea1,BurntOrange},
%{Linea2,BurntOrange},
{Linea1,Brown},
%{Linea2,Brown}
{Linea1,Magenta},
%{Linea2,Magenta}
{Linea1,Blue},
%{Linea2,Goldenrod}
}

\def\DATAFILEcAo{Graphics/CnTableBumpOur_p1.csv}
\def\DATAFILEcBo{Graphics/CnTableBumpOur_p2.csv}
\def\DATAFILEcCo{Graphics/CnTableBumpOur_p3.csv}
\def\DATAFILEcDo{Graphics/CnTableBumpOur_p4.csv}
\def\DATAFILEcAz{Graphics/CnTableBumpZank_p1.csv}
\def\DATAFILEcBz{Graphics/CnTableBumpZank_p2.csv}
\def\DATAFILEcCz{Graphics/CnTableBumpZank_p3.csv}
\def\DATAFILEcDz{Graphics/CnTableBumpZank_p4.csv}

\def\HeiFac{0.4}

\begin{sidewaysfigure}[htbp]
%\hfill
\begin{subfigure}[t]{0.475\linewidth}
\begin{tikzpicture}[trim axis right]
\begin{semilogyaxis}[
cycle list name=Lista1,
width=.97\linewidth,
height=\HeiFac\linewidth,
xmin=0,	xmax=2,
ymin=0.0000001, ymax=3,
ytick={1,0.01,0.0001,0.000001},
%ytick={1,10^(-3),10^(-6),10^(-9)},
%ytick={1,.01,.0001,.000001},
%title={$|\arg(c_n) - \arg(c_{h_n})|%_{\mathrm{mod}\,\pi}
%	$, \igastab},
title={\igastab},
ylabel={\fbox{\rotatebox{270}{$c_1$}}}, %ylabel style={align=right},
grid=both,
% legend columns=4, %8,
% legend pos=south east,
% legend style={at={(0.99,0.01)},anchor=south east},
% legend entries={$p=1$,$p=2$,$p=3$,$p=4$}
% legend entries={$c_1$,$c_2$,$c_3$,$c_5$,$c_6$,$c_7$,$c_9$,$c_{10}$}
]
\addplot table [x=phyT, y=c1, col sep=comma] {\DATAFILEcAo};
\addplot table [x=phyT, y=c1, col sep=comma] {\DATAFILEcBo};
\addplot table [x=phyT, y=c1, col sep=comma] {\DATAFILEcCo};
\addplot table [x=phyT, y=c1, col sep=comma] {\DATAFILEcDo};
\end{semilogyaxis}
\end{tikzpicture}
\end{subfigure}
\hfill
\begin{subfigure}[t]{0.475\linewidth}
\begin{tikzpicture}[trim axis right]
\begin{semilogyaxis}[
cycle list name=Lista1,
width=.97\linewidth,
height=\HeiFac\linewidth,
xmin=0,	xmax=2,
ymin=0.0000001, ymax=3,
ytick={1,0.01,0.0001,0.000001},
%ytick={1,10^(-2),10^(-4),10^(-6)},  % Works on Windows but not on Linux
%ytick={1,.01,.0001,.000001},
%title={$| \arg(c_n) - \arg(c_{h_n})|%_{\mathrm{mod}\,\pi}
%$, \femstab},
title={\femstab},
grid=both,
% legend columns=8,
% legend pos=south east,
% legend style={at={(0.99,0.01)},anchor=south east},
% legend entries={$c_1$,$c_2$,$c_3$,$c_5$,$c_6$,$c_7$,$c_9$,$c_{10}$}
]
\addplot table [x=phyT, y=c1, col sep=comma] {\DATAFILEcAz};
\addplot table [x=phyT, y=c1, col sep=comma] {\DATAFILEcBz};
\addplot table [x=phyT, y=c1, col sep=comma] {\DATAFILEcCz};
\addplot table [x=phyT, y=c1, col sep=comma] {\DATAFILEcDz};
\end{semilogyaxis}
\end{tikzpicture}
\end{subfigure}
%%%%%%%%%%%%%%%%%%%%%%%%%%%%%%%%%%%%%%%%%%%%%%%%%%%%%%%%%%%%%%%%%%%%%%%%%%%%%%%%%%%%%%%%%%%%%%%%%%%%%%
% p=2
%\hfill
\begin{subfigure}[t]{0.475\linewidth}
\begin{tikzpicture}[trim axis right]
\begin{semilogyaxis}[
cycle list name=Lista1,
width=.97\linewidth,
height=\HeiFac\linewidth,
xmin=0,	xmax=2,
ymin=0.0000001, ymax=3,
ytick={1,0.01,0.0001,0.000001},
ylabel={\fbox{\rotatebox{270}{$c_2$}}},%	ylabel={$| \arg(c_n) - \arg(c_{h_n}) |$},
grid=both,
% legend columns=8,
% legend pos=south east,
% legend style={at={(0.99,0.01)},anchor=south east},
% legend entries={$c_1$,$c_2$,$c_3$,$c_5$,$c_6$,$c_7$,$c_9$,$c_{10}$}
]
\addplot table [x=phyT, y=c2, col sep=comma] {\DATAFILEcAo};
\addplot table [x=phyT, y=c2, col sep=comma] {\DATAFILEcBo};
\addplot table [x=phyT, y=c2, col sep=comma] {\DATAFILEcCo};
\addplot table [x=phyT, y=c2, col sep=comma] {\DATAFILEcDo};
\end{semilogyaxis}
\end{tikzpicture}
\end{subfigure}
\hfill
\begin{subfigure}[t]{0.475\linewidth}
\begin{tikzpicture}[trim axis right]
\begin{semilogyaxis}[
cycle list name=Lista1,
width=.97\linewidth,
height=\HeiFac\linewidth,
xmin=0,	xmax=2,
ymin=0.0000001, ymax=3,
ytick={1,0.01,0.0001,0.000001},
%ytick={1,.01,.0001,.000001}, 
%ylabel={$\fbox{p=2}$},%	ylabel={$| \arg(c_n) - \arg(c_{h_n}) |$},
grid=both,
% legend columns=8,
% legend pos=south east,
% legend style={at={(0.99,0.01)},anchor=south east},
% legend entries={$c_1$,$c_2$,$c_3$,$c_5$,$c_6$,$c_7$,$c_9$,$c_{10}$}
]
\addplot table [x=phyT, y=c2, col sep=comma] {\DATAFILEcAz};
\addplot table [x=phyT, y=c2, col sep=comma] {\DATAFILEcBz};
\addplot table [x=phyT, y=c2, col sep=comma] {\DATAFILEcCz};
\addplot table [x=phyT, y=c2, col sep=comma] {\DATAFILEcDz};
\end{semilogyaxis}
\end{tikzpicture}
\end{subfigure}
%%%%%%%%%%%%%%%%%%%%%%%%%%%%%%%%%%%%%%%%%%%%%%%%%%%%%%%%%%%%%%%%%%%%%%%%%%%%%%%%%%%%%%%%%%%%%%%%%%%%%%
% p=3
%\hfill
\begin{subfigure}[t]{0.475\linewidth}
\begin{tikzpicture}[trim axis right]
\begin{semilogyaxis}[
cycle list name=Lista1,
width=.97\linewidth,
height=\HeiFac\linewidth,
xmin=0,	xmax=2,
ymin=0.0000001, ymax=3,
ytick={1,0.01,0.0001,0.000001},
ylabel={\fbox{\rotatebox{270}{$c_3$}}},%	ylabel={$| \arg(c_n) - \arg(c_{h_n}) |$},
grid=both,
]
\addplot table [x=phyT, y=c3, col sep=comma] {\DATAFILEcAo};
\addplot table [x=phyT, y=c3, col sep=comma] {\DATAFILEcBo};
\addplot table [x=phyT, y=c3, col sep=comma] {\DATAFILEcCo};
\addplot table [x=phyT, y=c3, col sep=comma] {\DATAFILEcDo};
\end{semilogyaxis}
\end{tikzpicture}
\end{subfigure}
\hfill
\begin{subfigure}[t]{0.475\linewidth}
\begin{tikzpicture}[trim axis right]
\begin{semilogyaxis}[
cycle list name=Lista1,
width=.97\linewidth,
height=\HeiFac\linewidth,
xmin=0,	xmax=2,
ymin=0.0000001, ymax=3,
ytick={1,0.01,0.0001,0.000001},
%ylabel={$\fbox{p=3}$},%	ylabel={$| \arg(c_n) - \arg(c_{h_n}) |$},
grid=both,
]
\addplot table [x=phyT, y=c3, col sep=comma] {\DATAFILEcAz};
\addplot table [x=phyT, y=c3, col sep=comma] {\DATAFILEcBz};
\addplot table [x=phyT, y=c3, col sep=comma] {\DATAFILEcCz};
\addplot table [x=phyT, y=c3, col sep=comma] {\DATAFILEcDz};
\end{semilogyaxis}
\end{tikzpicture}
\end{subfigure}
%%%%%%%%%%%%%%%%%%%%%%%%%%%%%%%%%%%%%%%%%%%%%%%%%%%%%%%%%%%%%%%%%%%%%%%%%%%%%%%%%%%%%%%%%%%%%%%%%%%%%%%%
%% p=4
%\hfill
\begin{subfigure}[t]{0.475\linewidth}
\begin{tikzpicture}[trim axis right]\begin{semilogyaxis}[
cycle list name=Lista1,
width=.97\linewidth,
height=\HeiFac\linewidth,
xlabel={$t$},
xmin=0,	xmax=2,
ymin=0.0000001, ymax=3,
ytick={1,0.01,0.0001,0.000001},
ylabel={\fbox{\rotatebox{270}{$c_4$}}},
grid=both,
]
\addplot table [x=phyT, y=c4, col sep=comma] {\DATAFILEcAo};
\addplot table [x=phyT, y=c4, col sep=comma] {\DATAFILEcBo};
\addplot table [x=phyT, y=c4, col sep=comma] {\DATAFILEcCo};
\addplot table [x=phyT, y=c4, col sep=comma] {\DATAFILEcDo};
\end{semilogyaxis}
\end{tikzpicture}
\end{subfigure}
\hfill
\begin{subfigure}[t]{0.475\linewidth}
\begin{tikzpicture}[trim axis right]\begin{semilogyaxis}[
cycle list name=Lista1,
width=.97\linewidth,
height=\HeiFac\linewidth,
xlabel={$t$},
xmin=0,	xmax=2,
ymin=0.0000001, ymax=3,
ytick={1,0.01,0.0001,0.000001},
grid=both,
legend columns=4, 
legend pos=south east,
legend style={at={(0.99,0.01)},anchor=south east},
legend entries={$p=1$,$p=2$,$p=3$,$p=4$}
]
\addplot table [x=phyT, y=c4, col sep=comma] {\DATAFILEcAz};
\addplot table [x=phyT, y=c4, col sep=comma] {\DATAFILEcBz};
\addplot table [x=phyT, y=c4, col sep=comma] {\DATAFILEcCz};
\addplot table [x=phyT, y=c4, col sep=comma] {\DATAFILEcDz};
\end{semilogyaxis}
\end{tikzpicture}
\end{subfigure}
\caption{Phase errors of the largest 4 Fourier coefficients ($c_1,c_2,c_3,c_4$) for the periodic problem with initial conditions \eqref{eq:bumpInitial} (bump), approximated with the \igastab\ method~\eqref{equivstabiga} (left) and the \femstab\ method~\eqref{stabZank} (right), $N_{\mathrm{dof}}=17\,424$ and $h_t \approx 2h_s$.
%\am{------SECOND VERSION------}
}
\label{fig:dispersion1234BumpRotByP}
\end{sidewaysfigure}

%% file: Chapters/Conclusions.tex
\section{Conclusions}
\label{sec:conclusions}
In this work, we devised a high-order unconditionally stable space--time isogeometric discretization for the linear acoustic wave equation, based on arbitrary-degree splines with maximal regularity. While for  a plain Galerkin space--time isogeometric formulation a CFL condition is required to ensure stability,  motivated by \cite[Lemma 4.2.26]{Z2020}, we have proposed a stabilization technique that extends the one proposed in \cite{SZ2019,Z2020}.

Our method is unconditionally stable for splines of maximal regularity $C^{p-1}$ in time, just as those presented in \cite{SZ2019,Z2020}  are unconditionally stable for $C^0$  finite elements. 
We observe numerically that our method delivers optimal convergence rates, energy conservation and desirable dissipation and dispersion properties.
The stabilization for splines of intermediate regularity $C^k$ with $1\le k\le p_t-1$, for example, enhancing the proposed stabilized term with additional penalty terms on the jump of higher-order derivatives, is also interesting and postponed to future works.

The proposed space-time formulation provides advantages in terms of stability, that holds without any CFL condition, and accuracy per number of degrees-of-freedom. However, it comes with a high computational cost when the linear system is solved using standard direct solvers, as is the case in this work. It is noteworthy that our method results in a global system in the $(d+1)$-dimensional space-time domain. The development of an efficient solver, following approaches such as~\cite{LMST2020}, is a crucial step and will be the focus of an upcoming work.

We remark that the proposed stabilized formulation is supported by a wide range of numerical experiments but no rigorous theoretical results is available yet, in particular for what concerns its unconditional stability.

%% file: Chapters/Appendix.tex
\section{Appendix}
\label{sec:appendix}

In this appendix, we collect the numerical experiments corresponding to the results presented in Table~\ref{table:stab_comparison}.
All plots refer to the relative errors incurred in the approximation of the smooth exact solution~\eqref{u1_ex}, measured in $L^2(Q)$-norm (left panels) and $H^1(Q)$-seminorm (right panels).
The relative errors are plotted either against $h_t$ with a fixed ratio $h_t/h_s=5$, or against the ratio $h_t/h_s$, with fixed $h_t=0.1562$, to highlight where a CFL condition is present.
\input{Graphics/ErrorNoStabRobustTot}
\input{Graphics/ErrorStabConvTot}
\input{Graphics/ErrorStabRobustTot}
\input{Graphics/ErrorStabZankPConvTot}
\input{Graphics/ErrorStabZankPRobustTot}

%% file: Graphics/ErrorNoStabRobustTot.tex
\tikzstyle{Linea1}=[dashed]
\tikzstyle{Linea2}=[thick]%,mark=*]

%\pgfplotsset{
	%    tick label style={font=\footnotesize},
	%    label style={font=\footnotesize},
	%    legend style={font=\footnotesize},
	%}
% 
% \pgfplotscreateplotcyclelist{Lista1}{%
	% 	{Linea2,Red},
	% 	%{Linea1,Red},
	% 	{Linea2,Green},
	% 	%{Linea1,Green},
	% 	{Linea2,Cyan},
	% 	%{Linea1,Cyan},
	% 	{Linea2,Violet},
	% 	%{Linea1,Violet}
	% }
\pgfplotscreateplotcyclelist{Lista1}{%
	{Linea2,Red,mark=*},
	% 	{Linea1,Red},
	{Linea2,Green,mark=triangle*},
	% 	{Linea1,Green},
	{Linea2,Cyan,mark=square*},
	% 	{Linea1,Cyan},
	{Linea2,Violet,mark=diamond*},
	% 	{Linea1,Violet}
}

\pgfplotscreateplotcyclelist{Lista2}{%
	%{Linea2,Red,mark=*},
	%{Linea1,Red},
	{Linea2,Green,mark=triangle*},
	%{Linea1,Green},
	{Linea2,Cyan,mark=square*},
	%{Linea1,Cyan},
	{Linea2,Violet,mark=diamond*},
	%{Linea1,Violet}
}

\def \DATAFILE {Graphics/ErrorTableNoStabRobustRST.csv}

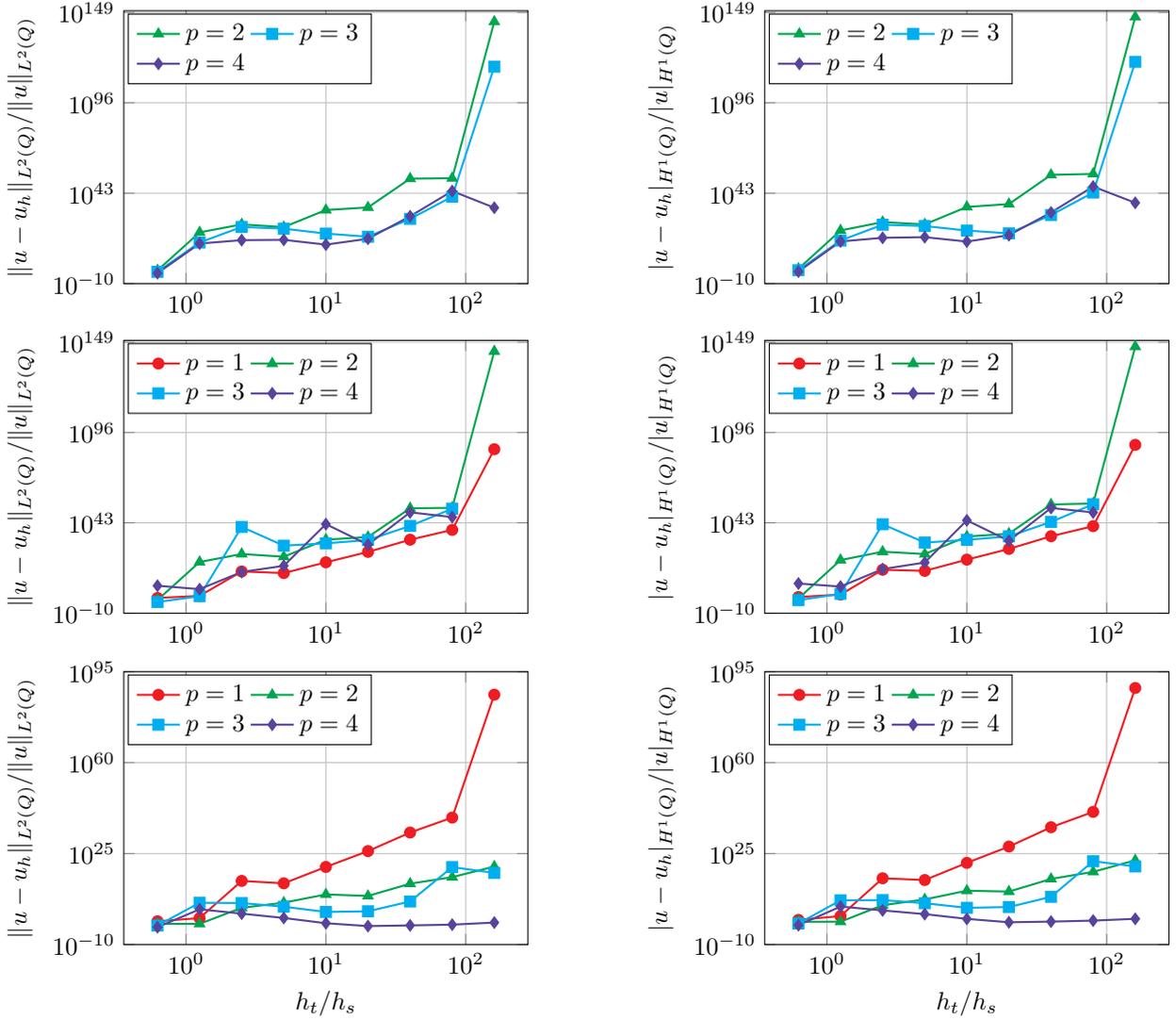
\begin{figure}[htbp]
	\centering
	\hspace*{\fill}
	\begin{subfigure}[t]{0.425\linewidth}
		\centering
		\begin{tikzpicture}[%font=\footnotesize, 
			trim axis left]
			\begin{loglogaxis}[
				cycle list name=Lista2,
				width=\linewidth,
				height=0.75\linewidth,
				ymin=1e-10,
				ymax=1e150,
				xminorticks=false,
				yminorticks=false,
				%max space between ticks=20,
				%xtick=data,
				%xticklabel={\pgfmathparse{exp(\tick)}\pgfmathprintnumber[fixed relative,precision=3]{\pgfmathresult}},
				%ytick={1,1e35,1e70,1e105,1e140}, 
				%yticklabels={$0.4$,$0.6$,$0.8$,$1$}, 
				ylabel={$\| u - u_h \|_{L^2(Q)}/\| u \|_{L^2(Q)}$},
				legend columns=2,
				%legend pos=south east,
				legend style={at={(0.01,0.99)},anchor= north west},
				xmajorgrids=true,
				ymajorgrids=true,
				legend entries={$p=2$,$p=3$,$p=4$}]
				%\addplot table [x=htfhs, y=L2P1R2, col sep=comma] {\DATAFILE};
				\addplot table [x=htfhs, y=L2P2R2, col sep=comma] {\DATAFILE};
				\addplot table [x=htfhs, y=L2P3R2, col sep=comma] {\DATAFILE};
				\addplot table [x=htfhs, y=L2P4R2, col sep=comma] {\DATAFILE};
			\end{loglogaxis}
		\end{tikzpicture}
	\end{subfigure}\hfill
	\begin{subfigure}[t]{0.425\linewidth}
		\centering
		\begin{tikzpicture}[%font=\footnotesize,
			trim axis right]
			\begin{loglogaxis}[
				cycle list name=Lista2,
				width=\linewidth,
				height=0.75\linewidth,
				ymin=1e-10,
				ymax=1e150,
				xminorticks=false,
				yminorticks=false,
				%max space between ticks=20,
				%xtick=data,
				%xticklabel={\pgfmathparse{exp(\tick)}\pgfmathprintnumber[fixed relative,precision=3]{\pgfmathresult}},
				%ytick={1,1e35,1e70,1e105,1e140}, 
				%yticklabels={$0.4$,$0.6$,$0.8$,$1$}, 
				ylabel={$| u - u_h |_{H^1(Q)}/| u |_{H^1(Q)}$},
				legend columns=2,
				%legend pos=south east,
				legend style={at={(0.01,0.99)},anchor= north west},
				xmajorgrids=true,
				ymajorgrids=true,
				legend entries={$p=2$,$p=3$,$p=4$}]
				%\addplot table [x=htfhs, y=H1P1R2, col sep=comma] {\DATAFILE};
				\addplot table [x=htfhs, y=H1P2R2, col sep=comma] {\DATAFILE};
				\addplot table [x=htfhs, y=H1P3R2, col sep=comma] {\DATAFILE};
				\addplot table [x=htfhs, y=H1P4R2, col sep=comma] {\DATAFILE};
			\end{loglogaxis}
		\end{tikzpicture}
	\end{subfigure}
	\hspace*{\fill}
	\\
	\hspace*{\fill}
		\begin{subfigure}[t]{0.425\linewidth}
			\centering
		\begin{tikzpicture}[%font=\footnotesize, 
			trim axis left]
			\begin{loglogaxis}[
				cycle list name=Lista1,
				width=\linewidth,
				height=0.75\linewidth,
				ymin=1e-10,
				ymax=1e150,
				xminorticks=false,
				yminorticks=false,
				%max space between ticks=20,
				%xtick=data,
				%xticklabel={\pgfmathparse{exp(\tick)}\pgfmathprintnumber[fixed relative,precision=3]{\pgfmathresult}},
				%ytick={1,1e35,1e70,1e105,1e140}, 
				%yticklabels={$0.4$,$0.6$,$0.8$,$1$}, 
				ylabel={$\| u - u_h \|_{L^2(Q)}/\| u \|_{L^2(Q)}$},
				legend columns=2,
				%legend pos=south east,
				legend style={at={(0.01,0.99)},anchor= north west},
				xmajorgrids=true,
				ymajorgrids=true,
				legend entries={$p=1$,$p=2$,$p=3$,$p=4$}]
				\addplot table [x=htfhs, y=L2P1R0, col sep=comma] {\DATAFILE};
				\addplot table [x=htfhs, y=L2P2R0, col sep=comma] {\DATAFILE};
				\addplot table [x=htfhs, y=L2P3R0, col sep=comma] {\DATAFILE};
				\addplot table [x=htfhs, y=L2P4R0, col sep=comma] {\DATAFILE};
			\end{loglogaxis}
		\end{tikzpicture}
	\end{subfigure}\hfill
	\begin{subfigure}[t]{0.425\linewidth}
		\centering
		\begin{tikzpicture}[%font=\footnotesize,
			trim axis right]
			\begin{loglogaxis}[
				cycle list name=Lista1,
				width=\linewidth,
				height=0.75\linewidth,
				ymin=1e-10,
				ymax=1e150,
				xminorticks=false,
				yminorticks=false,
				%max space between ticks=20,
				%xtick=data,
				%xticklabel={\pgfmathparse{exp(\tick)}\pgfmathprintnumber[fixed relative,precision=3]{\pgfmathresult}},
				%ytick={1,1e35,1e70,1e105,1e140},
				%yticklabels={$0.4$,$0.6$,$0.8$,$1$}, 
				ylabel={$| u - u_h |_{H^1(Q)}/| u |_{H^1(Q)}$},
				legend columns=2,
				%legend pos=south east,
				legend style={at={(0.01,0.99)},anchor= north west},
				xmajorgrids=true,
				ymajorgrids=true,
				legend entries={$p=1$,$p=2$,$p=3$,$p=4$}]
				\addplot table [x=htfhs, y=H1P1R0, col sep=comma] {\DATAFILE};
				\addplot table [x=htfhs, y=H1P2R0, col sep=comma] {\DATAFILE};
				\addplot table [x=htfhs, y=H1P3R0, col sep=comma] {\DATAFILE};
				\addplot table [x=htfhs, y=H1P4R0, col sep=comma] {\DATAFILE};
			\end{loglogaxis}
		\end{tikzpicture}
	\end{subfigure}
	\hspace*{\fill}
	\\
	\hspace*{\fill}
	\begin{subfigure}[t]{0.425\linewidth}
		\centering
		\begin{tikzpicture}[%font=\footnotesize, 
			trim axis left]
			\begin{loglogaxis}[
				cycle list name=Lista1,
				width=\linewidth,
				height=0.75\linewidth,
				xlabel={$h_t/h_s$},
				ymin=1e-10,
				ymax=1e95,
				xminorticks=false,
				yminorticks=false,
				%max space between ticks=20,
				%xtick=data,
				%xticklabel={\pgfmathparse{exp(\tick)}\pgfmathprintnumber[fixed relative,precision=3]{\pgfmathresult}},
				%ytick={1,1e20,1e41,1e61,1e81}, 
				%yticklabels={$0.4$,$0.6$,$0.8$,$1$}, 
				ylabel={$\| u - u_h \|_{L^2(Q)}/\| u \|_{L^2(Q)}$},
				legend columns=2,
				%legend pos=south east,
				legend style={at={(0.01,0.99)},anchor= north west},
				xmajorgrids=true,
				ymajorgrids=true,
				legend entries={$p=1$,$p=2$,$p=3$,$p=4$}]
				\addplot table [x=htfhs, y=L2P1R, col sep=comma] {\DATAFILE};
				\addplot table [x=htfhs, y=L2P2R, col sep=comma] {\DATAFILE};
				\addplot table [x=htfhs, y=L2P3R, col sep=comma] {\DATAFILE};
				\addplot table [x=htfhs, y=L2P4R, col sep=comma] {\DATAFILE};
			\end{loglogaxis}
		\end{tikzpicture}
	\end{subfigure}\hfill
	\begin{subfigure}[t]{0.425\linewidth}
		\centering
		\begin{tikzpicture}[%font=\footnotesize,
			trim axis right]
			\begin{loglogaxis}[
				cycle list name=Lista1,
				width=\linewidth,
				height=0.75\linewidth,
				xlabel={$h_t/h_s$},
				ymin=1e-10,
				ymax=1e95,
				xminorticks=false,
				yminorticks=false,
				%max space between ticks=20,
				%xtick=data,
				%xticklabel={\pgfmathparse{exp(\tick)}\pgfmathprintnumber[fixed relative,precision=3]{\pgfmathresult}},
				%ytick={1,1e20,1e41,1e61,1e81}, 
				%yticklabels={$0.4$,$0.6$,$0.8$,$1$}, 
				ylabel={$| u - u_h |_{H^1(Q)}/| u |_{H^1(Q)}$},
				legend columns=2,
				%legend pos=south east,
				legend style={at={(0.01,0.99)},anchor= north west},
				xmajorgrids=true,
				ymajorgrids=true,
				legend entries={$p=1$,$p=2$,$p=3$,$p=4$}]
				\addplot table [x=htfhs, y=H1P1R, col sep=comma] {\DATAFILE};
				\addplot table [x=htfhs, y=H1P2R, col sep=comma] {\DATAFILE};
				\addplot table [x=htfhs, y=H1P3R, col sep=comma] {\DATAFILE};
				\addplot table [x=htfhs, y=H1P4R, col sep=comma] {\DATAFILE};
			\end{loglogaxis}
		\end{tikzpicture}
	\end{subfigure}
	\hspace*{\fill}
	\caption{Relative errors of the non-stabilized isogeometric method~\eqref{igawave} with $C^{p-2}$-continuous splines (top), $C^{0}$-continuous splines (middle), splines of maximal regularity (bottom), plotted against the ratio $h_t/h_s$.We observe that the non-stabilized formulation is subject to a CFL condition in all cases.}
	\label{fig:NoStabRobustTot}
\end{figure}

%% file: Graphics/ErrorStabConvTot.tex
\tikzstyle{Linea1}=[dashed]
\tikzstyle{Linea2}=[thick]%,mark=*]

%\pgfplotsset{
%    tick label style={font=\footnotesize},
%    label style={font=\footnotesize},
%    legend style={font=\footnotesize},
%}

\pgfplotscreateplotcyclelist{Lista1}{%
%{Linea2,Red,mark=*},
%{Linea1,Red},
{Linea2,Green,mark=triangle*},
%{Linea1,Green},
{Linea2,Cyan,mark=square*},
%{Linea1,Cyan},
{Linea2,Violet,mark=diamond*},
%{Linea1,Violet}
}

\pgfplotscreateplotcyclelist{Lista2}{%
{Linea2,Red,mark=*},
%{Linea1,Red},
{Linea2,Green,mark=triangle*},
%{Linea1,Green},
{Linea2,Cyan,mark=square*},
%{Linea1,Cyan},
{Linea2,Violet,mark=diamond*},
%{Linea1,Violet}
}

\pgfplotscreateplotcyclelist{Lista3}{%
{Linea2,Red,mark=*},
{Linea1},
{Linea2,Green,mark=triangle*},
{Linea1},
{Linea2,Cyan,mark=square*},
{Linea1},
{Linea2,Violet,mark=diamond*},
{Linea1}
}

\def \DATAFILE {Graphics/ErrorTableStabConvRST.csv}

\begin{figure}[htbp]
\centering
\hspace*{\fill}
\begin{subfigure}[t]{0.425\linewidth}
\centering
\begin{tikzpicture}[%font=\footnotesize, 
trim axis left]
\begin{loglogaxis}[
yticklabel style={text width=2.125em, align=right},
cycle list name=Lista1,
width=\linewidth,
height=0.75\linewidth,
ymin=1e-10,
ymax=1e70,
xminorticks=false,
yminorticks=false,
%max space between ticks=20,
%xtick=data,
%xticklabel={\pgfmathparse{exp(\tick)}\pgfmathprintnumber[fixed relative,precision=3]{\pgfmathresult}},
%ytick={1e0,1e23,1e46,1e69,1e92}, 
%yticklabels={$0.4$,$0.6$,$0.8$,$1$}, 
ylabel={$\| u - u_h \|_{L^2(Q)}/\| u \|_{L^2(Q)}$},
legend columns=1,
%legend pos=south east,
legend style={at={(0.99,0.99)},anchor=north east},
xmajorgrids=true,
ymajorgrids=true,
legend entries={$p=2$,$p=3$,$p=4$}]
%\addplot table [x=ht, y=L2P1R2, col sep=comma] {\DATAFILE};
\addplot table [x=ht, y=L2P2R2, col sep=comma] {\DATAFILE};
\addplot table [x=ht, y=L2P3R2, col sep=comma] {\DATAFILE};
\addplot table [x=ht, y=L2P4R2, col sep=comma] {\DATAFILE};
\end{loglogaxis}
\end{tikzpicture}
\end{subfigure}\hfill
\begin{subfigure}[t]{0.425\linewidth}
\centering
\begin{tikzpicture}[%font=\footnotesize,
trim axis right]
\begin{loglogaxis}[
yticklabel style={text width=2.125em, align=right},
cycle list name=Lista1,
width=\linewidth,
height=0.75\linewidth,
ymin=1e-10,
ymax=1e70,
xminorticks=false,
yminorticks=false,
%max space between ticks=20,
%xtick=data,
%xticklabel={\pgfmathparse{exp(\tick)}\pgfmathprintnumber[fixed relative,precision=3]{\pgfmathresult}},
%ytick={1e0,1e23,1e46,1e69,1e92}, 
%yticklabels={$0.4$,$0.6$,$0.8$,$1$}, 
ylabel={$| u - u_h |_{H^1(Q)}/| u |_{H^1(Q)}$},
legend columns=1,
%legend pos=south east,
legend style={at={(0.99,0.99)},anchor=north east},
xmajorgrids=true,
ymajorgrids=true,
legend entries={$p=2$,$p=3$,$p=4$}]
%\addplot table [x=ht, y=H1sP1R2, col sep=comma] {\DATAFILE};
\addplot table [x=ht, y=H1sP2R2, col sep=comma] {\DATAFILE};
\addplot table [x=ht, y=H1sP3R2, col sep=comma] {\DATAFILE};
\addplot table [x=ht, y=H1sP4R2, col sep=comma] {\DATAFILE};
\end{loglogaxis}
\end{tikzpicture}
\end{subfigure}
\hspace*{\fill}
\\
\hspace*{\fill}
\begin{subfigure}[t]{0.425\linewidth}
\centering
\begin{tikzpicture}[%font=\footnotesize, 
trim axis left]
\begin{loglogaxis}[
yticklabel style={text width=2.125em, align=right},
cycle list name=Lista2,
width=\linewidth,
height=0.75\linewidth,
ymin=1e-10,
ymax=1e100,
xminorticks=false,
yminorticks=false,
%max space between ticks=20,
%xtick=data,
%xticklabel={\pgfmathparse{exp(\tick)}\pgfmathprintnumber[fixed relative,precision=3]{\pgfmathresult}},
%ytick={1e-3,1e18,1e40,1e61,1e82}, 
%yticklabels={$0.4$,$0.6$,$0.8$,$1$}, 
ylabel={$\| u - u_h \|_{L^2(Q)}/\| u \|_{L^2(Q)}$},
legend columns=1,
%legend pos=south east,
legend style={at={(0.99,0.99)},anchor=north east},
xmajorgrids=true,
ymajorgrids=true,
legend entries={$p=1$,$p=2$,$p=3$,$p=4$}]
\addplot table [x=ht, y=L2P1R0, col sep=comma] {\DATAFILE};
\addplot table [x=ht, y=L2P2R0, col sep=comma] {\DATAFILE};
\addplot table [x=ht, y=L2P3R0, col sep=comma] {\DATAFILE};
\addplot table [x=ht, y=L2P4R0, col sep=comma] {\DATAFILE};
\end{loglogaxis}
\end{tikzpicture}
\end{subfigure}\hfill
\begin{subfigure}[t]{0.425\linewidth}
\centering
\begin{tikzpicture}[%font=\footnotesize,
trim axis right]
\begin{loglogaxis}[
yticklabel style={text width=2.125em, align=right},
cycle list name=Lista2,
width=\linewidth,
height=0.75\linewidth,
ymin=1e-10,
ymax=1e100,
xminorticks=false,
yminorticks=false,
%max space between ticks=20,
%xtick=data,
%xticklabel={\pgfmathparse{exp(\tick)}\pgfmathprintnumber[fixed relative,precision=3]{\pgfmathresult}},
%ytick={1e-3,1e18,1e40,1e61,1e82}, 
%yticklabels={$0.4$,$0.6$,$0.8$,$1$}, 
ylabel={$| u - u_h |_{H^1(Q)}/| u |_{H^1(Q)}$},
legend columns=1,
%legend pos=south east,
legend style={at={(0.99,0.99)},anchor=north east},
xmajorgrids=true,
ymajorgrids=true,
legend entries={$p=1$,$p=2$,$p=3$,$p=4$}]
\addplot table [x=ht, y=H1sP1R0, col sep=comma] {\DATAFILE};
\addplot table [x=ht, y=H1sP2R0, col sep=comma] {\DATAFILE};
\addplot table [x=ht, y=H1sP3R0, col sep=comma] {\DATAFILE};
\addplot table [x=ht, y=H1sP4R0, col sep=comma] {\DATAFILE};
\end{loglogaxis}
\end{tikzpicture}
\end{subfigure}
\hspace*{\fill}
\\
\hspace*{\fill}
\begin{subfigure}[t]{0.425\linewidth}
\centering
\begin{tikzpicture}[%font=\footnotesize, 
trim axis left]
\begin{loglogaxis}[
yticklabel style={text width=2.125em, align=right},
cycle list name=Lista3,
width=\linewidth,
height=0.9\linewidth,
xlabel={$h_t$},
ymin=1e-9,
ymax=1e1,
xminorticks=false,
yminorticks=false,
%max space between ticks=20,
%xtick=data,
%xticklabel={\pgfmathparse{exp(\tick)}\pgfmathprintnumber[fixed relative,precision=3]{\pgfmathresult}},
%ytick={0.4,0.6,0.8,1}, 
%ytick={1e-7,1e-5,1e-3,1e-1,1e1}, 
ylabel={$\| u - u_h \|_{L^2(Q)}/\| u \|_{L^2(Q)}$},
legend columns=1,
%legend pos=south east,
legend style={at={(0.99,0.01)},anchor=south east},
xmajorgrids=true,
ymajorgrids=true,
legend entries={$p=1$,$O(h^2)$,$p=2$,$O(h^3)$,$p=3$,$O(h^4)$,$p=4$,$O(h^5)$}]
\addplot table [x=ht, y=L2P1RsR0, col sep=comma] {\DATAFILE};
\addplot table [x=ht, y=h2L2RsR0, col sep=comma] {\DATAFILE};
\addplot table [x=ht, y=L2P2RsR0, col sep=comma] {\DATAFILE};
\addplot table [x=ht, y=h3L2RsR0, col sep=comma] {\DATAFILE};
\addplot table [x=ht, y=L2P3RsR0, col sep=comma] {\DATAFILE};
\addplot table [x=ht, y=h4L2RsR0, col sep=comma] {\DATAFILE};
\addplot table [x=ht, y=L2P4RsR0, col sep=comma] {\DATAFILE};
\addplot table [x=ht, y=h5L2RsR0, col sep=comma] {\DATAFILE};
\end{loglogaxis}
\end{tikzpicture}
\end{subfigure}\hfill
\begin{subfigure}[t]{0.425\linewidth}
\centering
\begin{tikzpicture}[%font=\footnotesize,
trim axis right]
\begin{loglogaxis}[
yticklabel style={text width=2.125em, align=right},
cycle list name=Lista3,
width=\linewidth,
height=0.9\linewidth,
xlabel={$h_t$},
ymin=1e-9,
ymax=1e1,
xminorticks=false,
yminorticks=false,
%max space between ticks=20,
%xtick=data,
%xticklabel={\pgfmathparse{exp(\tick)}\pgfmathprintnumber[fixed relative,precision=3]{\pgfmathresult}},
%ytick={1e-7,1e-5,1e-3,1e-1,1e1},
%yticklabels={$0.4$,$0.6$,$0.8$,$1$}, 
ylabel={$| u - u_h |_{H^1(Q)}/| u |_{H^1(Q)}$},
legend columns=1,
%legend pos=south east,
legend style={at={(0.99,0.01)},anchor=south east},
xmajorgrids=true,
ymajorgrids=true,
legend entries={$p=1$,$O(h)$,$p=2$,$O(h^2)$,$p=3$,$O(h^3)$,$p=4$,$O(h^4)$}]
\addplot table [x=ht, y=H1sP1RsR0, col sep=comma] {\DATAFILE};
\addplot table [x=ht, y=h1H1sRsR0, col sep=comma] {\DATAFILE};
\addplot table [x=ht, y=H1sP2RsR0, col sep=comma] {\DATAFILE};
\addplot table [x=ht, y=h2H1sRsR0, col sep=comma] {\DATAFILE};
\addplot table [x=ht, y=H1sP3RsR0, col sep=comma] {\DATAFILE};
\addplot table [x=ht, y=h3H1sRsR0, col sep=comma] {\DATAFILE};
\addplot table [x=ht, y=H1sP4RsR0, col sep=comma] {\DATAFILE};
\addplot table [x=ht, y=h4H1sRsR0, col sep=comma] {\DATAFILE};
\end{loglogaxis}
\end{tikzpicture}
\end{subfigure}
\hspace*{\fill}
\caption{Relative errors of the \igastab\ method \eqref{equivstabiga} with $C^{p-2}$-continuous splines (top), $C^{0}$-continuous splines (middle), splines of maximal regularity w.r.t.\ time variable and $C^0$-continuous w.r.t.\ space variable (bottom), plotted against the time mesh-size $h_t$.
We observe that the key requirement for the stability and the optimal convergence of the IGA-Stab formulation is the maximal regularity in time of the basis functions.}
\label{fig:StabConvTot}
\end{figure}
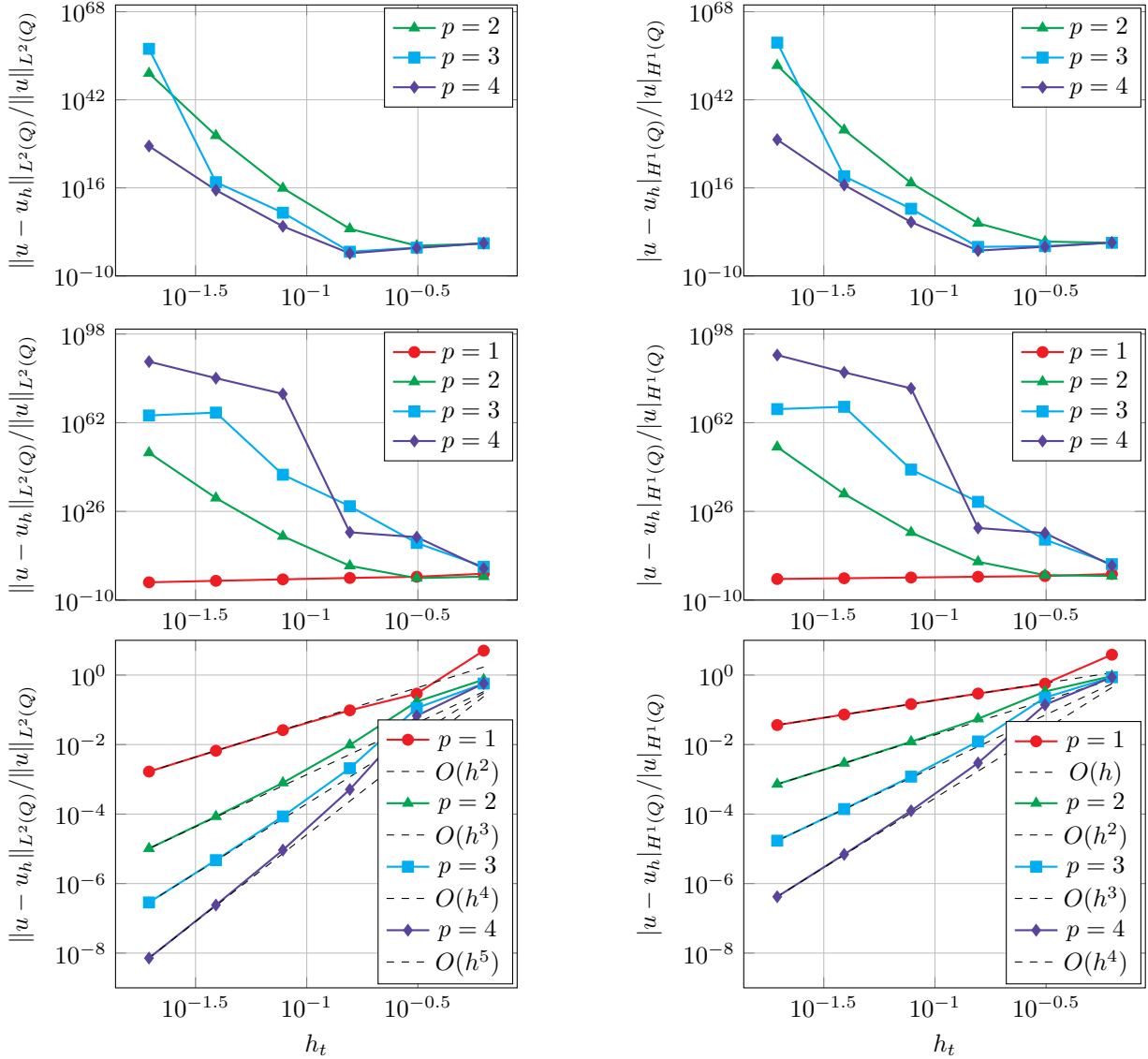

%% file: Graphics/ErrorStabRobustTot.tex
\tikzstyle{Linea1}=[dashed]
\tikzstyle{Linea2}=[thick]%,mark=*]

%\pgfplotsset{
%    tick label style={font=\footnotesize},
%    label style={font=\footnotesize},
%    legend style={font=\footnotesize},
%}
% 
% \pgfplotscreateplotcyclelist{Lista1}{%
% 	{Linea2,Red},
% 	%{Linea1,Red},
% 	{Linea2,Green},
% 	%{Linea1,Green},
% 	{Linea2,Cyan},
% 	%{Linea1,Cyan},
% 	{Linea2,Violet},
% 	%{Linea1,Violet}
% }
\pgfplotscreateplotcyclelist{Lista1}{%
%{Linea2,Red,mark=*},
% 	{Linea1,Red},
{Linea2,Green,mark=triangle*},
% 	{Linea1,Green},
{Linea2,Cyan,mark=square*},
% 	{Linea1,Cyan},
{Linea2,Violet,mark=diamond*},
% 	{Linea1,Violet}
}

\pgfplotscreateplotcyclelist{Lista2}{%
{Linea2,Red,mark=*},
%{Linea1,Red},
{Linea2,Green,mark=triangle*},
%{Linea1,Green},
{Linea2,Cyan,mark=square*},
%{Linea1,Cyan},
{Linea2,Violet,mark=diamond*},
%{Linea1,Violet}
}

\def \DATAFILE {Graphics/ErrorTableStabRobustRST.csv}

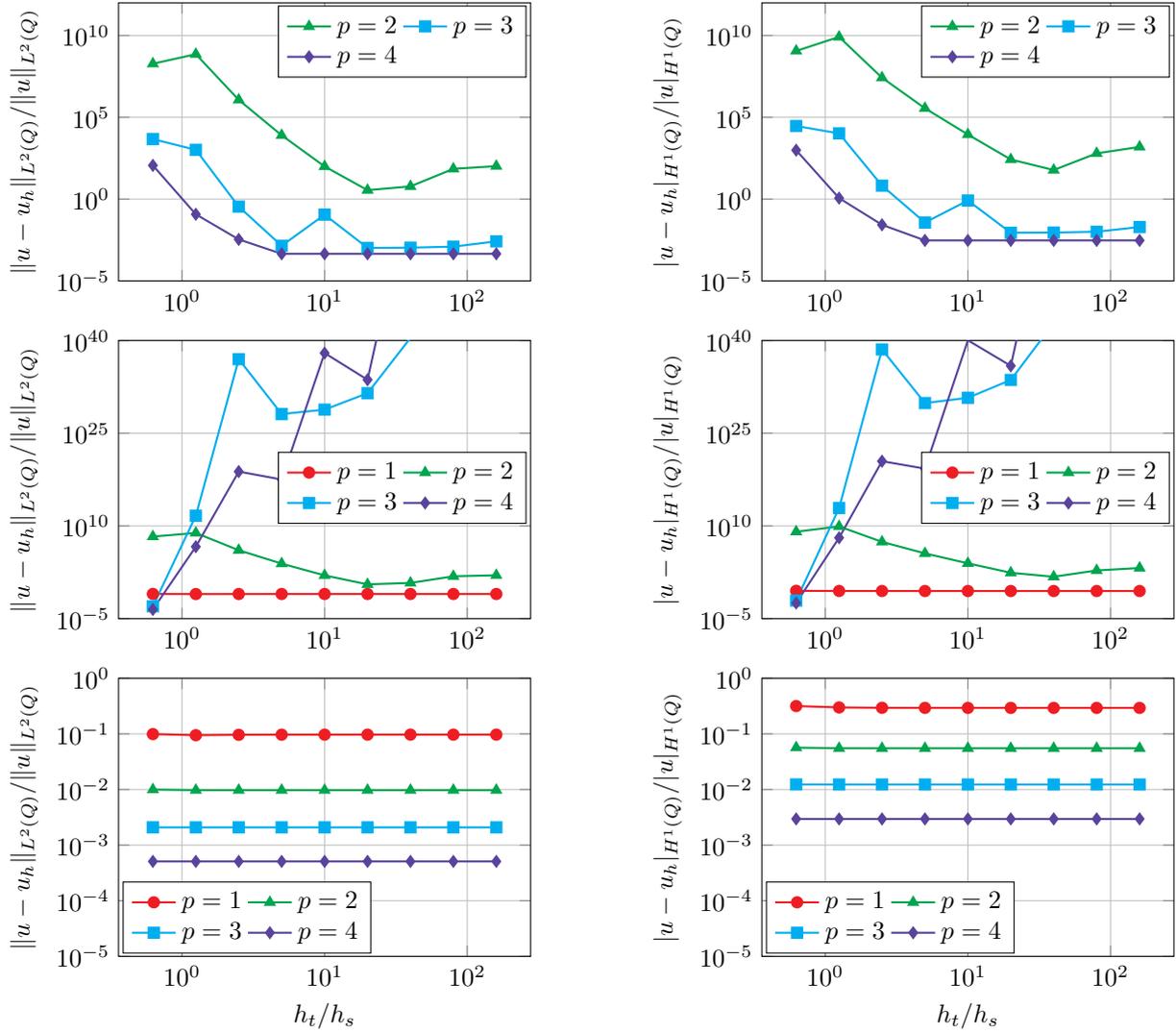
\begin{figure}[htbp]
\centering
\hspace*{\fill}
\begin{subfigure}[t]{0.425\linewidth}
\centering
\begin{tikzpicture}[%font=\footnotesize, 
trim axis left]
\begin{loglogaxis}[
cycle list name=Lista1,
width=\linewidth,
height=0.75\linewidth,
ymin=1e-5,
ymax=1e12,
xminorticks=false,
yminorticks=false,
%max space between ticks=20,
%xtick=data,
%xticklabel={\pgfmathparse{exp(\tick)}\pgfmathprintnumber[fixed relative,precision=3]{\pgfmathresult}},
%ytick={0.4,0.6,0.8,1}, 
%yticklabels={$0.4$,$0.6$,$0.8$,$1$}, 
ylabel={$\| u - u_h \|_{L^2(Q)}/\| u \|_{L^2(Q)}$},
legend columns=2,
%legend pos=south east,
legend style={at={(0.99,0.99)},anchor= north east},
xmajorgrids=true,
ymajorgrids=true,
legend entries={$p=2$,$p=3$,$p=4$}]
%\addplot table [x=htfhs, y=L2P1R2, col sep=comma] {\DATAFILE};
\addplot table [x=htfhs, y=L2P2R2, col sep=comma] {\DATAFILE};
\addplot table [x=htfhs, y=L2P3R2, col sep=comma] {\DATAFILE};
\addplot table [x=htfhs, y=L2P4R2, col sep=comma] {\DATAFILE};
\end{loglogaxis}
\end{tikzpicture}
\end{subfigure}\hfill
\begin{subfigure}[t]{0.425\linewidth}
\centering
\begin{tikzpicture}[%font=\footnotesize,
trim axis right]
\begin{loglogaxis}[
cycle list name=Lista1,
width=\linewidth,
height=0.75\linewidth,
ymin=1e-5,
ymax=1e12,
xminorticks=false,
yminorticks=false,
%max space between ticks=20,
%xtick=data,
%xticklabel={\pgfmathparse{exp(\tick)}\pgfmathprintnumber[fixed relative,precision=3]{\pgfmathresult}},
%ytick={0.4,0.6,0.8,1}, 
%yticklabels={$0.4$,$0.6$,$0.8$,$1$}, 
ylabel={$| u - u_h |_{H^1(Q)}/| u |_{H^1(Q)}$},
legend columns=2,
%legend pos=south east,
legend style={at={(0.99,0.99)},anchor= north east},
xmajorgrids=true,
ymajorgrids=true,
legend entries={$p=2$,$p=3$,$p=4$}]
%\addplot table [x=htfhs, y=H1P1R2, col sep=comma] {\DATAFILE};
\addplot table [x=htfhs, y=H1P2R2, col sep=comma] {\DATAFILE};
\addplot table [x=htfhs, y=H1P3R2, col sep=comma] {\DATAFILE};
\addplot table [x=htfhs, y=H1P4R2, col sep=comma] {\DATAFILE};
\end{loglogaxis}
\end{tikzpicture}
\end{subfigure}
\hspace*{\fill}
\\
\hspace*{\fill}
\begin{subfigure}[t]{0.425\linewidth}
\centering
\begin{tikzpicture}[%font=\footnotesize, 
trim axis left]
\begin{loglogaxis}[
cycle list name=Lista2,
width=\linewidth,
height=0.75\linewidth,
ymin=1e-5,
ymax=1e40,
xminorticks=false,
yminorticks=false,
%max space between ticks=20,
%xtick=data,
%xticklabel={\pgfmathparse{exp(\tick)}\pgfmathprintnumber[fixed relative,precision=3]{\pgfmathresult}},
%ytick={1e-2,1e0,1e10,1e42,1e91}, 
%yticklabels={$0.4$,$0.6$,$0.8$,$1$}, 
ylabel={$\| u - u_h \|_{L^2(Q)}/\| u \|_{L^2(Q)}$},
legend columns=2,
%legend pos=south east,
legend style={at={(0.99,0.6)},anchor= north east},
xmajorgrids=true,
ymajorgrids=true,
legend entries={$p=1$,$p=2$,$p=3$,$p=4$}]
\addplot table [x=htfhs, y=L2P1R0, col sep=comma] {\DATAFILE};
\addplot table [x=htfhs, y=L2P2R0, col sep=comma] {\DATAFILE};
\addplot table [x=htfhs, y=L2P3R0, col sep=comma] {\DATAFILE};
\addplot table [x=htfhs, y=L2P4R0, col sep=comma] {\DATAFILE};
\end{loglogaxis}
\end{tikzpicture}
\end{subfigure}\hfill
\begin{subfigure}[t]{0.425\linewidth}
\centering
\begin{tikzpicture}[%font=\footnotesize,
trim axis right]
\begin{loglogaxis}[
cycle list name=Lista2,
width=\linewidth,
height=0.75\linewidth,
ymin=1e-5,
ymax=1e40,
xminorticks=false,
yminorticks=false,
%max space between ticks=20,
%xtick=data,
%xticklabel={\pgfmathparse{exp(\tick)}\pgfmathprintnumber[fixed relative,precision=3]{\pgfmathresult}},
%ytick={0.4,0.6,0.8,1}, 
%yticklabels={$0.4$,$0.6$,$0.8$,$1$}, 
ylabel={$| u - u_h |_{H^1(Q)}/| u |_{H^1(Q)}$},
legend columns=2,
%legend pos=south east,
legend style={at={(0.99,0.6)},anchor= north east},
xmajorgrids=true,
ymajorgrids=true,
legend entries={$p=1$,$p=2$,$p=3$,$p=4$}]
\addplot table [x=htfhs, y=H1P1R0, col sep=comma] {\DATAFILE};
\addplot table [x=htfhs, y=H1P2R0, col sep=comma] {\DATAFILE};
\addplot table [x=htfhs, y=H1P3R0, col sep=comma] {\DATAFILE};
\addplot table [x=htfhs, y=H1P4R0, col sep=comma] {\DATAFILE};
\end{loglogaxis}
\end{tikzpicture}
\end{subfigure}
\hspace*{\fill}
\\
\hspace*{\fill}
\begin{subfigure}[t]{0.425\linewidth}
\centering
\begin{tikzpicture}[%font=\footnotesize, 
trim axis left]
\begin{loglogaxis}[
cycle list name=Lista2,
width=\linewidth,
height=0.75\linewidth,
xlabel={$h_t/h_s$},
ymin=1e-5,
ymax=1e0,
xminorticks=false,
yminorticks=false,
%max space between ticks=20,
%xtick=data,
%xticklabel={\pgfmathparse{exp(\tick)}\pgfmathprintnumber[fixed relative,precision=3]{\pgfmathresult}},
%ytick={1e-3,1e-2,1e-1}, 
%yticklabels={$0.4$,$0.6$,$0.8$,$1$}, 
ylabel={$\| u - u_h \|_{L^2(Q)}/\| u \|_{L^2(Q)}$},
legend columns=2,
%legend pos=south east,
legend style={at={(0.01,0.01)},anchor= south west},
xmajorgrids=true,
ymajorgrids=true,
legend entries={$p=1$,$p=2$,$p=3$,$p=4$}]
\addplot table [x=htfhs, y=L2P1RsR0, col sep=comma] {\DATAFILE};
\addplot table [x=htfhs, y=L2P2RsR0, col sep=comma] {\DATAFILE};
\addplot table [x=htfhs, y=L2P3RsR0, col sep=comma] {\DATAFILE};
\addplot table [x=htfhs, y=L2P4RsR0, col sep=comma] {\DATAFILE};
\end{loglogaxis}
\end{tikzpicture}
\end{subfigure}\hfill
\begin{subfigure}[t]{0.425\linewidth}
\centering
\begin{tikzpicture}[%font=\footnotesize,
trim axis right]
\begin{loglogaxis}[
cycle list name=Lista2,
width=\linewidth,
height=0.75\linewidth,
xlabel={$h_t/h_s$},
ymin=1e-5,
ymax=1e0,
xminorticks=false,
yminorticks=false,
%max space between ticks=20,
%xtick=data,
%xticklabel={\pgfmathparse{exp(\tick)}\pgfmathprintnumber[fixed relative,precision=3]{\pgfmathresult}},
%ytick={1e-3,1e-2,1e-1}, 
%yticklabels={$0.4$,$0.6$,$0.8$,$1$}, 
ylabel={$| u - u_h |_{H^1(Q)}/| u |_{H^1(Q)}$},
legend columns=2,
%legend pos=south east,
legend style={at={(0.01,0.01)},anchor= south west},
xmajorgrids=true,
ymajorgrids=true,
legend entries={$p=1$,$p=2$,$p=3$,$p=4$}]
\addplot table [x=htfhs, y=H1P1RsR0, col sep=comma] {\DATAFILE};
\addplot table [x=htfhs, y=H1P2RsR0, col sep=comma] {\DATAFILE};
\addplot table [x=htfhs, y=H1P3RsR0, col sep=comma] {\DATAFILE};
\addplot table [x=htfhs, y=H1P4RsR0, col sep=comma] {\DATAFILE};
\end{loglogaxis}
\end{tikzpicture}
\end{subfigure}
\hspace*{\fill}
\caption{Relative errors of the \igastab\ method \eqref{equivstabiga} with $C^{p-2}$-continuous splines (top), $C^0$-continuous splines (middle), splines of maximal regularity w.r.t.\ time variable and $C^0$-continuous w.r.t.\ space variable (bottom), plotted against the ratio $h_t/h_s$. These plots confirm the observations made in the caption of Figure~\ref{fig:StabConvTot}.}
\label{fig:StabRobustTot}
\end{figure}

%% file: Graphics/ErrorStabZankPConvTot.tex
\tikzstyle{Linea1}=[dashed]
\tikzstyle{Linea2}=[thick]%,mark=*]

%\pgfplotsset{
%    tick label style={font=\footnotesize},
%    label style={font=\footnotesize},
%    legend style={font=\footnotesize},
%}

\pgfplotscreateplotcyclelist{Lista1}{%
{Linea2,Red,mark=*},
%{Linea1,Red},
{Linea2,Green,mark=triangle*},
%{Linea1,Green},
{Linea2,Cyan,mark=square*},
%{Linea1,Cyan},
{Linea2,Violet,mark=diamond*},
%{Linea1,Violet}
}

\pgfplotscreateplotcyclelist{Lista2}{%
%{Linea2,Red,mark=*},
% 	{Linea1,Red},
{Linea2,Green,mark=triangle*},
% 	{Linea1,Green},
{Linea2,Cyan,mark=square*},
% 	{Linea1,Cyan},
{Linea2,Violet,mark=diamond*},
% 	{Linea1,Violet}
}

\def \DATAFILE {Graphics/ErrorTableStabZankPConvRST.csv}

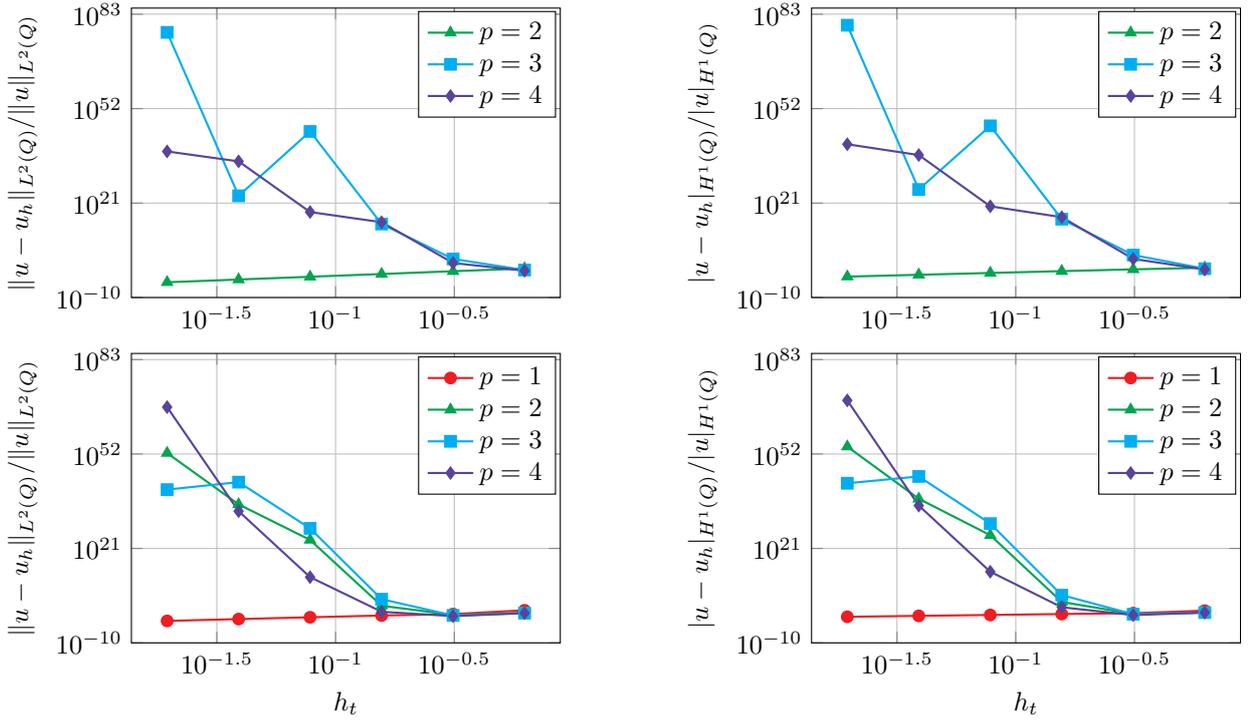
\begin{figure}[t]
\centering
\hspace*{\fill}
\begin{subfigure}[t]{0.425\linewidth}
\centering
\begin{tikzpicture}[%font=\footnotesize, 
trim axis left]
\begin{loglogaxis}[
cycle list name=Lista2,
width=\linewidth,
height=0.75\linewidth,
ymin=1e-10,
ymax=1e85,
xminorticks=false,
yminorticks=false,
%max space between ticks=20,
%xtick=data,
%xticklabel={\pgfmathparse{exp(\tick)}\pgfmathprintnumber[fixed relative,precision=3]{\pgfmathresult}},
%ytick={1e-5,1e14,1e33,1e52,1e71}, 
%yticklabels={$0.4$,$0.6$,$0.8$,$1$}, 
ylabel={$\| u - u_h \|_{L^2(Q)}/\| u \|_{L^2(Q)}$},
legend columns=1,
%legend pos=south east,
legend style={at={(0.99,0.99)},anchor=north east},
xmajorgrids=true,
ymajorgrids=true,
legend entries={$p=2$,$p=3$,$p=4$}]
%\addplot table [x=ht, y=L2P1R2, col sep=comma] {\DATAFILE};
\addplot table [x=ht, y=L2P2R2, col sep=comma] {\DATAFILE};
\addplot table [x=ht, y=L2P3R2, col sep=comma] {\DATAFILE};
\addplot table [x=ht, y=L2P4R2, col sep=comma] {\DATAFILE};
\end{loglogaxis}
\end{tikzpicture}
\end{subfigure}\hfill
\begin{subfigure}[t]{0.425\linewidth}
\centering
\begin{tikzpicture}[%font=\footnotesize,
trim axis right]
\begin{loglogaxis}[
cycle list name=Lista2,
width=\linewidth,
height=0.75\linewidth,
ymin=1e-10,
ymax=1e85,
xminorticks=false,
yminorticks=false,
%max space between ticks=20,
%xtick=data,
%xticklabel={\pgfmathparse{exp(\tick)}\pgfmathprintnumber[fixed relative,precision=3]{\pgfmathresult}},
%ytick={1e-5,1e14,1e33,1e52,1e71},
%yticklabels={$0.4$,$0.6$,$0.8$,$1$}, 
ylabel={$| u - u_h |_{H^1(Q)}/| u |_{H^1(Q)}$},
legend columns=1,
%legend pos=south east,
legend style={at={(0.99,0.99)},anchor=north east},
xmajorgrids=true,
ymajorgrids=true,
legend entries={$p=2$,$p=3$,$p=4$}]
%\addplot table [x=ht, y=H1sP1R2, col sep=comma] {\DATAFILE};
\addplot table [x=ht, y=H1sP2R2, col sep=comma] {\DATAFILE};
\addplot table [x=ht, y=H1sP3R2, col sep=comma] {\DATAFILE};
\addplot table [x=ht, y=H1sP4R2, col sep=comma] {\DATAFILE};
\end{loglogaxis}
\end{tikzpicture}
\end{subfigure}
\hspace*{\fill}
\\
\hspace*{\fill}
\begin{subfigure}[t]{0.425\linewidth}
\centering
\begin{tikzpicture}[%font=\footnotesize, 
trim axis left]
\begin{loglogaxis}[
cycle list name=Lista1,
width=\linewidth,
height=0.75\linewidth,
xlabel={$h_t$},
ymin=1e-10,
ymax=1e85,
xminorticks=false,
yminorticks=false,
%max space between ticks=20,
%xtick=data,
%xticklabel={\pgfmathparse{exp(\tick)}\pgfmathprintnumber[fixed relative,precision=3]{\pgfmathresult}},
%ytick={1e-5,1e14,1e33,1e52,1e71},
%yticklabels={$0.4$,$0.6$,$0.8$,$1$}, 
ylabel={$\| u - u_h \|_{L^2(Q)}/\| u \|_{L^2(Q)}$},
legend columns=1,
%legend pos=south east,
legend style={at={(0.99,0.99)},anchor=north east},
xmajorgrids=true,
ymajorgrids=true,
legend entries={$p=1$,$p=2$,$p=3$,$p=4$}]
\addplot table [x=ht, y=L2P1R, col sep=comma] {\DATAFILE};
\addplot table [x=ht, y=L2P2R, col sep=comma] {\DATAFILE};
\addplot table [x=ht, y=L2P3R, col sep=comma] {\DATAFILE};
\addplot table [x=ht, y=L2P4R, col sep=comma] {\DATAFILE};
\end{loglogaxis}
\end{tikzpicture}
\end{subfigure}\hfill
\begin{subfigure}[t]{0.425\linewidth}
\centering
\begin{tikzpicture}[%font=\footnotesize,
trim axis right]
\begin{loglogaxis}[
cycle list name=Lista1,
width=\linewidth,
height=0.75\linewidth,
xlabel={$h_t$},
ymin=1e-10,
ymax=1e85,
xminorticks=false,
yminorticks=false,
%max space between ticks=20,
%xtick=data,
%xticklabel={\pgfmathparse{exp(\tick)}\pgfmathprintnumber[fixed relative,precision=3]{\pgfmathresult}},
%ytick={1e-5,1e14,1e33,1e52,1e71},
%yticklabels={$0.4$,$0.6$,$0.8$,$1$}, 
ylabel={$| u - u_h |_{H^1(Q)}/| u |_{H^1(Q)}$},
legend columns=1,
legend pos=south east,
legend style={at={(0.99,0.99)},anchor=north east},
xmajorgrids=true,
ymajorgrids=true,
legend entries={$p=1$,$p=2$,$p=3$,$p=4$}]
\addplot table [x=ht, y=H1sP1R, col sep=comma] {\DATAFILE};
\addplot table [x=ht, y=H1sP2R, col sep=comma] {\DATAFILE};
\addplot table [x=ht, y=H1sP3R, col sep=comma] {\DATAFILE};
\addplot table [x=ht, y=H1sP4R, col sep=comma] {\DATAFILE};
\end{loglogaxis}
\end{tikzpicture}
\end{subfigure}
\hspace*{\fill}
\caption{Relative errors of the \femstab\ method \eqref{stabZank} with $C^{p-2}$-continuous splines (top), splines of maximal regularity (bottom), plotted against the time mesh-size $h_t$. We observe stability only in the lowest-regularity case.
}
\label{fig:StabZankPConvTot}
\end{figure}

%% file: Graphics/ErrorStabZankPRobustTot.tex
\tikzstyle{Linea1}=[dashed]
\tikzstyle{Linea2}=[thick]%,mark=*]

%\pgfplotsset{
%    tick label style={font=\footnotesize},
%    label style={font=\footnotesize},
%    legend style={font=\footnotesize},
%}
% 
% \pgfplotscreateplotcyclelist{Lista1}{%
% 	{Linea2,Red},
% 	%{Linea1,Red},
% 	{Linea2,Green},
% 	%{Linea1,Green},
% 	{Linea2,Cyan},
% 	%{Linea1,Cyan},
% 	{Linea2,Violet},
% 	%{Linea1,Violet}
% }
\pgfplotscreateplotcyclelist{Lista1}{%
{Linea2,Red,mark=*},
% 	{Linea1,Red},
{Linea2,Green,mark=triangle*},
% 	{Linea1,Green},
{Linea2,Cyan,mark=square*},
% 	{Linea1,Cyan},
{Linea2,Violet,mark=diamond*},
% 	{Linea1,Violet}
}

\pgfplotscreateplotcyclelist{Lista2}{%
%{Linea2,Red,mark=*},
% 	{Linea1,Red},
{Linea2,Green,mark=triangle*},
% 	{Linea1,Green},
{Linea2,Cyan,mark=square*},
% 	{Linea1,Cyan},
{Linea2,Violet,mark=diamond*},
% 	{Linea1,Violet}
}

\def \DATAFILE {Graphics/ErrorTableStabZankPRobustRST.csv}

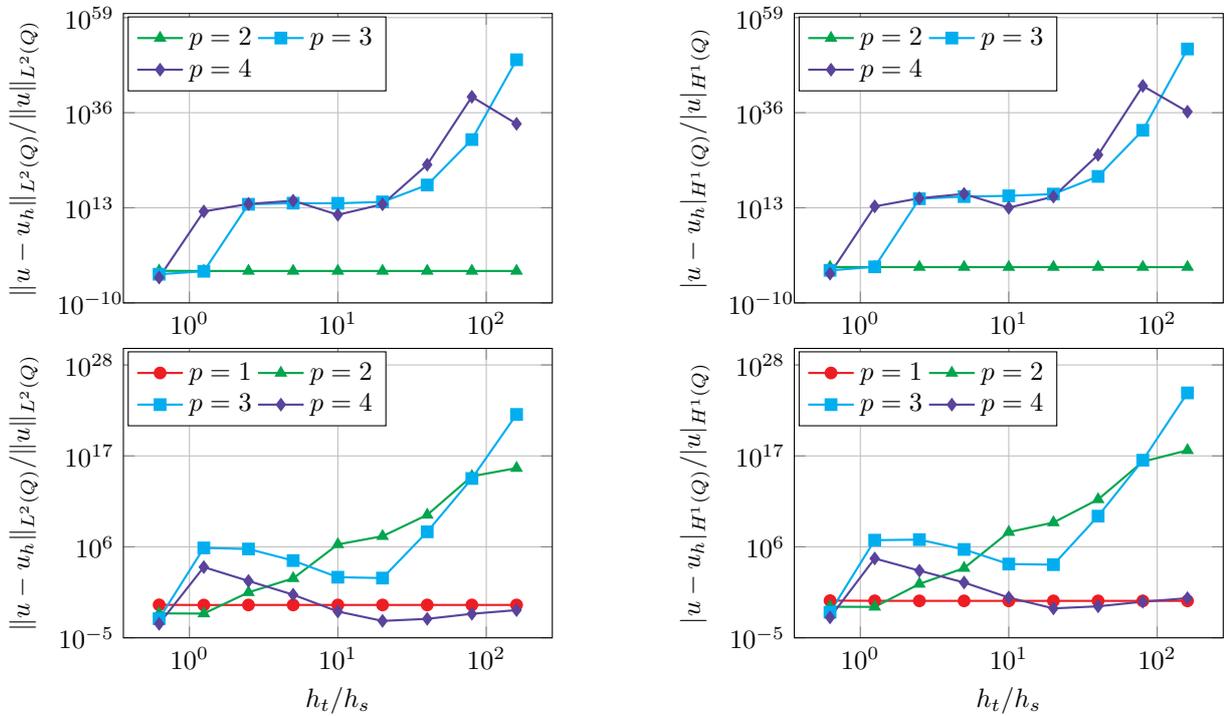
\begin{figure}[b]
\centering
\hspace*{\fill}
\begin{subfigure}[t]{0.425\linewidth}
\centering
\begin{tikzpicture}[%font=\footnotesize, 
trim axis left]
\begin{loglogaxis}[
cycle list name=Lista2,
yticklabel style={text width=2.125em, align=right},
width=\linewidth,
height=0.75\linewidth,
ymin=1e-10,
ymax=1e60,
xminorticks=false,
yminorticks=false,
%max space between ticks=20,
%xtick=data,
%xticklabel={\pgfmathparse{exp(\tick)}\pgfmathprintnumber[fixed relative,precision=3]{\pgfmathresult}},
%ytick={1e-9,1e6,1e21,1e36,1e51},
%yticklabels={$0.4$,$0.6$,$0.8$,$1$}, 
ylabel={$\| u - u_h \|_{L^2(Q)}/\| u \|_{L^2(Q)}$},
legend columns=2,
%legend pos=south east,
legend style={at={(0.01,0.99)},anchor= north west},
xmajorgrids=true,
ymajorgrids=true,
legend entries={$p=2$,$p=3$,$p=4$}]
%\addplot table [x=htfhs, y=L2P1R2, col sep=comma] {\DATAFILE};
\addplot table [x=htfhs, y=L2P2R2, col sep=comma] {\DATAFILE};
\addplot table [x=htfhs, y=L2P3R2, col sep=comma] {\DATAFILE};
\addplot table [x=htfhs, y=L2P4R2, col sep=comma] {\DATAFILE};
\end{loglogaxis}
\end{tikzpicture}
\end{subfigure}\hfill
\begin{subfigure}[t]{0.425\linewidth}
\centering
\begin{tikzpicture}[%font=\footnotesize,
trim axis right]
\begin{loglogaxis}[
yticklabel style={text width=2.125em, align=right},
cycle list name=Lista2,
width=\linewidth,
height=0.75\linewidth,
ymin=1e-10,
ymax=1e60,
xminorticks=false,
yminorticks=false,
%max space between ticks=20,
%xtick=data,
%xticklabel={\pgfmathparse{exp(\tick)}\pgfmathprintnumber[fixed relative,precision=3]{\pgfmathresult}},
%ytick={1e-9,1e6,1e21,1e36,1e51},
%yticklabels={$0.4$,$0.6$,$0.8$,$1$}, 
ylabel={$| u - u_h |_{H^1(Q)}/| u |_{H^1(Q)}$},
legend columns=2,
%legend pos=south east,
legend style={at={(0.01,0.99)},anchor= north west},
xmajorgrids=true,
ymajorgrids=true,
legend entries={$p=2$,$p=3$,$p=4$}]
%\addplot table [x=htfhs, y=H1P1R2, col sep=comma] {\DATAFILE};
\addplot table [x=htfhs, y=H1P2R2, col sep=comma] {\DATAFILE};
\addplot table [x=htfhs, y=H1P3R2, col sep=comma] {\DATAFILE};
\addplot table [x=htfhs, y=H1P4R2, col sep=comma] {\DATAFILE};
\end{loglogaxis}
\end{tikzpicture}
\end{subfigure}
\hspace*{\fill}
\\
\hspace*{\fill}
\begin{subfigure}[t]{0.425\linewidth}
\centering
\begin{tikzpicture}[%font=\footnotesize, 
trim axis left]
\begin{loglogaxis}[
yticklabel style={text width=2.125em, align=right},
cycle list name=Lista1,
width=\linewidth,
height=0.75\linewidth,
xlabel={$h_t/h_s$},
ymin=1e-5,
ymax=1e30,
xminorticks=false,
yminorticks=false,
%max space between ticks=20,
%xtick=data,
%xticklabel={\pgfmathparse{exp(\tick)}\pgfmathprintnumber[fixed relative,precision=3]{\pgfmathresult}},
%ytick={0.4,0.6,0.8,1}, 
%yticklabels={$0.4$,$0.6$,$0.8$,$1$}, 
%ytick={1e-5,1e2,1e9,1e16,1e23}, 
ylabel={$\| u - u_h \|_{L^2(Q)}/\| u \|_{L^2(Q)}$},
legend columns=2,
%legend pos=south east,
legend style={at={(0.01,0.99)},anchor= north west},
xmajorgrids=true,
ymajorgrids=true,
legend entries={$p=1$,$p=2$,$p=3$,$p=4$}]
\addplot table [x=htfhs, y=L2P1R, col sep=comma] {\DATAFILE};
\addplot table [x=htfhs, y=L2P2R, col sep=comma] {\DATAFILE};
\addplot table [x=htfhs, y=L2P3R, col sep=comma] {\DATAFILE};
\addplot table [x=htfhs, y=L2P4R, col sep=comma] {\DATAFILE};
\end{loglogaxis}
\end{tikzpicture}
\end{subfigure}\hfill
\begin{subfigure}[t]{0.425\linewidth}
\centering
\begin{tikzpicture}[%font=\footnotesize,
trim axis right]
\begin{loglogaxis}[
yticklabel style={text width=2.125em, align=right},
cycle list name=Lista1,
width=\linewidth,
height=0.75\linewidth,
xlabel={$h_t/h_s$},
ymin=1e-5,
ymax=1e30,
xminorticks=false,
yminorticks=false,
%max space between ticks=20,
%xtick=data,
%xticklabel={\pgfmathparse{exp(\tick)}\pgfmathprintnumber[fixed relative,precision=3]{\pgfmathresult}},
%ytick={1e-5,1e2,1e9,1e16,1e23}, 
%yticklabels={$0.4$,$0.6$,$0.8$,$1$}, 
ylabel={$| u - u_h |_{H^1(Q)}/| u |_{H^1(Q)}$},
legend columns=2,
%legend pos=south east,
legend style={at={(0.01,0.99)},anchor= north west},
xmajorgrids=true,
ymajorgrids=true,
legend entries={$p=1$,$p=2$,$p=3$,$p=4$}]
\addplot table [x=htfhs, y=H1P1R, col sep=comma] {\DATAFILE};
\addplot table [x=htfhs, y=H1P2R, col sep=comma] {\DATAFILE};
\addplot table [x=htfhs, y=H1P3R, col sep=comma] {\DATAFILE};
\addplot table [x=htfhs, y=H1P4R, col sep=comma] {\DATAFILE};
\end{loglogaxis}
\end{tikzpicture}
\end{subfigure}
\hspace*{\fill}
\caption{Relative errors of the \femstab\ method \eqref{stabZank} with $C^{p-2}$-continuous splines (top), splines of maximal regularity (bottom), plotted against the ratio $h_t/h_s$. The plots confirm what was observed in Figure~\ref{fig:StabZankPConvTot}.}
\label{fig:StabZankPRobustTot}
\end{figure}